\documentclass[12pt]{article}
\usepackage{amsmath,amsfonts,amscd,latexsym,amsxtra}

\begin{document}

\newcommand{\comment}[1]{{\bf #1}}

\newcommand{\tit}[1]{
\refstepcounter{section}
\vspace{5mm}
\begin{bf} \begin{center} \begin{Large}
\arabic{section}. #1
\end{Large}\end{center}\end{bf}}

\newcommand{\proof}{{\bf \underline{\underline{Proof}} : }}
\newcommand{\qed}{\hspace{10mm} \rule{2mm}{2mm}}

\newcommand{\bnull}{{\cal B}}

\newcommand{\cala}{{\cal A}}
\newcommand{\calb}{{\cal B}}
\newcommand{\calc}{{\cal C}}
\newcommand{\calf}{{\cal F}}
\newcommand{\calm}{{\cal M}}
\newcommand{\caln}{{\cal N}}
\newcommand{\calz}{{\cal Z}}
\newcommand{\calfin}{{\cal F\!I\!N}}
\newcommand{\rr}{{\mathbb R}}
\newcommand{\zz}{{\mathbb Z}}
\newcommand{\nn}{{\mathbb N}}
\newcommand{\zmod}{\zz/2}
\newcommand{\cc}{{\mathbb C}}
\newcommand{\qq}{{\mathbb Q}}

\newcommand{\map}{\operatorname{map}}
\newcommand{\colim}{\operatorname{colim}}
\newcommand{\tr}{\operatorname{tr}}
\newcommand{\pr}{\operatorname{pr}}
\newcommand{\id}{\operatorname{id}}
\newcommand{\im}{\operatorname{im}}
\newcommand{\cok}{\operatorname{cok}}
\newcommand{\clos}{\operatorname{clos}}
\newcommand{\defi}{\operatorname{def}}
\newcommand{\cone}{\operatorname{cone}}
\newcommand{\cyl}{\operatorname{cyl}}
\newcommand{\aut}{\operatorname{aut}}
\newcommand{\res}{\operatorname{res}}
\newcommand{\sing}{\operatorname{sing}}
\newcommand{\sign}{\operatorname{sign}}
\newcommand{\cell}{\operatorname{cell}}
\newcommand{\inte}{\operatorname{inte}}
\newcommand{\ind}{\operatorname{ind}}
\newcommand{\Tor}{\operatorname{Tor}}
\newcommand{\Ext}{\operatorname{Ext}}
\newcommand{\COMPLEXES}{\operatorname{COMPLEXES}}
\newcommand{\CHAINCOMPLEXES}{\operatorname{CHAIN-COMPLEXES}}
\newcommand{\Or}{\operatorname{Or}}
\newcommand{\supp}{\operatorname{supp}}
\newcommand{\virt}{\operatorname{virt}}
\newcommand{\class}{\operatorname{class}}
\newcommand{\con}{\operatorname{con}}
\newcommand{\consub}{\operatorname{consub}}
\newcommand{\ABEL}{\operatorname{ABEL}}
\newcommand{\vol}{\operatorname{vol}}
\newcommand{\ch}{\operatorname{ch}}
\newcommand{\HS}{\operatorname{HS}}
\newcommand{\cent}{\operatorname{cent}}
\newcommand{\limone}{\operatorname*{lim^1}}

\newcommand{\fgh}{\{\mbox{fin. gen. Hilb. } \cala \mbox{-mod.}\}}
\newcommand{\fgp}{\{\mbox{fin. gen. proj. } \cala \mbox{-mod.}\}}
\newcommand{\fp}{\{\mbox{fin. pres. } \cala \mbox{-mod.}\}}
\newcommand{\fgpf}{\{\mbox{fin. gen. proj. } \cala \mbox{-mod. with
inner prod.}\}}

\newcommand{\bfT}{\ensuremath{\mathbf{T}}}
\newcommand{\bfP}{\ensuremath{\mathbf{P}}}

\newcommand{\capac}{c}


\newcommand{\brn}{\begin{refnumber} \em }
\newcommand{\ern}{\em \end{refnumber}}

\newtheorem{theorem}{Theorem}[section]
\newtheorem{proposition}[theorem]{Proposition}
\newtheorem{lemma}[theorem]{Lemma}
\newtheorem{definition}[theorem]{Definition}
\newtheorem{example}[theorem]{Example}
\newtheorem{remark}[theorem]{Remark}
\newtheorem{notation}[theorem]{Notation}
\newtheorem{refnumber}[theorem]{}
\newtheorem{corollary}[theorem]{Corollary}
\newtheorem{assumption}[theorem]{Assumption}
\newtheorem{conjecture}[theorem]{Conjecture}
\newtheorem{question}[theorem]{Question}
{\catcode`@=11\global\let\c@equation=\c@theorem}
\renewcommand{\theequation}{\thetheorem}


\newcommand{\comsquare}[8]{
\begin{center}
$\begin{CD}
#1 @>#2>> #3\\
@V{#4}VV @VV{#5}V\\
#6 @>>#7> #8
\end{CD}$
\end{center}}

\newcommand{\textshortexactsequence}[5]{\mbox{$0\longrightarrow
#1\stackrel{#2}{\longrightarrow} #3\stackrel{#4}{\longrightarrow}
#5\longrightarrow 0$}}
\newcommand{\textshortexactsequenceone}[5]{\mbox{$1\longrightarrow
#1\stackrel{#2}{\longrightarrow} #3\stackrel{#4}{\longrightarrow}
#5\longrightarrow 1$}}



\typeout{--------------------- introduction -----------------------}
 \begin{center} \begin{LARGE} \begin{bf}
``Dimension theory of arbitrary modules
over finite von Neumann algebras and applications
to $L^2$-Betti numbers''\\[2mm]
by
\\[2mm]
Wolfgang L\"uck
\end{bf} \end{LARGE}\\[5mm]
\begin{Large}The paper is dedicated to Martin Kneser on the occasion of his
seventieth birthday
\end{Large}\end{center}

\begin{abstract}
We define for arbitrary modules over a finite von Neumann algebra $\cala$
a dimension taking values in $[0,\infty]$
which extends the classical notion of von Neumann dimension
for finitely generated projective $\cala$-modules and inherits
all its useful properties such as Additivity, Cofinality and Continuity.
This allows to define $L^2$-Betti numbers for arbitrary topological spaces
with an action of a discrete group $\Gamma$ extending the well-known
definition for regular coverings of compact manifolds. We show for an
amenable group $\Gamma$ that the $p$-th $L^2$-Betti number depends
only on the $\cc\Gamma$-module given by the
$p$-th singular homology. Using the generalized dimension
function we detect elements in $G_0(\cc\Gamma)$, provided that $\Gamma$
is amenable. We investigate the class of groups
for which the zero-th and first $L^2$-Betti numbers
resp. all $L^2$-Betti numbers
vanish. We study $L^2$-Euler characteristics and introduce
for a discrete group $\Gamma$ its Burnside group
extending the classical notions of Burnside ring and Burnside ring
congruences for finite $\Gamma$.
\\[3mm]
AMS-classification numbers: 55T99, 46L99\end{abstract}

\begin{bf} \begin{center} \begin{Large}
\vspace*{5mm}Introduction
\end{Large}\end{center}\end{bf}

Let us recall the original definition of
$L^2$-Betti numbers by Atiyah \cite{Atiyah (1976)}.
Let \mbox{$\overline{M} \longrightarrow M$} be a regular covering
of a closed Riemannian manifold $M$ with $\Gamma$ as group of deck
transformations. We lift the Riemannian metric to
a $\Gamma$-invariant Riemannian metric on $\overline{M}$.
Let \mbox{$L^2\Omega^p(\overline{M})$} be the Hilbert
space completion of the space
\mbox{$C^{\infty}_0\Omega^p(\overline{M})$}
of smooth $\rr$-valued $p$-forms on $\overline{M}$ with compact
support and the standard $L^2$-pre-Hilbert structure.
The Laplace operator $\Delta_p$ is
essentially selfadjoint in \mbox{$L^2\Omega^p(\overline{M})$}.
Let \mbox{$\Delta_p = \int \lambda dE_{\lambda}^p$} be the spectral
decomposition with right-continuous spectral family
\mbox{$\{E_{\lambda}^p \mid \lambda \in \rr \}$}. Let
\mbox{$E_{\lambda}^p(\overline{x},\overline{y})$} be the
Schwartz kernel of \mbox{$E_{\lambda}^p$}. Since
\mbox{$E_{\lambda}^p(\overline{x},\overline{x})$} is an endomorphism
of a finite-dimensional real vector space, its trace
\mbox{$\tr_{\rr}(E_{\lambda}^p(\overline{x},\overline{x})) \in\rr$}
is defined. Let $\calf$ be a fundamental domain for the $\Gamma$-action
on $\overline{M}$. Define the {\em analytic $L^2$-Betti number} by
\begin{eqnarray}
b_p^{(2)}(\overline{M}) & := & \int_{\calf}
\tr_{\rr}\left(E_{0}^p(\overline{x},\overline{x})\right)
d\!\vol_{\overline{x}}
\hspace{10mm} \in [0,\infty).
\label{analytic p- th L^2-Betti number}
\end{eqnarray}
By means of a Laplace transformation this can also be expressed
in terms of the heat kernel
\mbox{$e^{-t\Delta_p}(\overline{x},\overline{y})$}
on $\overline{M}$ by
\begin{eqnarray}
b_p^{(2)}(\overline{M}) & = & \lim_{t \to \infty} \int_{\calf}
\tr_{\rr}\left(e^{-t\Delta_p}(\overline{x},\overline{x})\right)
d\!\vol_{\overline{x}} \hspace{10mm} \in [0,\infty).
\label{L^2-Betti numbers and Laplace transform}
\end{eqnarray}
The $p$-th $L^2$-Betti number measures the size of the space
of smooth harmonic $L^2$-integrable $p$-forms on $\overline{M}$
and vanishes precisely if there is no such non-trivial form.
For a survey on $L^2$-Betti numbers and related invariants like
Novikov-Shubin invariants and $L^2$-torsion and their applications
and relations to geometry, spectral theory, group theory and $K$-theory
we refer for instance to
\cite[section 8]{Gromov (1993)},
\cite{Lott (1997)},
\cite{Lueck (1997)} and
\cite{Pansu (1996)}. In this paper, however, we will
not deal with the analytic side, but take an algebraic point of view.\par

The $L^2$-Betti numbers can also be defined in an algebraic manner.
Farber \cite{Farber (1996)},\cite{Farber(1996a)}
has shown that the category
of finitely generated Hilbert $\cala$-modules for a finite von Neumann
algebra $\cala$ can be embedded in an appropriate
abelian category and that one
can treat $L^2$-homology from a homological algebraic point of view.
Farber gives as an application
for instance an improvement of the Morse inequalites
of Novikov and Shubin \cite{Novikov-Shubin (1987)},
\cite{Novikov-Shubin (1987b)} in terms of $L^2$-Betti numbers
by taking the minimal number of generators into account. An equivalent
more algebra oriented approach is developed
in \cite{Lueck (1995)} where it is shown that the category
of finitely generated projective modules over $\cala$,
viewed just as a ring, is equivalent to the category
of finitely generated Hilbert $\cala$-modules
and that the category of finitely presented $\cala$-modules
is an abelian category. This allows to define for a
finitely generated projective
$\cala$-module $P$ {\em its von Neumann dimension}
\begin{eqnarray}
\label{eqn 0.0}
\dim_{\cala}(P) & \in & [0,\infty)
\end{eqnarray}
by using the classical definition
for finitely generated Hilbert $\cala$-modules
in terms of the von Neumann trace of a projector. This will be reviewed
in Section \ref{Review of von Neumann dimension}. \par

In Section
\ref{The generalized dimension function} we will prove the main
technical result of this paper that
this dimension can be extended to arbitrary $\cala$-modules
if one allows that the value may be infinite (what fortunately does
not happen in a lot of interesting situations). Moreover, this extension
inherits all good properties from the original definition
for finitely generated Hilbert $\cala$-modules such as
Additivity, Cofinality and Continuity and is
uniquely determined by these properties. More precisely, we will
introduce

\begin{definition} \label{new definition of dimension}
Define for a $\cala$-module $M$
$$\dim^{'}(M) ~ := ~
\sup\{\dim(P) \mid P \subset M \mbox{ finitely generated
projective } \cala\mbox{-submodule}\} \hspace{2mm} \in [0,\infty].$$
\qed
\end{definition}

Recall that the {\em dual module}
\mbox{$M^{\ast}$} of a left $\cala$-module is the left $\cala$-module
\mbox{$\hom_{\cala}(M,\cala)$} where the
$\cala$-multiplication is given by
\mbox{$(af)(x) = f(x)a^{\ast}$} for
$f \in M^{\ast}$, $x \in M$ and $a \in \cala$.

\begin{definition} \label{closure, K(M), P(M)}
Let $K$ be a $\cala$-submodule of the $\cala$-module $M$. Define
the {\em closure of $K$ in $M$}
to be the $\cala$-submodule of $M$
$$\overline{K} ~:=~ \{ x \in M ~ \mid ~ f(x) = 0 \mbox{ for all }
f \in M^{\ast}  \mbox{ with } K \subset \ker(f)\}.$$
For a finitely generated $\cala$-module $M$
define the $\cala$-submodule $\bfT M$ and
the $\cala$-quotient module $\bfP M$ by:
\begin{eqnarray*}
\bfT M & := &\{x \in M ~ \mid ~ f(x) = 0
\mbox{ for all } f \in M^{\ast}\};
\\
\bfP M  & := & M/\bfT M. \qed\end{eqnarray*}
\end{definition}

The notion of $\bfT M$ and $\bfP M$ corresponds in
\cite{Farber (1996)} to the torsion part and the projective part.
Notice that $\bfT M$ is the closure of the trivial submodule in $M$.
It is the same as the kernel of the canonical
map \mbox{$i(M) : M \longrightarrow (M^{\ast})^{\ast}$} which sends
$x \in M$ to the map
\mbox{$M^{\ast} \longrightarrow \cala \hspace{3mm}
f \mapsto f(x)^{\ast}$}. We will prove for a finite von Neumann algebra
$\cala$ in Section \ref{The generalized dimension function}

\begin{theorem} \label{properties of new definition}
\begin{enumerate}

\item $\cala$ is semi-hereditary, i.e. any finitely generated submodule
of a projective module is projective;

\item If $K \subset M$ is a submodule of the finitely generated
$\cala$-module $M$, then $M/\overline{K}$ is finitely generated and
projective and $\overline{K}$ is a direct summand in $M$;

\item If $M$ is a finitely generated $\cala$-module, then
$\bfP M$ is finitely generated projective and
$$M \cong \bfP M \oplus \bfT M;$$

\item The dimension $\dim^{'}$ has the following properties:

\begin{enumerate}
\item Continuity \par\noindent

If $K \subset M$ is a submodule of the finitely generated
$\cala$-module $M$, then:
$$\dim^{'}(K) ~ = ~ \dim^{'}(\overline{K});$$

\item Cofinality \par\noindent

Let $\{M_i\mid i \in I\}$ be a cofinal system of submodules of $M$,
i.e. \mbox{$M = \cup_{i \in I} M_i$} and for two indices $i$ and $j$
there is an index $k$ in $I$ satisfying \mbox{$M_i,M_j \subset M_k$}.
Then:
$$\dim^{'}(M) ~ = ~ \sup\{\dim^{'}(M_i) \mid i \in I\};$$

\item Additivity \par\noindent
If \mbox{$0 \longrightarrow M_0 \stackrel{i}{\longrightarrow} M_1
\stackrel{p}{\longrightarrow} M_2 \longrightarrow 0$} is an
exact sequence of $\cala$-modules, then:
$$\dim^{'}(M_1) ~ = ~\dim^{'}(M_0) + \dim^{'}(M_2),$$
where $r + s$ for \mbox{$r,s \in [0,\infty]$} is the ordinary sum of two
real numbers if both $r$ and $s$ are not $\infty$ and is $\infty$ otherwise;

\item Extension Property\par\noindent
If $M$ is finitely generated projective, then:
$$\dim{'}(M) ~ = ~ \dim(M);$$

\item If $M$ is a finitely generated $\cala$-module, then:
\begin{eqnarray*}
\dim^{'}(M) & = & \dim(\bfP M); \\
\dim^{'}(\bfT M) & = & 0;
\end{eqnarray*}

\item The dimension $\dim^{'}$ is uniquely determined by
Continuity, Cofinality, Additivity and the Extension Property. \qed

\end{enumerate}
\end{enumerate}
\end{theorem}

In the sequel we write $\dim$ instead of $\dim'$.
In Section \ref{Induction for group von Neumann algebras}
we will show for an inclusion \mbox{$i: \Delta \longrightarrow \Gamma$}
that the dimension
function is compatible with induction with the induced ring homomorphism
\mbox{$i: \caln(\Delta) \longrightarrow \caln(\Gamma)$}
and that $\caln(\Gamma)$ is faithfully flat over $\caln(\Delta)$.
(Theorem \ref{exactness and dimension under induction}).
This is important if one wants to relate the $L^2$-Betti numbers
of a regular covering to the ones of the universal covering.
We will prove that $\Gamma$ is non-amenable if and only if
\mbox{$\caln(\Gamma) \otimes_{\cc\Gamma} \cc$} is trivial
(Lemma \ref{invariants of caln(Gamma) otimes_{zzGamma} zz}.2). This
generalizes the result of Brooks \cite{Brooks (1981)}
(Remark \ref{alpha(X) = infty^+}).
\par

In Section \ref{L^2-invariants for arbitrary Gamma-spaces}
we use this generalized dimension function
to define for a (discrete) group $\Gamma$
and a $\Gamma$-space $X$ its
{\em $p$-th $L^2$-Betti number} by
\begin{eqnarray}
b_p^{(2)}(X;\caln(\Gamma)) & := &
\dim_{\caln(\Gamma)}\left(H_p^{\Gamma}(X;\caln(\Gamma))\right)
 \hspace{10mm} \in [0,\infty],
\label{L^2-Betti number of a Gamma-space}
\end{eqnarray}
where \mbox{$H_p^{\Gamma}(X;\caln(\Gamma))$} denotes the
$\caln(\Gamma)$-module given by the singular
homology of $X$ with coefficients in the
$\caln(\Gamma)$-$\zz\Gamma$-bimodule $\caln(\Gamma)$
(Definition \ref{L^2-Betti numbers for arbitrary spaces}).
This definition agrees with Atiyah's definition
\ref{analytic p- th L^2-Betti number} if $X$ is the total space
and $\Gamma$ the group of deck transformations of a regular covering
of a closed Riemannian manifold. We will compare our definition also
with the one of Cheeger and Gromov
\cite[section 2]{Cheeger-Gromov (1986)} (Remark
\ref{comparision with the definition of Cheeger and Gromov}).
In particular we can define for an
arbitrary (discrete) group $\Gamma$ its $p$-th $L^2$-Betti number
\begin{eqnarray}
b_p^{(2)}(\Gamma) & := & b_p^{(2)}(E\Gamma;\caln(\Gamma))
\hspace{10mm} \in [0,\infty],
\label{L^2-Betti numbers of a group}
\end{eqnarray}
where \mbox{$E\Gamma \longrightarrow B\Gamma$}
is the universal $\Gamma$-principal bundle. These generalizations
inherit all the useful properties from the original versions
and it pays off to have them at hand in this generality.
For instance if one is only interested in the $L^2$-Betti numbers
of a group $\Gamma$ for which $B\Gamma$ is
a $CW$-complex of finite type and hence the original (simplicial)
definition does apply,
it is important to have the more general definition available
because such a group $\Gamma$ may contain an interesting normal
subgroup $\Delta$ which is not even finitely generated.
A typical situation is when $\Gamma$ contains a normal infinite amenable
subgroup $\Delta$. Then all the $L^2$-Betti numbers of $B\Gamma$
are trivial by a result of Cheeger and Gromov
\cite[Theorem 0.2 on page 191]{Cheeger-Gromov (1986)}.
This result was the main motivation for our
attempt to construct the extensions of  dimension and of
$L^2$-Betti numbers described above.\par

In Section  \ref{Amenable groups}
we will get the theorem of Cheeger and Gromov mentioned above
as a corollary of the following result. If $\Gamma$ is amenable and
$M$ is a $\cc\Gamma$-module, then
\begin{eqnarray}
\dim_{\caln(\Gamma)}
\left(\Tor^{\cc\Gamma}_p(M,\caln( \Gamma))\right) & = & 0
\label{eqn 0.1}
\hspace{10mm} \mbox{ for }
p \ge 1,
\end{eqnarray}
where we consider $\caln(\Gamma)$ as a
$\caln(\Gamma)$-$\cc\Gamma$-bimodule
(Theorem \ref{dimension of higher Tor-s vanish in the amenable case}).
We get from \ref{eqn 0.1}
by a spectral sequence argument that the $L^2$-Betti numbers of a
$\Gamma$-space $X$ depend only on its singular homology
with complex coefficients viewed as $\cc\Gamma$-module, namely
\begin{eqnarray}
\label{eqn 0.2}
b_p^{(2)}(X;\caln(\Gamma)) & = &
\dim_{\caln(\Gamma)}\left(\caln(\Gamma) \otimes_{\cc\Gamma}
H_p^{\sing}(X;\cc)\right),
\end{eqnarray}
provided that $\Gamma$ is amenable
(Theorem \ref{L^2-Betti numbers and homology in the amenable case}).
The result of Cheeger and Gromov mentioned above follows from
\ref{eqn 0.2}
since the singular homology of $E\Gamma$ is
trivial in all dimension except
for dimension 0 where it is $\cc$.

Equation \ref{eqn 0.1}
will also play a cruircial role in detecting
non-trivial elements in the Grothendieck group $G_0(\cc\Gamma)$
of finitely generated $\cc\Gamma$-modules for amenable groups
$\Gamma$ which will be investigated in Section
\ref{Dimension functions and G_0(ccGamma)}. We will construct
for amenable $\Gamma$ a map
\begin{eqnarray}
G_0(\cc\Gamma)\otimes_{\zz}\cc & \longrightarrow & \class(\Gamma)_{cf},
\label{detecting elements in G_0(CGamma)}
\end{eqnarray}
where \mbox{$\class(\Gamma)_{cf}$} is the complex
vector space of functions from the set
$\con(\Gamma)_{cf}$ of finite conjugacy classes $(\gamma)$
of elements in $\Gamma$ to $\cc$
(Lemma \ref{map from G_0(CGamma) to K_0(N(Gamma)}
and Theorem \ref{K_0 and G_0 and HS and dim^u}).
This map is related to the Hattori-Stallings rank and the universal
center-valued trace and dimension
of $\caln(\Gamma)$ (Theorem \ref{K_0 and G_0 and HS and dim^u}).
In particular we will show that the class of
$\cc\Gamma$ in $G_0(\cc\Gamma)$ generates an infinite cyclic
direct summand in $G_0(\cc\Gamma)$ if $\Gamma$
is amenable and is trivial if $\Gamma$ contains a free group of rank 2
as subgroup
(Remark \ref{remark on amenbale groups and [CGamma] in G_0}).\par

We will investigate for \mbox{$d = 0,1 , \ldots$}
and \mbox{$d = \infty$} the class $\bnull_d$ of groups $\Gamma$
for which \mbox{$b_p(E\Gamma;\caln(\Gamma)) = 0$} for \mbox{$ p \le d$}
(Theorem \ref{properties of bnull})
and discuss applications in Section
\ref{Groups with vanishing L^2-Betti numbers}
(Theorem \ref{consequence from Gamma in bnull}).\par

We analyse $L^2$-Euler characteristics and the Burnside group
in Section \ref{L^2-Euler characteristics and the Burnside group}
generalizing the classical notions of Burnside ring,
Burnside ring congruences
and equivariant Euler characteristic for finite groups to infinite groups
(Theorem \ref{chi^{(2)} and cells},
Lemma \ref{chi_K^{Gamma}(X) for finite proper G-CW-complex X},
Lemma \ref{injectivity resp bijectivity of the global
character map},
Remark \ref{explicite inverse of the global L^2-character map}
and Lemma \ref{amenable groups and chi^{Gamma}_K(E(Gamma,calfin)))}).
In particular the $L^2$-Euler characteristic extends
the notion of virtual Euler characteristic of a group
to a larger class of groups (Corollary \ref{vanishing of chi_{virt}}).
\par

In Section \ref{Values of L^2-Betti numbers} we analyse the possible values
of the $L^2$-Betti numbers
(Theorem \ref{possible values of the L^2-Betti numbers}).
If there is no bound on the orders of finite subgroups of $\Gamma$,
then any non-negative real number can be realized as
$b_p^{(2)}(X;\caln(\Gamma))$ for $p \ge 3$
and a free $\Gamma$-$CW$-complex $X$. Otherwise
we show for the least
common multiple $d$ of the orders of finite subgroups that
\mbox{$d \cdot b_p^{(2)}(X;\caln(\Gamma))$} is an integer or infinite
for any $\Gamma$-space $X$ if this holds for any finite
free $\Gamma$-$CW$-complex $Y$
(Theorem \ref{possible values of the L^2-Betti numbers}).
The last conditions hold for instance for elementary
amenable groups and free groups by Linnell \cite{Linnell (1993)}.

The paper is organized as follows:\\[5mm]
\begin{tabular}{ll}
0. & Introduction
\\
\ref{Review of von Neumann dimension}.
& Review of von Neumann dimension
\\
\ref{The generalized dimension function}.
& The generalized dimension function
\\
\ref{Induction for group von Neumann algebras}.
& Induction for group von Neumann algebras
\\
\ref{L^2-invariants for arbitrary Gamma-spaces}.
& $L^2$-invariants for arbitrary $\Gamma$-spaces
\\
\ref{Amenable groups}. & Amenable groups
\\
\ref{Dimension functions and G_0(ccGamma)}. &
Dimension functions and $G_0(\cc\Gamma)$
\\
\ref{Groups with vanishing L^2-Betti numbers}. &
Groups with vanishing $L^2$-Betti numbers
\\
\ref{L^2-Euler characteristics and the Burnside group}. &
$L^2$-Euler characteristics and the Burnside group
\\
\ref{Values of L^2-Betti numbers}. &
Values of $L^2$-Betti numbers
\\ & References
\end{tabular}


\setcounter{section}{0}

\typeout{--------------------- section 1 -----------------------}

\tit{Review of von Neumann dimension}
\label{Review of von Neumann dimension}

In this section we recall some basic facts
about finitely generated Hilbert-modules
and finitely generated projective
modules over a finite von Neumann algebra.
We fix for the sequel
\begin{notation}
\label{basic notation}
Let $\cala$ be a finite von Neumann algebra and
\mbox{$\tr: \cala \longrightarrow \cc $} be a normal finite faithful
trace. Denote by $\Gamma$ an (arbitrary) discrete group.
Let $\caln(\Gamma)$ be the group von Neumann algebra
with the standard trace $\tr_{\caln(\Gamma)}$.
\par

Module means always left-module and group actions on spaces are from
the left unless explicitly stated differently.
We will always work in the category of
compactly generated spaces (see \cite{Steenrod (1967)} and
\cite[I.4]{Whitehead (1978)}). \qed\end{notation}

Next we recall our main example for $\cala$ and $\tr$, namely
{\em the group von Neumann algebra}
$\caln(\Gamma)$ with {\em the standard trace}. The reader
who is not familiar with the general concept of finite von
Neumann algebras may always think of this example.
Let $l^2(\Gamma)$ be the Hilbert space of formal sums
\mbox{$\sum_{\gamma \in \Gamma} \lambda_{\gamma} \cdot \gamma$}
with complex coefficients  $\lambda_{\gamma}$
which are square-summable, i.e.
\mbox{$\sum_{\gamma \in \Gamma}
|\lambda_{\gamma}|^2 < \infty $}. Define {\em the group
von Neumann algebra} and {\em the standard trace} by
\begin{eqnarray}
\caln(\Gamma) & := & \calb(l^2(\Gamma),l^2(\Gamma))^{\Gamma};
\label{definition of group von Neumann algebra}
\\
\tr_{\caln(\Gamma)}(a)  & := & \langle a(e),e\rangle_{l^2(\Gamma)};
\label{definition of standard trace}
\end{eqnarray}
where \mbox{$\calb(l^2(\Gamma),l^2(\Gamma))^{\Gamma}$}
is the space of
bounded $\Gamma$-equivariant operators
from $l^2(\Gamma)$ to itself,
\mbox{$a \in \caln (\Gamma)$} and
\mbox{$e \in \Gamma \subset l^2(\Gamma)$} is
the unit element. The given trace on $\cala$ extends to a trace on
square-matrices over $\cala$
in the usual way
\begin{eqnarray}
\tr: M(n,n,\cala) \longrightarrow \cc & \hspace{10mm} &
A ~ \mapsto ~ \sum_{i = 1}^n \tr(A_{i,i}).
\label{trace for matrices}
\end{eqnarray}
Taking adjoints induces the structure of a
{\em ring with involution} on $\cala$, i.e. we obtain a map
\mbox{$\ast : \cala \longrightarrow \cala \hspace{3mm}
a \mapsto a^{\ast}$},
which satisfies
\mbox{$(a+b)^{\ast} = a^{\ast} + b^{\ast}$},
\mbox{$(ab)^{\ast} = b^{\ast}a^{\ast}$} and
\mbox{$(a^{\ast})^{\ast} = a$} and \mbox{$1^{\ast} = 1$}
for all \mbox{$a,b \in \cala$}.
This involution induces an involution on matrices
\begin{eqnarray}
\ast : M(m,n,\cala) \longrightarrow M(n,m,\cala) & \hspace{3mm} &
A = (A_{i,j}) \mapsto A^{\ast} = (A_{j,i}^{\ast}).
\end{eqnarray}

\begin{definition}
\label{von Neumann dimension of a finitely generated projective
cala-module}
Let $P$ be a finitely generated projective $\cala$-module.
Let \mbox{$A \in M(n,n,\cala)$} be a matrix such that
\mbox{$A = A^{\ast}$}, \mbox{$A^2 = A$} and the image of the
$\cala$-linear map \mbox{$A: \cala^n \longrightarrow \cala^n$}
induced by right multiplication with $A$ is $\cala$-isomorphic to $P$.
Define the {\em von Neumann dimension} of $P$
$$\dim(P) = \dim_{\cala}(P) :=
\tr_{\cala}(A) \hspace{10mm} \in [0,\infty).
\qed$$
\end{definition}

It is not hard to check that this definition is independent
of the choice of $A$ and depends only on the isomorphism class of
$P$. Moreover the dimension is faithful, i.e. \mbox{$\dim(P) = 0$}
implies \mbox{$P = 0$}, is additive under direct sums and satisfies
\mbox{$\dim(\cala^n) = n$}.\par

We recall that we have defined
$\overline{K}$, $\bfT M$ and $\bfP M$ for \mbox{$K \subset M$}
in Definition \ref{closure, K(M), P(M)}. A sequence
\mbox{$L \stackrel{f}{\longrightarrow} M
\stackrel{g}{\longrightarrow} N$}
of $\cala$-modules is {\em weakly exact}
resp. {\em exact at $M$} if
\mbox{$\overline{\im(f)} = \ker(g)$} resp.
\mbox{$\im(f) =\ker(g)$} holds. A morphism
\mbox{$f: M \longrightarrow N$} of $\cala$-modules is called
a {\em weak isomorphism} if its kernel is trivial
and the closure of its image is $N$.\par

Next we explain how these concepts above correspond to their
analogues for finitely generated Hilbert $\cala$-modules.
Let $l^2(\cala)$ be the
Hilbert space completion of $\cala$ which is viewed as a
pre-Hilbert space by the inner product
\mbox{$\langle a , b \rangle = \tr(ab^{\ast})$}. A
{\em finitely generated Hilbert $\cala$-module} $V$
is a Hilbert space $V$ together with a left
operation of $\cala$ by $\cc$-linear maps such that
there exists a unitary $\cala$-embedding of $V$ in
\mbox{$\oplus_{i=1}^n l^2(\cala)$} for some $n$.
A morphism of finitely generated
Hilbert $\cala$-modules is a bounded  $\cala$-equivariant operator.
Denote by $\fgh$ the category of finitely generated
Hilbert $\cala$-modules. A sequence
\mbox{$U \xrightarrow{f} V \xrightarrow{g} W$} of
finitely generated Hilbert $\cala$-modules is
{\em exact resp. weakly exact at $V$} if
\mbox{$\im(f) = \ker(g)$} resp. \mbox{$\overline{\im(f)} = \ker(g)$}
holds. A morphism \mbox{$f: V \longrightarrow W$} is a
{\em weak isomorphism} if its kernel is trivial and its image is dense.
For a survey on finite von Neumann
algebras and Hilbert $\cala$-modules
we refer for instance
to \cite[section 1]{Lott-Lueck (1995)},
\cite[section 1]{Lueck-Rothenberg (1991)}. \par

The right regular representation
\mbox{$\cala  \longrightarrow  \calb(l^2(\cala),l^2(\cala))^{\cala}$}
from $\cala$ into the space of
bounded $\cala$-equivariant operators from $l^2(\cala)$ to itself
sends \mbox{$a \in \cala$} to the extension of the map
\mbox{$\cala \longrightarrow \cala \hspace{3mm}
b \mapsto ba^{\ast}$}. It is known
to be bijective
\cite[Theorem 1 in I.5.2 on page 80, Theorem 2 in I.6.2 on page 99]
{Dixmier (1981)}. Hence we obtain a bijection
\begin{eqnarray}
\nu: M(m,n,\cala) & \longrightarrow &
\calb(l^2(\cala)^m,l^2(\cala)^n)^{\cala},
\label{nu for cala^n}
\end{eqnarray}
which is compatible with the $\cc$-vector space structures,
the involutions and composition.\par

The details of the following theorem
and its proof can be found in \cite[section 2]{Lueck (1995)}.
It is essentially a consequence of \ref{nu for cala^n}
and the construction of the idempotent completion
of a category.
It allows us to forget the Hilbert-module-structures
and simply work with the von Neumann algebra as a plain ring.
An equivalent
approach is given by Farber  \cite{Farber (1996)},\cite{Farber(1996a)}
and is identified with the one here in
\cite[Theorem 0.9]{Lueck (1995)}.
An {\em inner product} on a finitely generated projective
$\cala$-module $P$ is a map
\mbox{$\mu: P \times P \longrightarrow \cala$}
which is linear in the first variable, symmetric in the sense
\mbox{$\mu(x,y) = \mu(y,x)^{\ast}$} and positive in the sense
\mbox{$\mu(x,x) > 0 \Longleftrightarrow x \not= 0$}  such that
the induced map \mbox{$P \longrightarrow P^{\ast}$} sending
\mbox{$y \in P$} to \mbox{$\mu(-,y)$} is bijective.

\begin{theorem} \label{functor nu}
\begin{enumerate}
\item
There is a functor
$$\nu : \fgpf \longrightarrow \fgh$$
which is  an equivalence of $\cc$-categories with involutions;

\item Any finitely generated projective
$\cala$-module has an inner product.
Two finitely generated projective
$\cala$-modules with inner product
are unitarily $\cala$-isomorphic if and only if
the underlying $\cala$-modules are $\cala$-isomorphic;

\item
Let $\nu^{-1}$ be an inverse of $\nu$ which is
well-defined up to unitary natural equivalence.
The composition of $\nu^{-1}$ with the forgetful functor induces an
equivalence of $\cc$-categories
$$\fgh \longrightarrow \fgp;$$

\item $\nu$ and $\nu^{-1}$ preserve
weak exactness and exactness. \qed
\end{enumerate}
\end{theorem}

Of course Definition
\ref{von Neumann dimension of a finitely generated projective
cala-module} of $\dim(P)$ for a finitely generated projective
$\cala$-module agrees with the usual von Neumann dimension of the
associated Hilbert $\cala$-module $\nu(P)$
after any choice of inner product on $P$.


\typeout{--------------------- section 2 -----------------------}

\tit{The generalized dimension function}
\label{The generalized dimension function}

In this section we give the proof of
Theorem \ref{properties of new definition} and investigate the
behaviour of the dimension under colimits.
We recall that we have introduced
\mbox{$\dim^{'}(M)$} for an arbitrary $\cala$-module $M$
in Definition \ref{new definition of dimension} and
$\overline{K}$, $\bfT M$ and $\bfP M$ for \mbox{$K \subset M$}
in Definition \ref{closure, K(M), P(M)}. We begin with the proof
of Theorem \ref{properties of new definition}.\par

\proof 1.) is proven in \cite[Corollary 2.4]{Lueck (1995)} for finite von 
Neumann algebras. However, Pardo pointed out to us that 
any von Neumann algebra is semi-hereditary. This follows from
the facts that any von Neumann algebra is a Baer $*$-ring and hence in 
particular a Rickart $C^*$-algebra 
\cite[Definition 1, Definition 2 and Proposition 9 in 
Chapter 1.4]{Berberian (1972)} and that a $C^*$-algebra is 
semi-hereditary if and only if it is
Rickart \cite[Corollary 3.7 on page 270]{Ara and Goldstein (1993)}.
\par\noindent
2.) and 4.a) in the special case that $M=P$
for a finitely generated projective $\cala$-module $P$.\\[2mm]
Let \mbox{${\cal P} = \{P_i \mid
i \in I\}$} be the directed system of finitely generated projective
$\cala$-modules of $K$. Notice that
\mbox{${\cal P}$} is indeed directed by inclusion since the
submodule of $P$ generated by two finitely generated projective
submodules is again finitely generated and hence by 1.)
finitely generated projective.
Let \mbox{$j_i: P_i \longrightarrow P$}
be the inclusion. Equip $P$ and each $P_i$ with a fixed inner product
and let
\mbox{$\pr_i: \nu (P) \longrightarrow \nu (P)$} be the orthogonal
projection satisfying \mbox{$\im (\pr_i) = \overline{\im (\nu (j_i))}$}
and \mbox{$\pr : \nu (P) \longrightarrow \nu (P)$} be the orthogonal
projection satisfying
\mbox{$\im (\pr) =  \overline{\cup_{i \in I} \im (\pr_i)}$}.
Next we show
\begin{eqnarray}
\im (\nu^{-1}(\pr)) & = & \overline{K}. \label{eqn 2.0}
\end{eqnarray}
Let \mbox{$f: P \longrightarrow \cala$} be a $\cala$-map with
\mbox{$K \subset \ker (f)$}. Then
\mbox{$f \circ j_i = 0$} and therefore
\mbox{$\nu (f) \circ \nu (j_i) = 0$} for all \mbox{$i \in I$}.
 We get
\mbox{$\im (\pr_i) \subset \ker(\nu (f))$} for all
\mbox{$i \in I$}. Because the kernel of $\nu(f)$ is closed
we conclude
\mbox{$\im (\pr) \subset \ker (\nu (f))$}.
This shows \mbox{$\im (\nu^{-1}(\pr)) \subset \ker (f)$}
and hence \mbox{$\im (\nu^{-1}(\pr)) \subset \overline{K}$}.
As \mbox{$K \subset \ker (\id - \nu^{-1}(\pr)) = \im (\nu^{-1}(\pr))$},
we conclude  \mbox{$\overline{K} \subset \im (\nu^{-1}(\pr))$}.
This finishes the proof of \ref{eqn 2.0}
and of 2.) in the special case $M=P$.
\par
Next we prove
\begin{eqnarray}
\dim'(K) & = & \dim(\overline{K}).
\label{eqn 2.1}
\end{eqnarray}
The inclusion $j_i$ induces a weak isomorphism
\mbox{$\nu (P_i) \longrightarrow \im(\pr_i)$}
of finitely generated Hilbert $\cala$-modules.
If we apply the Polar Decompostion Theorem to it
we obtain a unitary $\cala$-isomorphism
from \mbox{$\nu (P_i)$} to \mbox{$\im(\pr_i)$}.
This implies \mbox{$\dim(P_i) = \tr(\pr_i)$}.
Therefore it remains to prove
\begin{eqnarray}
\tr (\pr)  & := & \sup \{\tr(\pr_i) \mid i \in I\}.
\label{eqn 2.-1}
\end{eqnarray}
As $\tr$ is normal, it suffices to show for
\mbox{$x \in \nu (P)$} that the net \mbox{$\{\pr_i(x) \mid i \in I\}$}
converges to $\pr(x)$. Let $\epsilon > 0$ be given. Choose
$i(\epsilon) \in I$ and \mbox{$x_{i(\epsilon)} \in  \im(\pr_{i(\epsilon)})$}
with
\mbox{$||\pr(x) - x_{i(\epsilon )}|| \le ~\epsilon/2$}.
We conclude for all $i \ge i(\epsilon)$
\begin{eqnarray*}
||\pr(x) - \pr_i(x) || & \le &
||\pr(x) - \pr_{i(\epsilon)}(x)||
\\
& \le & ||\pr(x) - \pr_{i(\epsilon)}(x_{i(\epsilon)})|| +
||\pr_{i(\epsilon)}(x_{i(\epsilon)}) - \pr_{i(\epsilon)}(x)||
\\
& \le & ||\pr(x) - x_{i(\epsilon)}|| +
||\pr_{i(\epsilon)}(x_{i(\epsilon)} - \pr(x))||
\\
& \le & ||\pr(x) - x_{i(\epsilon)}|| +
||\pr_{i(\epsilon)}|| \cdot ||x_{i(\epsilon)} - \pr(x)||
\\
& \le & 2 \cdot ||\pr(x) - x_{i(\epsilon)}||
\\
& \le & \epsilon.
\end{eqnarray*}
Now \ref{eqn 2.-1} and hence \ref{eqn 2.1} follow.
In particular we get from \ref{eqn 2.1}
for any finitely generated projective submodule
$Q_0$ of a finitely generated projective $\cala$-module $Q$
\begin{eqnarray}
\dim(Q_0) & \le & \dim(Q), \label{eqn 2.-2}
\end{eqnarray}
since by definition \mbox{$\dim(Q_0) \le \dim'(Q_0)$}
and \mbox{$\dim(\overline{Q_0}) \le \dim(Q)$} follows
from additivity of $\dim$ under direct sums and that we have already
proven that $\overline{Q_0}$ is a direct summand in $Q$. This implies
for a finitely generated projective $\cala$-module $Q$
\begin{eqnarray}
\dim(Q) & = & \dim'(Q).
\label{eqn 2.3}
\end{eqnarray}
Now \ref{eqn 2.1} and \ref{eqn 2.3} imply
4.a.) in the special case
that $M = P$ for a finitely generated projective $\cala$-module
$P$.\par\noindent
4.d.) has been already proven in \ref{eqn 2.3}.
\par\noindent
4.b.) If $P \subset M$ is a finitely generated projective submodule,
then there is an index \mbox{$i \in I$} with \mbox{$P \subset M_i$} by
cofinality.
\par\noindent
4.c.) Let $P \subset M_2$ be a finitely generated projective submodule.
We obtain an exact sequence
\mbox{$0 \longrightarrow M_0 \longrightarrow p^{-1}(P) \longrightarrow
P \longrightarrow 0$}. Since \mbox{$p^{-1}(P) \cong M_0 \oplus P$},
we conclude
$$\dim^{'}(M_0) + \dim(P) \le \dim^{'}(p^{-1}(P)) \le \dim^{'}(M_1).$$
Since this holds for all finitely generated projective submodules
$P \subset M_2$, we get
\begin{eqnarray}
\label{eqn 2.4}
\dim^{'}(M_0) + \dim^{'}(M_2) & \le & \dim^{'}(M_1).
\end{eqnarray}
Let $Q \subset M_1$ be finitely generated projective.
We obtain exact  sequences
$$\begin{array}{ccccccccc}
0 & \longrightarrow & i(M_0) \cap Q & \longrightarrow & Q &
\longrightarrow & p(Q) & \longrightarrow & 0;
\\
0 & \longrightarrow & \overline{i(M_0) \cap Q} & \longrightarrow & Q &
\longrightarrow & Q/\overline{i(M_0) \cap Q} & \longrightarrow & 0.
\end{array}$$
By the special case of 2.) which we have already proven above
$\overline{i(M_0) \cap Q}$ is a direct summand in $Q$. We conclude
\begin{eqnarray*}
\dim(Q) & = & \dim(\overline{i(M_0) \cap Q}) +
\dim(Q/\overline{i(M_0) \cap Q}).
\end{eqnarray*}
From the special case 4.a.) we have already proven above,
4.d.) and the fact
that there is an epimorphism from $p(Q)$ onto the finitely generated
projective $\cala$-module $Q/ \overline{i(M_0) \cap Q}$, we conclude
\begin{eqnarray*}
\dim(\overline{i(M_0) \cap Q}) & = & \dim^{'}(i(M_0) \cap Q);
\\
\dim(Q/\overline{i(M_0) \cap Q}) & \le & \dim^{'}(p(Q)).
\end{eqnarray*}
Since obviously \mbox{$\dim^{'}(M) \le  \dim^{'}(N)$} holds for
$\cala$-modules $M$ and $N$ with $M \subset N$, we get
\begin{eqnarray*}
\dim(Q) & = & \dim(\overline{i(M_0) \cap Q}) +
\dim(Q/\overline{i(M_0) \cap Q})
\\
& \le &  \dim^{'}(i(M_0) \cap Q) + \dim^{'}(p(Q))
\\
& \ \le & \dim^{'}(M_0) + \dim^{'}(M_2).
\end{eqnarray*}
Since this holds for all finitely generated projective submodules
$Q \subset M_1$, we get
\begin{eqnarray}
\label{eqn 2.5}
\dim^{'}(M_1) & \le &
\dim^{'}(M_0) + \dim^{'}(M_2).
\end{eqnarray}
Now 4.c.) follows from \ref{eqn 2.4} and \ref{eqn 2.5}.
\par\noindent
2.) and 4.a.) Choose a finitely generated free $\cala$-module
$F$ together with an epimorphism
\mbox{$q: F \longrightarrow M$}.
One easily checks that \mbox{$q^{-1}(\overline{K})$} is
\mbox{$\overline{q^{-1}(K)}$} and that
\mbox{$F/q^{-1}(\overline{K})$} and
\mbox{$M/\overline{K}$} are isomorphic.
From the special case of 2.) and 4.) a.)
which we have already proven above
we conclude that
\mbox{$F/\overline{q^{-1}(K)}$} and hence
\mbox{$M/\overline{K}$} are finitely generated projective and
$$\dim^{'}(q^{-1}(K)) ~ = ~
\dim^{'}(\overline{q^{-1}(K)}) ~ = ~ \dim^{'}(q^{-1}(\overline{K})).$$
If $L$ is the kernel of $q$, we conclude from Additivity
\begin{eqnarray*}
\dim^{'}(q^{-1}(\overline{K}))& = &
\dim^{'}(L) + \dim^{'}(\overline{K});
\\
\dim^{'}(q^{-1}(K))& = &
\dim^{'}(L) + \dim^{'}(K).
\end{eqnarray*}
Now 2.) and 4.a) follow in general.
\par\noindent
3.) follows from 2.), as \mbox{$\overline{\{0\}} = \bfT M$}
and \mbox{$M/\bfT M = \bfP M$} by definition.
\par\noindent
4.e.) From 2.),  4.c.) and 4.d.) we get:
\mbox{$\dim^{'}(M) = \dim^{'}(\bfT M) + \dim(\bfP M)$}.
If we apply 4.a) to
\mbox{$\{0\} \subset M$} we get \mbox{$\dim^{'}(\bfT M)  =  0$}
because of
\mbox{$\overline{\{0\}} = \bfT M$}. \par\noindent
4.f.) Let $\dim^{''}$ be another function satisfying
Continuity, Cofinality, Additivity and the Extension Property.
We want to show for a
$\cala$-module $M$
$$\dim^{''}(M) = \dim^{'}(M).$$
Since 4.e.) is a consequence of Continuity, Additivity and
the Extension Property alone,
this is obvious provided $M$ is finitely generated. Since the system
of finitely generated submodules of a module is cofinal, the
claim follows from Cofinality.
This finishes the proof of Theorem
\ref{properties of new definition}. \qed
\par

\begin{notation}
\label{dim 0 dim^{prime}}
In view of Theorem \ref{properties of new definition} we
will not distinguish between $\dim^{'}$ and $\dim$
in the sequel. \qed
\end{notation}

Next we investigate the behaviour
of  dimension under colimits
indexed by a directed set.  We mention that
colimit is sometimes called in the literatur also
inductive limit or direct limit. The harder case of inverse limits which
is not needed in this paper will be treated at a different place
(see also \cite[Appendix]{Cheeger-Gromov (1986)}).

\begin{theorem}
\label{dimension of colimits for arbitrary index sets}
Let $I$ be a category such that between two objects
there is at most one morphisms and for two objects $i_1$
and $i_2$ there is an object $i_0$ with \mbox{$i_1 \le i_0$}
and \mbox{$i_2 \le i_0$} where we
write \mbox{$i \le k$} for two objects
$i$ and $k$ if and only if there is a morphism from $i$ to $k$.
Let $M_i$ be a covariant functor from $I$
to the category of $\cala$-modules.
For \mbox{$i \le j$} let \mbox{$\phi_{i,j}: M_i \longrightarrow M_j$}
be the associated morphism of $\cala$-modules.
For \mbox{$i \in I$} let \mbox{$\psi_{i}: M_i \longrightarrow \colim_I M_i$}
be the canonical morphism of $\cala$-modules. Then:
\begin{enumerate}

\item We get for the dimension of the $\cala$-module given by
the colimit $\colim_{I} M_i$
$$\dim\left(\colim_{I} M_i\right) ~ = ~
\sup\left\{\dim (\im(\psi_i)) \mid i \in I\right\};$$

\item Suppose
for each \mbox{$i \in I$} that there is \mbox{$i_0 \in I$} with
\mbox{$i \le i_0$} such that \mbox{$\dim(\im(\phi_{i,i_0})) < \infty$}
holds. Then:
$$\dim\left(\colim_{I} M_i\right) ~ = ~
\sup\left\{\inf\left\{
\dim(\im(\phi_{i,j}: M_i \longrightarrow M_j)) \mid j \in I, i \le j
\right\}\mid i \in I\right\}.$$
\end{enumerate}
\end{theorem}
\proof
1.) Recall that $\colim_{I} M_i$ can be constructed
as \mbox{$\coprod_{i \in I} M_i/\sim $} for the equivalence relation
for which \mbox{$x \in M_i \sim y \in M_j$} holds precisely
if there is \mbox{$k \in I $} with \mbox{$i \le k$} and
\mbox{$j \le k$} with the property
\mbox{$\phi_{i,k}(x) = \phi_{j,k}(y)$}. With this description
one easily checks
$$\colim_{I} M_i ~ = ~
\cup_{i \in I} \im(\psi_i: M_i \longrightarrow \colim_{I} M_i)$$
Now apply Cofinality of $\dim$ (see Theorem
\ref{properties of new definition}.4).\par\noindent
2.) It remains to show for \mbox{$i \in I $}
\begin{eqnarray}
\dim (\im(\psi_i)) & = &
\inf\left\{
\dim(\im(\phi_{i,j}: M_i \longrightarrow M_j)) \mid j \in I, i \le j
\right\}.
\label{eqn 2.660}
\end{eqnarray}
By assumption there is \mbox{$i_0 \in I $} with \mbox{$i \le i_0$}
such that \mbox{$\dim(\im(\phi_{i,i_0}))$} is finite.
Let $K_{i_0,j}$ be the kernel of the map
\mbox{$\im(\phi_{i,i_0}) \longrightarrow \im(\phi_{i,j})$}
induced by \mbox{$\phi_{i_0,j}$} for
\mbox{$i_0 \le j$} and $K_{i_0}$
be the kernel of the map \mbox{$\im(\phi_{i,i_0}) \longrightarrow
\im(\psi_i)$} induced by
\mbox{$\psi_{i_0}$}. Then
\mbox{$K_{i_0} ~ = ~ \cup_{j \in I, i_0 \le j} K_{i_0,j}$}
and hence by Cofinality (see Theorem
\ref{properties of new definition}.4)
$$\dim(K_{i_0}) ~ = ~
\sup\{\dim(K_{i_0,j}) \mid j \in I, i_0 \le j\}.$$
Since \mbox{$\dim(\im(\phi_{i,i_0}))$} is finite, we get
from Additivity (see Theorem
\ref{properties of new definition}.4)
\begin{eqnarray}
\dim(\im(\psi_i)) & = &
\dim\left(\im\left(\psi_{i_0}|_{\im(\phi_{i,i_0})}: \im(\phi_{i,i_0})
\longrightarrow \colim_I M_i\right)\right) \nonumber
\\
& = & \dim(\im(\phi_{i,i_0})) - \dim(K_{i_0}) \nonumber
\\
& = & \dim(\im(\phi_{i,i_0})) -
\sup\{\dim(K_{i_0,j}) \mid j \in I, i_0 \le j\} \nonumber
\\
& = &\inf\left\{ \dim(\im(\phi_{i,i_0})) -
\dim(K_{i_0,j}) \mid j \in I, i_0 \le j \right\} \nonumber
\\
& = & \inf\left\{\dim\left(\im(\phi_{i_0,j}|_{\im(\phi_{i,i_0})}:
\im(\phi_{i,i_0}) \longrightarrow
\im(\phi_{i,j}))\right)\mid j \in I, i_0 \le j \right\} \nonumber
\\
& = & \inf\left\{\dim\left(\im(\phi_{i,j})\right)
\mid j \in I, i_0 \le j \right\}
\label{eqn 2.678}
\end{eqnarray}
Given \mbox{$j_0 \in J$} with \mbox{$i \le j_0$}, there is
\mbox{$j \in I$} with \mbox{$i_0 \le j$} and \mbox{$j_0 \le j$}
and hence with
$$\dim(\im(\phi_{i,j_0})) \ge \dim(\im(\phi_{i,j})).$$
This implies
\begin{eqnarray}
\inf\{\dim(\im(\phi_{i,j})) \mid j \in J, i \le j\} & = &
\inf\{\dim(\im(\phi_{i,j})) \mid j \in J, i_0 \le j\}.
\label{eqn 2.699}
\end{eqnarray}
Now \ref{eqn 2.660} follows from \ref{eqn 2.678} and \ref{eqn 2.699}.
This finishes the proof of Theorem
\ref{dimension of colimits for arbitrary index sets}. \qed\par

\begin{example} \label{conditions in Theorem
 ref(dimension of colimits for arbitrary index sets) are necessary}
\em The condition in Theorem
\ref{dimension of colimits for arbitrary index sets}.2
that for each \mbox{$i \in I$}  there is \mbox{$i_0 \in I$} with
\mbox{$i \le i_0$} with \mbox{$\dim(\im(\phi_{i,i_0})) < \infty$}
is necessary as the following example shows.
Take \mbox{$I = \nn$}. Define
\mbox{$M_j = \oplus_{n = j}^{\infty} \cala$}
and
\mbox{$\phi_{j,k}:  \oplus_{m=j}^{\infty} \cala \longrightarrow
\oplus_{m=k}^{\infty} \cala $} to be the projection.
Then \mbox{$\dim(\im(\phi_{j,k})) = \infty$} for all \mbox{$j \le k$},
but \mbox{$\colim_I M_i$} is
trivial and hence has dimension zero. \qed \em
\end{example}

\begin{remark} \label{axiomatic point of view} \em
From an axiomatic point of view we have only needed
the following basic properties of $\cala$.
Namely, let $R$ be an associative ring with unit which has the
following properties

\begin{enumerate}
\item There is a dimension function $\dim$ which assigns to any
finitely generated projective $R$-module $P$ an element
$$\dim(P) \in [0,\infty)$$
such that \mbox{$\dim(P \oplus Q) = \dim(P) + \dim(Q)$} holds
and \mbox{$\dim(P)$} depends only on the isomorphism class of $P$;

\item If $K \subset P$ is a submodule of the finitely generated
projective $\cala$-module $P$, then $\overline{K}$ is a
direct summand in $P$. Moreover
$$\dim(\overline{K}) ~ = ~ \sup\{\dim(P) \mid P \subset K
\mbox{ finitely generated projective  }R\mbox{-submodule}\}.$$

\end{enumerate}

Then with Definition \ref{new definition of dimension}
Theorem \ref{properties of new definition} carries over to $R$.
One has essentially to copy the part of the proof which begins with
\ref{eqn 2.-2}. \par

An easy example where these axioms are satisfied is the case
where $R$ is a principal ideal domain and $\dim$ is
the usual rank of a finitely generated free $R$-module.
Then the extended dimension for
a $R$-module $M$ is just the dimension
of the rational vector space \mbox{$F \otimes_{R} M$}
for $F$ the quotient field of $R$. Notice that the case
of a von Neumann algebra $R = \cala$ is harder since
$\cala$ is not noetherian in general.\par

In Definition \ref{closure, K(M), P(M)} we have defined
$\bfT M$ and $\bfP M$ only for finitely generated $\cala$-modules
$M$ although the definition makes sense in general. The reason is
that the following definition for arbitrary $\cala$-modules
seems to be more appropriate
\begin{eqnarray}
\bfT M & := & \cup \{N \subset M \mid \dim(N) = 0\};
\label{new definition of TM}
\\
\bfP M & := & M/\bfT M.
\label{new definition of PM}
\end{eqnarray}
One easily checks using Theorem \ref{properties of new definition}.4
that $\bfT M$ is the largest submodule of $M$ with
trivial dimension and that these definitions
\ref{new definition of TM} and \ref{new definition of PM}
agree with Definition \ref{closure, K(M), P(M)} if
$M$ is finitely generated. One can show by example
that they do not agree if one applies them to arbitrary
$\cala$-modules. In the case of a principal ideal domain $R$
the torsion submodule
of a $R$-module $M$ is just $\bfT M$ in the sense
of definition \ref{new definition of TM}. \qed\em
\end{remark}


\typeout{--------------------- section 3 -----------------------}

\tit{Induction for group von Neumann algebras}
\label{Induction for group von Neumann algebras}

Next we investigate how the dimension behaves under induction.
Let \mbox{$i : \Delta \longrightarrow \Gamma$} be an
inclusion of groups. We claim that associated to $i$
there is a ring homomorphism of the group von Neumann algebras,
also denoted by
\begin{eqnarray}
i: \caln(\Delta) & \longrightarrow & \caln(\Gamma).
\label{i on group von Neumann algebra}
\end{eqnarray}
Recall that $\caln (\Delta)$ is the same as
the ring \mbox{$\calb(l^2(\Delta),l^2(\Delta))^{\Delta}$}
of bounded $\Delta$-equivariant operators
\mbox{$f : l^2(\Delta) \longrightarrow l^2(\Delta)$}.
Notice that
\mbox{$\cc\Gamma \otimes_{\cc\Delta} l^2(\Delta)$} can be
viewed as a dense subspace of $l^2(\Gamma)$ and that
$f$ defines a $\cc\Gamma$-homomorphism
\mbox{$\id \otimes_{\cc\Delta} f  :
\cc\Gamma \otimes_{\cc\Delta} l^2(\Delta)
\longrightarrow
\cc\Gamma \otimes_{\cc\Delta} l^2(\Delta)$}
which is bounded with respect to
the pre-Hilbert structure induced on
\mbox{$\cc\Gamma \otimes_{\cc\Delta} l^2(\Gamma)$}
from $l^2(\Gamma)$.
Hence \mbox{$\id \otimes_{\cc\Delta} f $} extends
to a $\Gamma$-equivariant bounded operator
\mbox{$i(f) : l^2(\Gamma) \longrightarrow l^2(\Gamma)$}.\par

Given a $\caln(\Delta)$-module $M$, define
{\em the induction with $i$} to be the $\caln(\Gamma)$-module
\begin{eqnarray}
i_{\ast}(M) &  := & \caln(\Gamma) \otimes_{\caln(\Delta)} M.
\label{induction of modules}
\end{eqnarray}
Obviously $i_{\ast}$ is a covariant functor from the category
of $\caln(\Delta)$-modules to the category of
$\caln(\Gamma)$-modules,
preserves direct sums and the
properties finitely generated and projective and sends
$\caln(\Delta)$ to $\caln(\Gamma)$.

\begin{theorem} \label{exactness and dimension under induction}
Let \mbox{$i : \Delta \longrightarrow \Gamma$} be an
injective group homomorphisms. Then:
\begin{enumerate}

\item Induction with $i$ is a faithfully flat functor
from the category of $\caln(\Delta)$-modules to the category
of $\caln(\Gamma)$-modules, i.e. a  sequence
of $\caln(\Delta)$-modules
\mbox{$M_0 \longrightarrow M_1 \longrightarrow M_2$} is exact
at $M_1$ if and only if the induced
sequence of $\caln(\Gamma)$-modules
\mbox{$i_{\ast}M_0 \longrightarrow i_{\ast}M_1
\longrightarrow i_{\ast}M_2$} is exact at $i_{\ast}M_1$;

\item For any $\caln(\Delta)$-module $M$ we have:
$$\dim_{\caln(\Delta)}(M) ~ = ~ \dim_{\caln(\Gamma)}(i_{\ast}M).$$

\end{enumerate}
\end{theorem}
\proof The proof consists of the following steps.
\par\noindent
Step 1: \hspace{5mm} $\dim_{\caln(\Delta)}(M) =
\dim_{\caln(\Gamma)}(i_{\ast}(M))$, provided $M$
is a finitely generated projective $\caln(\Delta)$-module.
\par
Let \mbox{$A \in M(n,n,\caln(\Delta))$} be a matrix such that
\mbox{$A = A^{\ast}$}, \mbox{$A^2 = A$} and the image of the
$\caln(\Delta)$-linear map
\mbox{$A: \caln(\Delta)^n \longrightarrow \caln(\Delta)^n$}
induced by right multiplication with $A$ is $\caln(\Delta)$-isomorphic to
$M$. Let $i(A)$ be the matrix in \mbox{$M(n,n,\caln(\Gamma))$} obtained
from $A$ by applying $i$ to each entry. Then
\mbox{$i(A) = i(A)^{\ast}$}, \mbox{$i(A)^2 = i(A)$} and the image of the
$\caln(\Gamma)$-linear map
\mbox{$i(A): \caln(\Gamma)^n \longrightarrow \caln(\Gamma)^n$}
induced by right multiplication with $i(A)$ is
$\caln(\Gamma)$-isomorphic to $i_{\ast}M$. Hence
we get from Definition
\ref{von Neumann dimension of a finitely generated projective
cala-module}
\begin{eqnarray*}
\dim_{\caln(\Delta)}(M) & = & \tr_{\caln(\Delta)}(A);
\\
\dim_{\caln(\Gamma)}(i_{\ast}M) & = & \tr_{\caln(\Gamma)}(i(A)).
\end{eqnarray*}
Therefore it suffices to show
\mbox{$\tr_{\caln(\Gamma)}(i(a))  =  \tr_{\caln(\Delta)}(a)$}
for \mbox{$a \in \caln(\Delta)$}.
This is an easy consequence of the Definition
\ref{definition of standard trace} of the standard trace.

\par\noindent
Step 2: \hspace{5mm} If $M$ is finitely presented $\caln(\Delta)$-module,
then
\begin{eqnarray*}
\dim_{\caln(\Delta)}(M) & = &
\dim_{\caln(\Gamma)}(i_{\ast}(M));
\\
\Tor_1^{\caln(\Delta)}(\caln(\Gamma),M) & = & 0.
\end{eqnarray*}
\par

Since $M$ is finitely presented, it splits  as
\mbox{$M = \bfT M \oplus \bfP M$} where
\mbox{$\bfP M$} is finitely generated projective and
there is an exact sequence
\mbox{$0 \longrightarrow \caln(\Delta)^n
\stackrel{f}{\longrightarrow}\caln(\Delta)^n
\longrightarrow \bfT M \longrightarrow 0$} with \mbox{$f^{\ast} = f$}
\cite[Theorem 1.2, Lemma 3.4]{Lueck (1995)}.
If we apply the
right exact functor induction with $i$ to it,
we get an exact sequence
\mbox{$\caln(\Gamma)^n \stackrel{i_{\ast}f}{\longrightarrow}
\caln(\Gamma)^n
\longrightarrow i_{\ast}\bfT M \longrightarrow 0$}
with \mbox{$(i_{\ast}f)^{\ast} = i_{\ast}f$}.
Because of Step 1, Additivity
(see Theorem \ref{properties of new definition}.4) and the definition of
$\Tor$ it suffices to show that $i_{\ast}f$ is injective.
Let $\nu$ be the functor introduced in
Theorem \ref{functor nu} or \cite[section 2]{Lueck (1995)}.
Then $i(\nu(f))$ is $\nu(i_{\ast}f)$.
Because $\nu$ respects  weak exactness
(see Theorem \ref{functor nu} or
\cite[Lemma 2.3]{Lueck (1995)}) $\nu(f)$
has dense image since
\mbox{$\caln(\Delta)^n \stackrel{i_{\ast}f}{\longrightarrow}
\caln(\Delta)^n \longrightarrow 0$} is weakly exact.
Then one easily checks that
\mbox{$\nu(i_{\ast}f) = i_{\ast}(\nu(f))$}
has dense image
since \mbox{$\cc\Gamma \otimes_{\cc\Delta} l^2(\Delta)$} is a dense
subspace of $l^2(\Gamma)$.
Since the kernel of a bounded operator of Hilbert spaces
is the orthogonal complement of the image of its adjoint
and $\nu(i_{\ast}f)$ is selfadjoint,
\mbox{$\nu(i_{\ast}f)$} is injective.
Since  $\nu^{-1}$ respects  exactness
(see Theorem \ref{functor nu} or
\cite[Lemma 2.3]{Lueck (1995)})
$i_{\ast}f$ is injective.
\par\noindent
Step 3:  \hspace{5mm}
$\Tor_1^{\caln(\Delta)}(\caln(\Gamma),M) = 0$
provided, $M$ is a finitely generated $\caln(\Delta)$-module.\par

Choose an exact sequence \mbox{$0 \longrightarrow K
\stackrel{g}{\longrightarrow} P
\longrightarrow M \longrightarrow 0$}
such that $P$ is a finitely generated
projective $\caln (\Delta)$-module.
The associated long exact sequence of
$\Tor$-groups shows
that \mbox{$\Tor_1^{\caln(\Delta)}(\caln(\Gamma),M)$}
is trivial if and  only if
\mbox{$i_{\ast}g : i_{\ast}K \longrightarrow i_{\ast}P$}
is injective. For each element $x$ in $i_{\ast}K$ there is
a finitely generated submodule $K^{\prime} \subset K$
such that $x$ lies in the image of the map
\mbox{$i_{\ast}K^{\prime} \longrightarrow i_{\ast}K$}
induced by the inclusion. Hence it suffices to show
for any finitely generated submodule
\mbox{$K^{\prime} \subset P$} that the inclusion induces
an injection \mbox{$i_{\ast}K^{\prime} \longrightarrow i_{\ast}P$}.
This follows since Step 2 applied to the finitely presented module
$P/K$ shows
\mbox{$\Tor_1^{\caln(\Delta)}(\caln(\Gamma),P/K) = 0$}.
\par\noindent
Step 4: \hspace{5mm} $i_{\ast}$ is an exact functor.\par

By standard homological algebra we have to show that
$\Tor^1_{\caln(\Delta)}(\caln(\Gamma),M) = 0$ is trivial
for all $\caln(\Delta)$-modules $M$. Notice that
$M$ is the colimit of the directed system of its
finitely generated submodules (directed by inclusion)
and that the functor $\Tor$ commutes in both variables
with colimits over directed systems
\cite[Proposition VI.1.3. on page 107]{Cartan-Eilenberg (1956)}.
Now the claim follows from Step 3.
\par\noindent
Step 5: \hspace{5mm} Let \mbox{$\{M_i \mid i \in I\}$}
be the directed system
of finitely generated submodules of
the $\caln(\Delta)$-module $M$. Then
\begin{eqnarray*}
\dim_{\caln(\Delta)}(M) & = &
\sup\{\dim_{\caln(\Delta)}(M_i) \mid i \in I\};
\\
\dim_{\caln(\Gamma)}(i_{\ast}M) & = &
\sup\{\dim_{\caln(\Gamma)}(i_{\ast}M_i) \mid i \in I\}.
\end{eqnarray*}
Because of Step 4 we can view $i_{\ast}M_i$
as a submodule of $i_{\ast}M$.
Now apply Cofinality
(see Theorem \ref{properties of new definition}.4).
\par\noindent
Step 6: \hspace{5mm} The second assertion of
Theorem \ref{exactness and dimension under induction} is true.

Because of Step 5 it suffices to prove the claim in the case
that $M$ is finitely generated because any module
is the colimit of the directed
system of its finitely generated submodules.
Choose an exact sequence \mbox{$0 \longrightarrow K
\stackrel{g}{\longrightarrow} P \longrightarrow M \longrightarrow
0$}
such that $P$ is a finitely generated
projective $\caln (\Delta)$-module. Because of Step 4 and
Additivity (see Theorem \ref{properties of new definition}) we get
\begin{eqnarray*}
\dim_{\caln(\Delta)}(M)  & = &
\dim_{\caln(\Delta)}(P) - \dim_{\caln(\Delta)}(K);
\\
\dim_{\caln(\Gamma)}(i_{\ast}M)  & = &
\dim_{\caln(\Gamma)}(i_{\ast}P) - \dim_{\caln(\Gamma)}(i_{\ast}K).
\end{eqnarray*}
Because of Step 1 it remains to prove
\begin{eqnarray*}
\dim_{\caln(\Delta)}(K) & = &
\dim_{\caln(\Gamma)}(i_{\ast}K).
\end{eqnarray*}
Because of Step 5 it suffices to treat the case where
$K \subset P$ is finitely generated. Since $\caln(\Delta)$
is semi-hereditary (see Theorem \ref{properties of new definition}.1)
$K$ is finitely generated projective and the claim follows
from Step 1.
\par\noindent
Step 7: \hspace{5mm} The first assertion of
Theorem \ref{exactness and dimension under induction} is true.
\par

Because we know already from Step 4 that $i_{\ast}$
is exact, it remains to prove for a $\caln(\Delta)$-module $M$
$$i_{\ast}M = 0 ~ \Longleftrightarrow ~ M = 0.$$
\par
Suppose \mbox{$i_{\ast}M = 0$}. In order to show \mbox{$M = 0$}
we have to prove for any
$\caln(\Delta)$-map \mbox{$f: \caln(\Delta) \longrightarrow M$}
that it is trivial. Let $K$ be the kernel of $f$.
Because $i_{\ast}$ is exact by
Step 4 and \mbox{$i_{\ast}M = 0$} by assumption,
the inclusion induces an isomorphism
\mbox{$i_{\ast}K \longrightarrow i_{\ast}\caln(\Delta) $}.
Since \mbox{$i_{\ast}K$} is a finitely generated $\caln(\Gamma)$-module
and $i_{\ast}$ is exact by Step 4,
there is a finitely generated submodule \mbox{$K' \subset K$}
such that the inclusion induces an isomorphism
\mbox{$i_{\ast}K' \longrightarrow i_{\ast}\caln(\Delta)$}. Let
\mbox{$\caln(\Delta)^m \longrightarrow K'$} be an epimorphism.
Let \mbox{$g: \caln(\Delta)^m \longrightarrow \caln(\Delta)$}
be the obvious composition. Because $i_{\ast}$ is exact by
Step 4 the induced map
\mbox{$i_{\ast}g: i_{\ast}\caln(\Delta)^m \longrightarrow
i_{\ast}\caln(\Delta)$} is surjective. Hence it remains
to prove that $g$ itself is surjective because then
\mbox{$K' = \caln(\Delta)$} and
the map \mbox{$f: \caln(\Delta) \longrightarrow M$} is trivial.
Since the functors $\nu^{-1}$ and $\nu$
of Theorem \ref{functor nu}
are exact we have to show for a $\Delta$-equivariant bounded
operator \mbox{$h: l^2(\Delta)^m\longrightarrow l^2(\Delta)$}
that $h$ is surjective if
\mbox{$i(h): l^2(\Gamma)^m\longrightarrow l^2(\Gamma)$} is surjective.
Let \mbox{$\{E_{\lambda} \mid \lambda \ge 0\}$}
be the spectral family of the positive operator
\mbox{$h\circ h^{\ast}$}.
Then \mbox{$\{i(E_{\lambda}) \mid \lambda \ge 0\}$}
is the spectral family of the positive operator
\mbox{$i(h)\circ i(h)^{\ast}$}.
Notice that $h$ resp. $i(h)$ is surjective if and only if
\mbox{$E_{\lambda} = 0$} resp. \mbox{$i(E_{\lambda}) = 0$}
for some \mbox{$\lambda > 0$}. Because \mbox{$E_{\lambda} = 0$}
resp. \mbox{$i(E_{\lambda}) = 0$} is equivalent to
\mbox{$\dim_{\caln(\Delta)}(\im(\nu^{-1}(E_{\lambda}))) = 0$}
resp. \mbox{$\dim_{\caln(\Gamma)}(\im(\nu^{-1}(i(E_{\lambda})))) = 0$} and
\mbox{$\im(\nu^{-1}(i(E_{\lambda}))) = i_{\ast}\im(\nu^{-1}(E_{\lambda}))$},
the claim follows from Step 6.
This finishes the proof of Theorem
\ref{exactness and dimension under induction}. \qed\par

The proof of the first two assertions of Theorem
\ref{exactness and dimension under induction}
would be obvious
if we would know
that $\caln (\Gamma)$ viewed as a $\caln (\Delta)$-module is
projective. Notice that this is a stronger statement than
proven in Theorem \ref{exactness and dimension under induction}.
One would have to show that the higher Ext-groups instead of the
Tor-groups vanish to get this stronger statement.
However, the proof for the Tor-groups does not go through
directly since the Ext-groups are not compatible with colimits.

\begin{lemma}
\label{invariants of caln(Gamma) otimes_{zzGamma} zz}
Let $H\subset \Gamma$ be a subgroup. Then
\begin{enumerate}
\item $\dim\left(\caln(\Gamma) \otimes_{\cc\Gamma} \cc[\Gamma/H]\right) =
|H|^{-1}$, where
$|H|^{-1}$ is defined to be zero if $H$ is infinite;

\item $\caln(\Gamma) \otimes_{\cc\Gamma} \cc[\Gamma/H]$ is trivial
if and only if $H$ is non-amenable;

\item If $\Gamma$ is infinite and $V$ is a $\cc\Gamma$-module
which is finite-dimensional over $\cc$, then
$$\dim(\caln(\Gamma) \otimes_{\cc\Gamma} V) ~ = ~ 0.$$
\end{enumerate}
\end{lemma}
\proof
3.) Since $V$ is finitely generated as $\cc\Gamma$-module
\mbox{$\caln(\Gamma)\otimes_{\cc\Gamma} V$} is a finitely generated
$\caln(\Gamma)$-module. Because of
Theorem \ref{properties of new definition} it suffices to show
that there is no $\caln(\Gamma)$-homomorphism from
\mbox{$\caln(\Gamma)\otimes_{\cc\Gamma} V$} to $\caln(\Gamma)$.
This is equivalent to the claim that there is no
$\cc\Gamma$-homomorphism from $V$ to $\caln(\Gamma)$. Since the map
\mbox{$\caln(\Gamma) \longrightarrow l^2(\Gamma)$} given by evaluation
at the unit element \mbox{$e \in \Gamma \subset l^2(\Gamma)$} is
$\Gamma$-equivariant and injective it suffices to show
that $l^2(\Gamma)$ contains no $\Gamma$-invariant linear subspace $W$
which is finite-dimensional as complex vector space.
Since any finite-dimensional topological vector space is
complete, $W$ is a Hilbert $\caln(\Gamma)$-submodule.
Let \mbox{$\pr : l^2(\Gamma) \longrightarrow l^2(\Gamma)$}
be an orthogonal $\Gamma$-equivariant projection onto
$W$. Then we get for any \mbox{$\gamma \in \Gamma$}
\begin{eqnarray}
\dim(W) & = & \langle \pr(\gamma),\gamma\rangle.
\label{eqn 3.758483}
\end{eqnarray}
Let \mbox{$\{v_1,v_2, \ldots, v_r\}$} be an orthonormal basis
for the Hilbert subspace \mbox{$W \subset l^2(\Gamma)$}.
For \mbox{$\gamma \in \Gamma$} we write
\mbox{$\pr(\gamma) = \sum_{i = 1}^r \lambda_i(\gamma) \cdot v_i$}.
We get from \mbox{$||\pr(\gamma)||^2 \le 1$}
\begin{eqnarray}
|\lambda_i(\gamma)| & \le  & 1.
\label{eqn 3.57122}
\end{eqnarray}
Given $\epsilon > 0$, we can choose $\gamma(\epsilon)$
satisfying
\begin{eqnarray}
\langle v_i,\gamma(\epsilon) \rangle & \le & r^{-1} \cdot \epsilon
\hspace{10mm} \mbox{ for } i= 1,2, \ldots, r.
\label{eqn 3.58002}
\end{eqnarray}
Now \ref{eqn 3.57122} and \ref{eqn 3.58002} imply
\begin{eqnarray}
\langle \pr(\gamma(\epsilon)), \gamma(\epsilon) \rangle & \le & \epsilon.
\label{eqn 3.68022}
\end{eqnarray}
Since \ref{eqn 3.68022} holds for all $\epsilon > 0$,
we conclude \mbox{$\dim(W) = 0$} and hence \mbox{$W = 0$} from
equation \ref{eqn 3.758483}.
\par\noindent
1. and 2.) If $i: H \longrightarrow \Gamma$ is the inclusion, then
\mbox{$i_{\ast}\left(\caln(H) \otimes_{\cc H} \cc\right)$} and
\mbox{$\caln(\Gamma) \otimes_{\cc\Gamma} \cc[\Gamma/H]$} are
isomorphic as $\caln(\Gamma)$-modules. Because of
Theorem \ref{exactness and dimension under induction}
it remains to treat the special case \mbox{$\Gamma = H$} for the
first two assertions. The first assertion follows from the third
for infinite $\Gamma$ and is obvious for finite $\Gamma$.
Next we prove the second assertion.\par

Let $S$ be a set of generators of $\Gamma$.
Then
\mbox{$\oplus_{s \in S} \zz\Gamma
\xrightarrow{\oplus_{s \in S} r_{s-1}}
\cc\Gamma \xrightarrow{\epsilon} \cc \longrightarrow 0$}
is exact where
\mbox{$\epsilon(\sum_{\gamma \in \Gamma} \lambda_{\gamma} \cdot \gamma) =
\sum_{\gamma \in \Gamma} \lambda_{\gamma}$}.
Since the tensor product is right exact, we obtain
an exact sequence
\mbox{$\oplus_{s \in S} \caln(\Gamma) \xrightarrow{\oplus_{s \in S}r_{s-1}}
\caln(\Gamma) \xrightarrow{\epsilon}
\caln(\Gamma) \otimes_{\cc\Gamma} \cc \longrightarrow 0$}.
Hence $\caln(\Gamma) \otimes_{\cc\Gamma} \cc$  is trivial
if and only if
\mbox{$\oplus_{s \in S} \caln(\Gamma) \xrightarrow{\oplus_{s \in S}r_{s-1}}
\caln(\Gamma)$}
is surjective. This is
equivalent to the existence of a finite subset $T \subset S$
such that
\mbox{$\oplus_{s \in T} \caln(\Gamma) \xrightarrow{\oplus_{s \in T}r_{s-1}}
\caln(\Gamma)$}
is surjective. Let $\Delta \subset \Gamma$ be the subgroup generated by
$T$. Then the map above is the induction with the inclusion
of $\Delta$ in $\Gamma$ applied to
\mbox{$\oplus_{t \in T} \caln(\Delta)
\xrightarrow{\oplus_{t \in T}r_{t-1}}
\caln(\Delta)$}.
Hence we conclude from Theorem
\ref{exactness and dimension under induction}.1
that $\caln(\Gamma) \otimes_{\cc\Gamma} \cc$ is trivial
if $\caln(\Delta) \otimes_{\cc\Delta} \cc$ is trivial
for some finitely generated subgroup $\Delta \subset \Gamma$.
Since $\Gamma$ is amenable if and only if each of its finitely generated
subgroups is amenable
\cite[Proposition 0.16 on page 14]{Paterson (1988)},
we can assume without loss of generality that
$\Gamma$ is finitely generated, i.e. $S$ is finite.
We can also assume that $S$ is symmetric, i.e. $s \in S$ implies
$s^{-1} \in S$. \par

Because the functor $\nu$ of Theorem
\ref{functor nu} is exact,
$\caln(\Gamma) \otimes_{\cc\Gamma} \cc$
is trivial if and only if the operator
\mbox{$f: \oplus_{s \in S} l^2(\Gamma) \xrightarrow{\oplus_{s \in S}r_{s-1}}
l^2(\Gamma)$}
is surjective. This is equivalent to the bijectivity of  the operator
$$\frac{1}{2\cdot|S|} f\circ f^*: ~
l^2(\Gamma) \xrightarrow{\id - \sum_{s \in S} \frac{1}{|S|}\cdot r_{s}}
l^2(\Gamma).$$
It is bijective if and only if the spectral radius of the operator
\mbox{$l^2(\Gamma) \xrightarrow{\sum_{s \in S} \frac{1}{|S|}\cdot r_{s}}
l^2(\Gamma)$}
is different from $1$. Since this operator is convolution with
a probability distribution
whose support contains $S$, namely
$$P: \Gamma \longrightarrow [0,1] \hspace{10mm} \gamma \mapsto
\left\{\begin{array}{ll} |S|^{-1} & \gamma \in S\\ 0 & \gamma \notin S
\end{array}\right.$$
the spectral radius is $1$ precisely if $\Gamma$ is amenable
\cite{Kesten (1959)}.  This finishes the proof
of Lemma  \ref{invariants of caln(Gamma) otimes_{zzGamma} zz}.\qed


\typeout{--------------------- section 4 -----------------------}

\tit{$L^2$-invariants for arbitrary $\Gamma$-spaces}
\label{L^2-invariants for arbitrary Gamma-spaces}

In this section we extend the notion of $L^2$-Betti numbers
for regular coverings of $CW$-complexes of finite type
(,i.e. with finite skeletons) with $\Gamma$ as group of deck
transformations to (compactly generated) topological spaces
with action of a (discrete) group $\Gamma$. We will continue with
using Notation \ref{basic notation}.

\begin{definition} \label{L^2-Betti numbers for arbitrary spaces}
Let $X$ be a (left) $\Gamma$-space and
$V$ be a $\cala$-$\zz\Gamma$-bimodule. Let
\mbox{$H_p^{\Gamma}(X;V)$} be {\em the singular
homology of $X$ with coefficients
in $V$}, i.e. the $\cala$-module given by the homology
of the $\cala$-chain complex
\mbox{$V \otimes_{\zz\Gamma} C_{\ast}^{\sing}(X)$}, where
$C_{\ast}^{\sing}(X)$ denotes the singular $\zz\Gamma$-chain complex of
$X$. Define the {\em $p$-th $L^2$-Betti number of $X$ with
coefficients in $V$} by
$$b_p^{(2)}(X;V) ~ :=
\dim_{\cala}\left(H_p^{\Gamma}(X;V)\right)
\hspace{5mm} \in [0,\infty ]. $$
and the {\em $p$-th $L^2$-Betti number of the group $\Gamma$} by
$$b_p^{(2)}(\Gamma) := b_p^{(2)}(E\Gamma;\caln(\Gamma)). \qed$$
\end{definition}

Next we compare cellular and singular chain complexes
and show that it does not matter whether we use
singular or cellular chain complexes
in the case that $X$ is a $\Gamma$-$CW$-complex.
For basic definitions and facts about
$\Gamma$-$CW$-complexes we refer for instance to
\cite[sections II.1 and II.2]{tom Dieck (1987)},
\cite[sections 1 and 2]{Lueck (1989)}.

\begin{lemma}
\label{comparision of cellular and singular chain complex}
Let $X$ be a $\Gamma$-$CW$-complex. Then there is a up to
$\zz\Gamma$-homotopy unique and in $X$ natural
$\zz\Gamma$-chain homotopy equivalence
$$f(X) : C^{\cell}_{\ast}(X) \longrightarrow C_{\ast}^{\sing}(X).$$
In particular we get for any $\cala$-$\zz\Gamma$-bimodule $V$
an in $X$ and $V$ natural isomorphism
$$H_p(V \otimes_{\zz\Gamma} C^{\cell}_{\ast}(X)) \xrightarrow{\cong}
H_p(V \otimes_{\zz\Gamma}C_{\ast}^{\sing}(X)).\qed$$
\end{lemma}
\proof Obviously the second assertion follows from the first
assertion which is proven as follows.\par

Let $Y$ be a $CW$-complex with cellular $\zz$-chain complex
$C^{\cell}_{\ast}$ and singular $\zz$-chain complex $C^{\sing}_{\ast}$.
We define a third (intermediate)
$\zz$-chain complex $C_{\ast}^{\inte}(Y)$
as the subcomplex of $C_{\ast}^{\sing}$ whose $n$-th chain module
is the kernel of
$$C_n^{\sing}(Y_n) \xrightarrow{c_n^{\sing}} C_{n-1}^{\sing}(Y_n)
\longrightarrow C_{n-1}^{\sing}(Y_n,Y_{n-1}).$$
There are an in $Y$ natural inclusion
and an in $Y$ natural epimorphism
of $\zz$-chain complexes
\begin{eqnarray}
i(Y) : C^{\inte}_{\ast}(Y) & \longrightarrow & C_{\ast}^{\sing}(Y);
\label{inclusion of C^{int} in C^{sing}}
\\
p(Y) : C^{\inte}_{\ast}(Y) & \longrightarrow & C_{\ast}^{\cell}(Y);
\label{epimorphism of C^{int} onto C^{cell}}
\end{eqnarray}
which induce isomorphisms on homology
\cite[page 263]{Lueck (1989)}. \par

If $\Gamma$ acts freely on the $\Gamma$-$CW$-complex $X$,
then $C^{\cell}_{\ast}$ and $C^{\sing}_{\ast}$ are free
$\zz\Gamma$-chain complexes and we get a $\zz\Gamma$-chain homotopy
equivalence well-defined up to $\zz\Gamma$-homotopy
from the fundamental theorem of homological algebra
and the fact that the chain maps
\ref{inclusion of C^{int} in C^{sing}}
and \ref{epimorphism of C^{int} onto C^{cell}} induce
isomorphisms on homology. In the general case one has to
go to the orbit category $\Or(\Gamma)$ and apply module theory
over this category instead of over $\zz\Gamma$.\par

The orbit category $\Or(\Gamma)$ has as objects homogenous spaces
and as morphisms $\Gamma$-maps. The $\Gamma$-$CW$-complex
$X$ defines a contravariant functor
$$\underline{X}: \Or(\Gamma) \longrightarrow \{CW-\COMPLEXES\}
\hspace{10mm} \Gamma/H \mapsto \map(\Gamma/H,X)^\Gamma = X^H.$$
Its composition with the functor
$C^{\cell}_{\ast}$, $C^{\inte}_{\ast}$
resp. $C^{\sing}_{\ast}$ from the category of
$CW$-complexes to the category of chain complexes yields
$\zz\Or(\Gamma)$-chain complexes, i.e. contravariant
functors
\begin{eqnarray*}
C^{\cell}_{\ast}(\underline{X}) : \Or(\Gamma) & \longrightarrow &
\{\zz-\CHAINCOMPLEXES\};
\\
C^{\inte}_{\ast}(\underline{X}) : \Or(\Gamma) & \longrightarrow &
\{\zz-\CHAINCOMPLEXES\};
\\
C^{\sing}_{\ast}(\underline{X}) : \Or(\Gamma) & \longrightarrow &
\{\zz-\CHAINCOMPLEXES\}.
\end{eqnarray*}
We obtain natural transformations from the natural chain maps
\ref{inclusion of C^{int} in C^{sing}} and
\ref{epimorphism of C^{int} onto C^{cell}}
\begin{eqnarray}
i(\underline{X}) : C^{\inte}_{\ast}(\underline{X}) &
\longrightarrow & C_{\ast}^{\sing}(\underline{X});
\label{inclusion of C^{int} in C^{sing} over the orbit category}
\\
p(\underline{X}) : C^{\inte}_{\ast}(\underline{X}) &
\longrightarrow &C_{\ast}^{\cell}(\underline{X});
\label{epimorphism of C^{int}
 onto C^{cell} over the orbit category}
\end{eqnarray}
which induce isomorphisms on homology. Hence it suffices
to show that $C^{\cell}_{\ast}(\underline{X})$ and
$C^{\sing}_{\ast}(\underline{X})$ are free and hence projective
in the sense of \cite[Definition 9.17]{Lueck (1989)}
because then we obtain a homotopy equivalence
of $\zz\Or(\Gamma)$-chain complexes from
$C^{\cell}_{\ast}(X)$ to $C_{\ast}^{\sing}(X)$
\cite[Lemma 11.3]{Lueck (1989)}
whose evaluation at $\Gamma/1$
is the desired $\zz\Gamma$-chain homotopy equivalence.
The proofs that these two chain complexes are free are
simple versions of the arguments in
\cite[Lemma 13.2]{Lueck (1989)}. Notice that in \cite{Lueck (1989)}
the $\Gamma$-$CW$-complex is required to be proper, but this
condition is needed there only because there $\Gamma$ is assumed to be a
Lie group and universal coverings are built in,
and can be dropped in the discrete case.
\qed

\begin{remark}
\label{comparision with all definitions}
\em Originally the $L^2$-Betti numbers of a regular covering
\mbox{$\overline{M} \longrightarrow M$} of a closed
Riemannian manifold with group of deck transformations
$\Gamma$ were defined by Atiyah \cite{Atiyah (1976)}
in terms of the heat kernel as explained in
\ref{L^2-Betti numbers and Laplace transform}
in the introduction. It follows from the $L^2$-Hodge-deRham
theorem \cite{Dodziuk (1977)} that this analytic definition
agrees with the combinatorial
definition of $b_p^{(2)}(\overline{X})$
in terms of the associated cellular
$L^2$-chain complex and the von
Neumann dimension of finitely generated
Hilbert $\caln(\Gamma)$-modules for
a triangulation $X$ of $M$.
Because of Lemma
\ref{comparision of cellular and singular chain complex}
this combinatorial definition agrees with the Definition
\ref{L^2-Betti numbers for arbitrary spaces}. \par

Analogously to the case of $L^2$-Betti numbers
we will extend the notion of Novikov-Shubin invariants
for regular coverings of compact Riemannian manifolds
to arbitary $\Gamma$-spaces and prove that they are positive for
the universal covering of a
aspherical closed manifold with elementary-amenable fundamental group
in another paper.
\qed
\end{remark}

The next results are well-known in the case where
$X$ is a regular covering of a $CW$-complex of finite type.
We call a map \mbox{$g: Y \longrightarrow Z$}
{\em homologically $n$-connected for $n \ge 1$}
if the map induced on singular homology with complex coefficients
\mbox{$g_{\ast}: H^{\sing}_k(Y;\cc) \longrightarrow H^{\sing}_k(Y;\cc)$}
is bijective for $k < n$ and surjective for $k = n$.
The map $g$ is called a {\em weak homolopy equivalence} if
it is $n$-connected for all $n \ge 1$.

\begin{lemma}
\label{weak homotopy invariance}
Let \mbox{$f: X \longrightarrow Y$} be a $\Gamma$-map
and let $V$ be a $\cala$-$\zz\Gamma$-bimodule.
\begin{enumerate}
\item Suppose for $n \ge 1$
that for each subgroup $H \subset \Gamma$ the induced
map \mbox{$f^H: X^H \longrightarrow Y^H$} is homologically $n$-connected.
Then the map induced by $f$
$$f_*: H_p^{\Gamma}(X;V) \longrightarrow H_p^{\Gamma}(Y;V)$$
is bijective for $p < n$ and surjective for $p =n$
and we get
\begin{eqnarray*}
b_p^{(2)}(X;V) & =   & b_p^{(2)}(Y;V)
\hspace{10mm} \mbox{ for } p < n;
\\
b_p^{(2)}(X;V) & \ge & b_p^{(2)}(Y;V)
\hspace{10mm} \mbox{ for } p = n;
\end{eqnarray*}

\item Suppose
such that for each subgroup $H \subset \Gamma$ the induced
map \mbox{$f^H: X^H \longrightarrow Y^H$}
is a weak homology equivalence.
Then for all $p \ge 0$ the map induced by $f$
$$f_*: H_p^{\Gamma}(X;V) \longrightarrow H_p^{\Gamma}(Y;V)$$
is bijective and we get
$$b_p^{(2)}(X;V) ~ = ~ b_p^{(2)}(Y;V). \qed$$
\end{enumerate}
\end{lemma}
\proof We give only the proof of the second assertion,
the one of the first assertion is an elementary modification.
The map $f$ induces a homotopy equivalence
of $\zz\Or(\Gamma)$-chain complexes
\mbox{$C^{\sing}_{\ast}(\underline{f}): C^{\sing}_{\ast}(\underline{X})
\longrightarrow C^{\sing}_{\ast}(\underline{Y})$}
in the notation of the proof of Lemma
\ref{comparision of cellular and singular chain complex}
since the singular $\zz\Or(\Gamma)$-chain complexes of $\underline{X}$
and $\underline{Y}$ are free in the sense of
\cite[Definition 9.17]{Lueck (1989)}(see
\cite[Lemma 11.3]{Lueck (1989)}). Its evaluation at $\Gamma/1$
is a $\zz\Gamma$-chain equivalence. Hence $f$ induces
a chain equivalence
$$V \otimes_{\zz\Gamma} C^{\sing}_{\ast}(f) :
V \otimes_{\zz\Gamma} C^{\sing}_{\ast}(X)
\longrightarrow V \otimes_{\zz\Gamma} C^{\sing}_{\ast}(Y)$$
and Lemma \ref{weak homotopy invariance} follows.\qed

We get as a direct consequence from Theorem
\ref{exactness and dimension under induction}

\begin{theorem}
\label{induction for spaces}
Let \mbox{$i: \Delta\longrightarrow \Gamma$} be an inclusion of groups
and let $X$ be a $\Delta$-space.
Then
\begin{eqnarray*}
H_p^{\Gamma}(\Gamma\times_{\Delta} X;\caln(\Gamma)) & = &
i_*H_p^{\Delta}(X;\caln(\Delta));
\\
b_p^{(2)}(\Gamma\times_{\Delta} X;\caln(\Gamma)) & = &
b_p^{(2)}(X;\caln(\Delta)). \qed
\end{eqnarray*}
\end{theorem}

\begin{theorem}
\label{L^2-invariants in dimension 0}
Let $X$ be a path-connected $\Gamma$-space. Then
\begin{enumerate}
\item There is an isomorphism of $\caln(\Gamma)$-modules
\mbox{$H_0^{\Gamma}(X;\caln(\Gamma)) ~ \cong
~ \caln(\Gamma)\otimes_{\zz\Gamma} \cc$};

\item $b^{(2)}_0(X;\caln(\Gamma)) = |\Gamma|^{-1}$, where
$|\Gamma|^{-1}$ is defined to be zero if the order $|\Gamma|$ of
$\Gamma$ is infinite;

\item $H_0^{\Gamma}(X;\caln(\Gamma))$ is trivial if and only if
$\Gamma$ is non-amenable.
\end{enumerate}
\end{theorem}
\proof
The first assertion follows from the fact that
\mbox{$C_1^{\sing}(X) \longrightarrow C_0^{\sing}(X) \longrightarrow \zz
\longrightarrow 0$}
is an exact sequence of $\zz\Gamma$-modules
and the tensor product is right exact.
The other two assertions follow from Lemma
\ref{invariants of caln(Gamma) otimes_{zzGamma} zz}. \qed\par

\begin{remark} \label{alpha(X) = infty^+}\em
Let \mbox{$\widetilde{M} \longrightarrow M$} be the universal covering of
a closed Riemannian manifold with fundamental group $\pi$.
Brooks \cite{Brooks (1981)} has shown that the analytic
Laplace operator $\Delta_0$ on $\widetilde{M}$ in dimension zero
has zero not in its spectrum if and only if
$\pi$ is non-amenable. Now $\Delta_0$ has zero not in its spectrum
if and only if
\mbox{$H_0^{\pi}(\widetilde{M},\caln(\pi))$} is trivial because of
\cite[paragraph after Definition 3.11, Theorem 6.1]{Lueck (1995)}
and the fact that the analytic and combinatorial spectral density function
are dilatationally equivalent \cite{Efremov(1991)}.
Hence Theorem \ref{L^2-invariants in dimension 0}
generalizes the result of Brooks. Notice that both
Brook's and our proof use \cite{Kesten (1959)}. Compare also with
the result  \cite[Corollary III.2.4 on page 188]{Guichardet (1980)} 
that a group $\Gamma$ is non-amenable if and 
only if $H^1(\Gamma,l^2(\Gamma))$ is Hausdorff.\qed\em
\end{remark}

\begin{remark} \em
\label{comparision with the definition of Cheeger and Gromov}
Next we compare our approach with the one in
\cite[section 2]{Cheeger-Gromov (1986)}. We begin with the case
of a countable simplicial complex $X$ with free simplicial
$\Gamma$-action. Then for any exhaustion
\mbox{$X_0 \subset X_1 \subset X_2 \subset \ldots X$}
by $\Gamma$-equivariant simplicial subcomplexes for which
$X/\Gamma$ is compact, the $p$-th $L^2$-Betti number
in the sense and notation of
\cite[2.8 on page 198]{Cheeger-Gromov (1986)} is given by
$$b_p^{(2)}(X:\Gamma) ~ = ~
\lim_{j \to \infty} \lim_{k \to \infty}
\dim_{\caln(\Gamma}\left(\im(\overline{H}^p_{(2)}(X_k:\Gamma)
\xrightarrow{i_{j,k}^{\ast}} \overline{H}^p_{(2)}(X_j:\Gamma)\right),
$$
where \mbox{$i_{j,k}: X_j \longrightarrow X_k$} is the inclusion for
\mbox{$j \le k$}.
We get from \cite[Lemma 1.3]{Lueck (1995)} and
Lemma \ref{comparision of cellular and singular chain complex}
\begin{center}
$\dim_{\caln(\Gamma)}\left(\im(\overline{H}^p_{(2)}(X_k:\Gamma)
\xrightarrow{i_{j,k}^{\ast}} \overline{H}^p_{(2)}(X_j:\Gamma)\right)$
\\
$~ = ~
\dim_{\caln(\Gamma)}\left(\im(H_p^{\Gamma}(X_j;\caln(\Gamma))
\xrightarrow{(i_{j,k})_{\ast}}
H_p^{\Gamma}(X_k;\caln(\Gamma))\right).$
\end{center}
Hence we conclude from Theorem
\ref{dimension of colimits for arbitrary index sets}
that the definitions in
\cite[2.8 on page 198]{Cheeger-Gromov (1986)}
and in
\ref{L^2-Betti numbers for arbitrary spaces}
agree:
\begin{eqnarray}
b_p^{(2)}(X:\Gamma) &  = & b_p^{(2)}(X;\caln(\Gamma)).
\label{Cheeger-Gromov equals ours in the free case}
\end{eqnarray}
If $\Gamma$ is countable and $X$ is a countable simplicial
complex with simplicial $\Gamma$-action, then
by \cite[Proposition 2.2 on page 198]{Cheeger-Gromov (1986)}
and by \ref{Cheeger-Gromov equals ours in the free case}
\begin{eqnarray}
b_p^{(2)}(X:\Gamma) &  = & b_p^{(2)}(E\Gamma \times X:\Gamma);
\\
b_p^{(2)}(X:\Gamma) &  = & b_p^{(2)}(E\Gamma \times X;\caln(\Gamma)).
\label{Cheeger-Gromov equals ours}
\end{eqnarray}

Cheeger and Gromov \cite[Section 2]{Cheeger-Gromov (1986)}
define $L^2$-cohomology and $L^2$-Betti numbers of a
$\Gamma$-space $X$ by considering the category
whose objects are $\Gamma$-maps \mbox{$f: Y \longrightarrow X$} for
a simplicial complex $Y$ with cocompact free simplicial
$\Gamma$-action and then using inverse limits to
extend the classical notions for
finite free $\Gamma$-$CW$-complexes such as $Y$
to $X$. Our approach avoids the technical difficulties concerning
inverse limits and is closer to standard notions, the only non-standard
part is the verification of the properties of the extended
dimension function (Theorem \ref{properties of new definition}
and Theorem \ref{exactness and dimension under induction}).
\em \end{remark}


\typeout{--------------------- section 5  -----------------------}

\tit{Amenable groups}
\label{Amenable groups}

In this section we investigate amenable groups.
For information about amenable
groups we refer for instance to \cite{Paterson (1988)}.
The main technical result
of this section is the next lemma whose proof uses ideas
of the proof of \cite[Lemma 3.1 on page 203]{Cheeger-Gromov (1986)}.

\begin{theorem}
\label{dimension of higher Tor-s vanish in the amenable case}
Let $\Gamma$ be amenable and
$M$ be a $\cc\Gamma$-module. Then
\begin{eqnarray*}
\dim_{\caln(\Gamma)}\left(\Tor^{\cc\Gamma}_p(M,\caln( \Gamma))\right)
& = & 0
\hspace{10mm} \mbox{ for } p \ge 1,
\end{eqnarray*}
where we consider $\caln(\Gamma)$ as a
$\caln(\Gamma)$-$\cc\Gamma$-bimodule.
\end{theorem}
\proof
Step 1: \hspace{1mm}
If $M$ is a finitely presented $\cc\Gamma$-module, then
\mbox{$\dim_{\caln(\Gamma)}\left(\Tor^{\cc\Gamma}_1
(M,\caln( \Gamma))\right) = 0$}.
\par
Choose a finite presentation
$$\oplus_{i=1}^m \cc\Gamma \xrightarrow{f} \oplus_{i=1}^n \cc\Gamma
\xrightarrow{p} M \longrightarrow 0.$$
For an element \mbox{$u = \sum_{\gamma} \lambda_{\gamma} \cdot \gamma$}
in \mbox{$l^2(\Gamma)$} define its {\em support}
\begin{eqnarray*}
\supp(u) & := &\{\gamma \in \Gamma \mid \lambda_{\gamma} \not= 0\}
\hspace{5mm}
\subset \Gamma.
\end{eqnarray*}
Let \mbox{$B \in M(m,n,\cc\Gamma)$} be the matrix describing $f$, i.e.
the component \mbox{$f_{i,j} : \cc\Gamma \longrightarrow \cc\Gamma$}
is given by right multiplication with $b_{i,j}$. Define
the finite subset $S$ by
\begin{eqnarray*}
S & := & \{\gamma \mid \gamma \mbox{ or } \gamma^{-1} \in \bigcup_{i,j}
\supp(b_{i,j})\}.
\end{eqnarray*}
Let \mbox{$f^{(2)}: \oplus_{i=1}^m l^2(\Gamma) \longrightarrow
\oplus_{i=1}^n l^2(\Gamma)$} be
the bounded $\Gamma$-equivariant operator
induced by $f$. Denote by $K$ the $\Gamma$-invariant linear subspace
of \mbox{$\oplus_{i=1}^m l^2(\Gamma)$} which is the image of
the kernel of $f$ under the canonical inclusion
\mbox{$k: \oplus_{i=1}^m \cc\Gamma \longrightarrow
\oplus_{i=1}^m l^2(\Gamma)$}.
Next we show for the closure $\overline{K}$ of $K$
\begin{eqnarray}
\overline{K} & = & \ker(f^{(2)}).
\label{eqn 5.1}
\end{eqnarray}
Let \mbox{$\pr: \oplus_{i=1}^m l^2(\Gamma) \longrightarrow
\oplus_{i=1}^m l^2(\Gamma) $} be the orthogonal projection onto
the closed $\Gamma$-invariant
subspace \mbox{$\overline{K}{^\perp} \cap \ker(f^{(2)})$}.
The von Neumann dimension of \mbox{$\im(\pr)$} is zero
if and only if $\pr$ itself is zero. Hence \ref{eqn 5.1}
will follow if we can prove
\begin{eqnarray}
\tr_{\caln(\Gamma)}(\pr) & = & 0.
\label{eqn 5.2}
\end{eqnarray}
Let $\epsilon > 0$ be given. Since $\Gamma$ is amenable,
there is a finite non-empty subset \mbox{$A \subset \Gamma$} satisfying
\cite[Theorem F.6.8 on page 308]{Benedetti-Petronio(1992)}
\begin{eqnarray}
\frac{|\partial_SA|}{|A|} & \le & \epsilon,
\label{eqn 5.3}
\end{eqnarray}
where \mbox{$\partial_SA$} is defined by
\mbox{$\{a \in A \mid \mbox{ there is } s \in S \mbox{ with }
as \notin A\}$}. Define 
$$\Delta ~ := ~ 
\left\{\gamma \in \Gamma \mid \gamma \in \partial_SA \mbox{ or }
\gamma s \in \partial_SA \mbox{ for some } s \in S\right\} ~ = ~ 
\partial_SA \cup \left(\cup_{s \in S} (\partial_SA)s\right).$$

Let \mbox{$\pr_A: l^2(\Gamma) \longrightarrow \l^2(\Gamma) $} be the
projection sending
\mbox{$\sum_{\gamma\in\Gamma} \lambda_{\gamma} \cdot \gamma$}
to \mbox{$\sum_{\gamma \in A} \lambda_{\gamma} \cdot \gamma$}.
Define
\mbox{$\pr_{\Delta}: l^2(\Gamma) \longrightarrow \l^2(\Gamma)$}
analogously.
Next we show for \mbox{$s \in S$} and \mbox{$u \in l^2(\Gamma)$}
\begin{eqnarray}
\pr_A \circ r_s(u) & = & r_s \circ \pr_A(u)
\hspace{10mm} \mbox{ ,if } \pr_{\Delta}(u) = 0,
\label{eqn 5.4}
\end{eqnarray}
where \mbox{$r_s: l^2(\Gamma) \longrightarrow l^2(\Gamma)$} is right
multiplication with $s$. Since \mbox{$s \in S$} implies
\mbox{$s^{-1}\in S$},
we get the following equality of subsets of $\Gamma$
\begin{eqnarray*}
\{\gamma \in \Gamma \mid \gamma s \in A,
\gamma \notin \Delta \} & = &
\{\gamma \in \Gamma \mid \gamma \in A, \gamma \notin \Delta\}.
\end{eqnarray*}
Now \ref{eqn 5.4} follows from the following calculation
for \mbox{$u = \sum_{\gamma \in \Gamma, \gamma \notin \Delta}~
\lambda_{\gamma} \cdot \gamma ~ \in l^2(\Gamma)$}
\begin{eqnarray*}
\pr_A \circ r_s(u) & = &
\sum_{\gamma s \in A, \gamma \notin \Delta} ~
\lambda_{\gamma} \cdot \gamma s
\\
 & = &
\sum_{\gamma\in A, \gamma \notin \Delta} ~
\lambda_{\gamma} \cdot \gamma s
\\
& = & \left(\sum_{\gamma\in A, \gamma \notin \Delta} ~
\lambda_{\gamma} \cdot \gamma\right)\cdot  s
\\
 & = &
r_s \circ \pr_{A}(u).
\end{eqnarray*}
We have defined $S$ such that each entry in the matrix $B$ describing
$f$ is a linear combination of elements in $S$.
Hence \ref{eqn 5.4} implies
\begin{eqnarray*}
\left(\oplus_{j=1}^n \pr_A\right) \circ f^{(2)}(u) & = &
f^{(2)} \circ \left(\oplus_{i=1}^m \pr_A\right)(u)
\hspace{10mm} \mbox{ ,if } \pr_{\Delta}(u_i) = 0
\mbox{ for } i=1,2 \ldots ,m.
\end{eqnarray*}
Notice that the image of \mbox{$\oplus_{i=1}^m \pr_A$}
lies in \mbox{$\oplus_{i=1}^m \cc\Gamma$}.
We conclude
\begin{eqnarray*}
\oplus_{i=1}^m \pr_A(u) & \in & K
\hspace{10mm} \mbox{, if } u \in \ker(f^{(2)}),
\pr_{\Delta}(u_i) = 0
\mbox{ for } i=1,2 \ldots ,m.
\end{eqnarray*}
This shows
\begin{eqnarray*}
\left(\pr \circ \oplus_{i=1}^m \pr_A\right)\left(\ker(f^{(2)}) \cap
\oplus_{i=1}^m \ker(\pr_{\Delta})\right) & = & 0.
\end{eqnarray*}
Since \mbox{$\ker(\pr_{\Delta})$} has complex codimension
$|\Delta|$ in $l^2(\Gamma)$ and $|\Delta|  \le  (|S| +1) \cdot |\partial_SA|$,
we conclude for the complex dimension $\dim_{\cc}$
of complex vector spaces
\begin{eqnarray}
\dim_{\cc}\left(\left(\pr \circ \oplus_{i=1}^m
\pr_A\right)\left(\ker(f^{(2)})\right)\right) & \le &
m \cdot (|S| +1) \cdot |\partial_SA|.
\label{eqn 5.6}
\end{eqnarray}
Since \mbox{$\pr \circ\pr_A$} is an endomorphism of  Hilbert spaces
with finite-dimensional image, it is trace-class and its trace
$\tr_{\cc}(\pr \circ\pr_A)$ is defined. We get
\begin{eqnarray}
\tr_{\caln(\Gamma)}(\pr) & \le &
\frac{\tr_{\cc}\left(\pr \circ \oplus_{i=1}^m\pr_A\right)}{|A|}
\label{eqn 5.7}
\end{eqnarray}
from the following computation for
\mbox{$e \in \Gamma \subset l^2(\Gamma)$}
the unit element
\begin{eqnarray*}
\tr_{\caln(\Gamma )}(\pr) & = &
\sum_{i=1}^m \langle\pr_{i,i}(e),e\rangle
\\
& = &
\frac{1}{|A|} \sum_{i=1}^m
|A|\cdot \langle\pr_{i,i}(e),e\rangle
\\
& = &
\frac{1}{|A|} \sum_{i=1}^m
\sum_{\gamma \in A} \langle\pr_{i,i}(\gamma),\gamma\rangle
\\
& = &
\frac{1}{|A|} \sum_{i=1}^m \sum_{\gamma \in A}
\langle \pr_{i,i} \circ \pr_A(\gamma),\gamma\rangle
\\
& = &
\frac{1}{|A|} \sum_{i=1}^m \sum_{\gamma \in \Gamma}
\langle \pr_{i,i} \circ \pr_A(\gamma),\gamma\rangle
\\
& = &
\frac{1}{|A|} \sum_{i=1}^m \tr_{\cc}(\pr_{i,i} \circ \pr_A)
\\
& = &
\frac{1}{|A|}
\tr_{\cc}\left(\pr \circ (\oplus_{i=1}^m\pr_A)\right).
\end{eqnarray*}
If $H$ is a Hilbert space and
\mbox{$f: H \longrightarrow H$} is a bounded operator with
finite-dimensional image, then
\mbox{$\tr_{\cc}(f) \le ||f|| \cdot \dim_{\cc}\left(f(\im(f))\right)$}.
Since the image of $\pr$ is contained in
$\ker(f^{(2)})$ and $\pr$ and $\pr_A$
have operator norm $1$, we conclude
\begin{eqnarray}
\tr_{\cc}\left(\pr \circ \oplus_{i=1}^m\pr_A\right)
& \le &
\dim_{\cc}\left(\left(\pr \circ \oplus_{i=1}^m
\pr_A\right)\left(\ker(f^{(2)})\right)\right).
\label{eqn 5.8}
\end{eqnarray}
Equations \ref{eqn 5.3}, \ref{eqn 5.6}, \ref{eqn 5.7} and \ref{eqn 5.8}
imply
\begin{eqnarray*}
\tr_{\caln(\Gamma)}(\pr) \le m \cdot (|S| +1) \cdot \cdot \epsilon.
\end{eqnarray*}
Since this holds for all $\epsilon > 0$, we get \ref{eqn 5.2}
and hence \ref{eqn 5.1} is true.\par

Let \mbox{$\pr_{\overline{K}}: \oplus_{i=1}^m l^2(\Gamma)
\longrightarrow \oplus_{i=1}^m l^2(\Gamma)$}
be the projection onto $\overline{K}$.
Let \mbox{$i: \ker(f) \longrightarrow \oplus_{i=1}^m \cc\Gamma$}
be the inclusion. It induces a map
$$j: \caln(\Gamma) \otimes_{\cc\Gamma} \ker(f)
\xrightarrow{\id \otimes_{\cc\Gamma} i}
\caln(\Gamma) \otimes_{\cc\Gamma} \oplus_{i=1}^m \cc\Gamma
\xrightarrow{\cong}  \oplus_{i=1}^m \caln(\Gamma).$$
Next we want to show
\begin{eqnarray}
\im(\nu^{-1}(\pr_{\overline{K}})) & = &
\overline{\im(j)}.
\label{eqn 5.9}
\end{eqnarray}
Let \mbox{$x \in \ker(f)$}. Then
\begin{eqnarray}
(\id - \nu^{-1}(\pr_{\overline{K}}))\circ j(1 \otimes x)
& = &
(\id - \pr_{\overline{K}})\circ k \circ i(x), \label{eqn 5.364655}
\end{eqnarray}
where \mbox{$k: \oplus_{i=1}^m \cc\Gamma
\longrightarrow  \oplus_{i=1}^m l^2(\Gamma)$} is the inclusion.
Since \mbox{$(\id - \pr_{\overline{K}})$} is trivial on $K$
we get \mbox{$(\id - \pr_{\overline{K}})\circ k \circ i = 0$}.
Now we conclude from \ref{eqn 5.364655} that
\mbox{$\im(j)\subset \ker(\id - \nu^{-1}(\pr_{\overline{K}}))$}
and hence
\mbox{$\im(j) \subset \im(\nu^{-1}(\pr_{\overline{K}}))$} holds.
This shows
\mbox{$\overline{\im(j)} \subset \im(\nu^{-1}(\pr_{\overline{K}}))$}.
It remains to prove for any $\caln(\Gamma)$-map
\mbox{$g: \oplus_{i=1}^m \caln(\Gamma) \longrightarrow
\oplus_{i=1}^m \caln(\Gamma)$} with
\mbox{$\im(j) \subset \ker(g)$} that
\mbox{$g \circ \nu^{-1}(\pr_{\overline{K}})$} is trivial.
Obviously \mbox{$K \subset \ker(\nu(g))$}. Since
\mbox{$\ker(\nu(g))$} is a closed subspace,
we get \mbox{$\overline{K} \subset \ker(\nu(g))$}.
We conclude \mbox{$\nu(g)\circ \pr_{\overline{K}}  = 0$}
and hence \mbox{$g \circ \nu^{-1}(\pr_{\overline{K}}) = 0$}.
This finishes the proof of \ref{eqn 5.9}.\par

Since $\nu^{-1} $ preserves exactness by Theorem \ref{functor nu}
and \mbox{$\id \otimes_{\cc\Gamma} f = \nu^{-1}(f^{(2)})$}, we conclude from
\ref{eqn 5.1} and \ref{eqn 5.9} that the sequence
$$\caln(\Gamma) \otimes_{\cc\Gamma} \ker(f)
\xrightarrow{\id \otimes_{\cc\Gamma} i}
\caln(\Gamma) \otimes_{\cc\Gamma} \oplus_{i=1}^m \cc\Gamma
\xrightarrow{\id \otimes_{\cc\Gamma} f}
\caln(\Gamma) \otimes_{\cc\Gamma} \oplus_{i=1}^m \cc\Gamma$$
is weakly exact. Continuity of the dimension function
(see Theorem \ref{properties of new definition}.4)
implies
\begin{eqnarray*}
\dim_{\caln(\Gamma)}\left(\ker(\id \otimes_{\cc\Gamma} f)/
\im(\id \otimes_{\cc\Gamma} i)\right) & = & 0.
\end{eqnarray*}
Since \mbox{$\Tor_1^{\cc\Gamma}(M;\caln(\Gamma)) =
\ker(\id \otimes_{\cc\Gamma} f)/\im(\id \otimes_{\cc\Gamma} i)$}
holds, Step 1 follows.
\par\noindent
Step 2:  \hspace{5mm}
If $M$ is a $\cc\Gamma$-module, then
\mbox{$\dim_{\caln(\Gamma)}
\left(\Tor^{\cc\Gamma}_1(M,\caln( \Gamma))\right) = 0$}.
\par
Obviously $M$ is the union of its finitely generated submodules.
Any finitely generated module $M$ is a colimit over a directed
system of finitely presented modules, namely, choose
an epimorphism from a finitely generated free
module F to M with kernel K. Since K is the union of
its finitely generated submodules, M is the colimit of the directed
system $F/L$ where $L$ runs over the finitely generated submodules
of $K$. The functor $\Tor$ commutes in both variables
with colimits over directed systems
\cite[Proposition VI.1.3. on page 107]{Cartan-Eilenberg (1956)}.
Now the claim follows from
Step 1 and Theorem \ref{dimension of colimits for arbitrary index sets}.
\par\noindent
Step 3:  \hspace{5mm} If $M$ is a $\cc\Gamma$-module, then
\mbox{$\dim_{\caln(\Gamma)}
\left(\Tor^{\cc\Gamma}_p(M,\caln( \Gamma))\right) = 0$}
for all \mbox{$ p \ge 1$}
\par
We use induction over $p \ge 1$. The induction begin
is already done in Step 2. Choose an exact sequence
\textshortexactsequence{N}{}{F}{}{M}
of $\caln(\Gamma)$-modules such that $F$ is free.
Then we obtain an isomorphism
\mbox{$\Tor_p^{\cc\Gamma}(M,\caln(\Gamma)) \cong
\Tor_{p-1}^{\cc\Gamma}(N,\caln(\Gamma))$}
and the induction step follows. This finishes
the proof of Theorem
\ref{dimension of higher Tor-s vanish in the amenable case}. \qed

\begin{theorem} \label{L^2-Betti numbers and homology in the amenable case}
Let $\Gamma$ be an amenable group and $X$ a $\Gamma$-space. Then
$$b_p^{(2)}(X;\caln(\Gamma)) ~ = ~
\dim_{\caln(\Gamma)}\left(\caln(\Gamma) \otimes_{\cc\Gamma}
H_p^{\sing}(X;\cc)\right)$$
where \mbox{$H_p^{\sing}(X;\cc)$} is the
$\cc\Gamma$-module given by the singular homology of $X$
with complex coefficients.
In particular $b_p^{(2)}(X;\caln(\Gamma))$
depends only on the $\cc\Gamma$-module
\mbox{$H_p^{\sing}(X;\cc)$}.
\end{theorem}
\proof
We have to show for a $\cc\Gamma$-chain complex $C_*$
with \mbox{$C_p = 0$} for \mbox{$p < 0$}
\begin{eqnarray}
\dim\left(H_p(\caln(\Gamma) \otimes_{\cc\Gamma} C_*)\right)
& = &
\dim\left(\caln(\Gamma) \otimes_{\cc\Gamma} H_p(C_*)\right).
\label{5.3255}
\end{eqnarray}
We begin with the case where $C_*$ is projective. Then
there is a universal coefficient spectral sequence converging to
\mbox{$H_{p+q}(\caln(\Gamma)\otimes_{\cc\Gamma} C_*)$}
\cite[Theorem 5.6.4 on page 143]{Weibel(1994)}
whose $E^2$-term is
\mbox{$E^2_{p,q} ~ = ~ \Tor^{\cc\Gamma}_p(H_q(C_*),\caln(\Gamma))$}.
Now Additivity of $\dim_{\caln(\Gamma)}$
(see Theorem \ref{properties of new definition}.4)
together with Theorem
\ref{dimension of higher Tor-s vanish in the amenable case}
imply \ref{5.3255} if $C_*$ is projective.\par

Next we prove \ref{5.3255} in the case where $C_*$ is acyclic.
One reduces the claim to two-dimensional
$C_*$ and then checks this special case
using long exact $\Tor$-sequences,
Additivity (see Theorem \ref{properties of new definition}.4)
and Theorem
\ref{dimension of higher Tor-s vanish in the amenable case}.\par

In the general case one chooses a projective $\cc\Gamma$-chain complex
$P_*$ together with a $\cc\Gamma$-chain map
\mbox{$f_*: P_* \longrightarrow C_*$}
which induces an isomorphism on homology. Since
the mapping cylinder is $\cc\Gamma$-chain homotopy
equivalent to $C_*$, the mapping cone
of $f_*$ is acyclic and hence \ref{5.3255} is true for $P_*$
and the mapping cone, we get \ref{5.3255} for $C_*$ from
Additivity (see Theorem \ref{properties of new definition}.4).
This finishes the proof of Theorem
\ref{L^2-Betti numbers and homology in the amenable case}.\qed\par

We obtain as an immediate corollary from
Theorem \ref{L^2-invariants in dimension 0}
and Theorem
\ref{L^2-Betti numbers and homology in the amenable case} (cf.
\cite[Theorem 0.2 on page 191]{Cheeger-Gromov (1986)}).

\begin{corollary}
\label{vanishing of L^2-Betti numbers for amenable groups}
If $\Gamma$ is infinite amenable, then for all \mbox{$p \ge 0$}
$$b_p^{(2)}(\Gamma) ~ = ~ 0. \qed$$
\end{corollary}

\begin{remark} \label{characterization of amenbale groups by dim(Tor)}
\em It is likely that Theorem
\ref{dimension of higher Tor-s vanish in the amenable case}
characterizes amenable groups. Namely, if $\Gamma$
contains a free group $F$ of rank two, then $\Gamma$
is non-amenable and
Theorem \ref{dimension of higher Tor-s vanish in the amenable case}
becomes false beause Theorem
\ref{exactness and dimension under induction}
implies
\begin{eqnarray*}
\dim_{\caln(\Gamma)}\left(
\Tor^{\cc\Gamma}_1(\cc[\Gamma/F];\caln(\Gamma))\right)
& = &
\dim_{\caln(\Gamma)}\left(\Tor^{\cc F}_1(\cc;\caln(\Gamma))\right))
\\
& = &
\dim_{\caln(\Gamma)}\left((\caln(\Gamma) \otimes_{\caln(F)}
\Tor^{\cc F}_1(\cc;\caln(F))\right)
\\
& = &
\dim_{\caln(F)}\left(\Tor^{\cc F}_1(\cc;\caln(F))\right))
\\
& = &
b_1^{(2)}(F)
\\
& = &
-\chi(BF)
\\
& = & 1. \qed
\end{eqnarray*}
\em \end{remark}

\begin{remark} \em
\label{flatness of caln(Gamma) over CGamma}
In view of Theorem
\ref{dimension of higher Tor-s vanish in the amenable case}
the question arises when $\caln(\Gamma)$ is flat over $\cc\Gamma$.
Except for virtually cyclic groups, i.e. groups
which are finite or contain an infinite cyclic subgroup of finite index,
we know no examples of finitely presented groups $\Gamma$
such that $\caln(\Gamma)$ is flat over $\cc\Gamma$.
If $\caln(\Gamma)$ is flat over $\cc\Gamma$, then
$\caln(\Delta)$ is flat over $\cc\Delta$ for any subgroup
\mbox{$\Delta \subset \Gamma$} by Theorem
\ref{exactness and dimension under induction}.1.
Moreover, $\caln(\Gamma)$ is flat over $\cc\Gamma$ if and only if
$\caln(\Delta)$ is flat over $\cc\Delta$ for any finitely generated
subgroup \mbox{$\Delta \subset \Gamma$}.
This follows from Theorem
\ref{exactness and dimension under induction}.1 and the facts
that the functor $\Tor$ commutes in both variables
with colimits over directed systems
\cite[Proposition VI.1.3. on page 107]{Cartan-Eilenberg (1956)},
any $\cc\Gamma$-module is the colimit of its fintely generated
submodules, any finitely generated
$\cc\Gamma$-module is the colimit of a directed system of
finitely presented $\cc\Gamma$-modules and any
finitely presented $\cc\Gamma$-submodule is obtained by induction
from a finitely presented $\cc\Delta$-module for a finitely generated
subgroup $\Delta \subset \Gamma$.
It is not hard to check
that $\caln(\Gamma)$ is flat over $\cc\Gamma$, if
$\caln(\Delta)$ is flat over $\cc\Delta$ for some subgroup
\mbox{$\Delta \subset \Gamma$} of finite index and
that $\caln(\zz)$ is flat over $\cc\zz$ using the fact that $\cc G$
is semi-simple for finite $G$ and $\cc\zz$ is a principal ideal domain.
In particular $\caln(\Gamma)$ is flat over $\cc\Gamma$
if $\Gamma$ is virtually cyclic. \par

Now suppose that $\caln(\Gamma)$ is flat over $\cc\Gamma$.
If $B\Gamma$ is a $CW$-complex of finite type,
then \mbox{$b_p^{(2)}(E\Gamma;\caln(\Gamma)) = 0$} for
\mbox{$p \ge 1$} and the $p$-th
Novikov-Shubin invariant satisfies
\mbox{$\alpha_p(B\Gamma) = \infty^+$} for
\mbox{$p \ge 2$}. This implies for instance
that $\Gamma$ does not contain a
subgroup which is isomorphic to $\zz \ast \zz$
(Remark \ref{characterization of amenbale groups by dim(Tor)})
or $\zz \times \zz$
\cite[Proposition 39 on page 494]{Lott(1992a)}.
If $\Gamma$ is non-amenable and $B\Gamma$ is a finite $CW$-complex,
then $B\Gamma$ is a counterexample to the zero-in-the-spectrum-conjecture
\cite{Lott (1997)}, \cite[section 11]{Lueck (1997)}.
\em \qed
\end{remark}


\typeout{--------------------- section 6  -----------------------}

\tit{Dimension functions and $G_0(\cc\Gamma)$}
\label{Dimension functions and G_0(ccGamma)}

Let $G_0(\cc\Gamma)$ be the abelian group which has as set
of generators the isomorphism classes
of finitely generated (not necessarily projective)
$\cc\Gamma$-modules and has for
each exact sequence of finitely generated $\cc\Gamma$-modules
\mbox{$0 \longrightarrow M_0 \longrightarrow M_1 \longrightarrow
M_2 \longrightarrow 0$} the relation
\mbox{$[M_0] - [M_1] + [M_2] = 0$}. Given a finitely generated
$\cc\Gamma$-module $M$, the $\caln(\Gamma)$-module
\mbox{$\caln(\Gamma) \otimes_{\cc\Gamma} M$} is
a finitely generated $\caln(\Gamma)$-module. We have defined
\mbox{$\bfT \caln(\Gamma) \otimes_{\cc\Gamma} M$} and
\mbox{$\bfP \caln(\Gamma) \otimes_{\cc\Gamma} M$} in Definition
\ref{closure, K(M), P(M)}. Recall from Theorem
\ref{properties of new definition}.3 that
\mbox{$\bfP \caln(\Gamma) \otimes_{\cc\Gamma} M$}
is a finitely generated projective $\caln(\Gamma)$-module.
Define maps
\begin{eqnarray}
i: K_0(\cc\Gamma) \longrightarrow G_0(\cc\Gamma) & & [P] \mapsto [P];
\label{map from K_0 to G_0}
\\
k: K_0(\cc\Gamma) \longrightarrow K_0(\caln(\Gamma))
& &
[P] \mapsto [\caln(\Gamma)\otimes_{\cc\Gamma} P].
\label{map from K_0(CGamma) to K_0(N(Gamma))}
\end{eqnarray}

\begin{lemma} \label{map from G_0(CGamma) to K_0(N(Gamma)}
If $\Gamma$ is amenable, the map
$$j: G_0(\cc\Gamma) \longrightarrow K_0(\caln(\Gamma))
\hspace*{10mm} [M] \mapsto [\bfP\caln(\Gamma) \otimes_{\cc\Gamma} M]$$
is a well-defined homomorphism.
The composition \mbox{$j \circ i$} agrees with $k$ for the maps
$i$ and $k$ defined in \ref{map from K_0 to G_0}
and \ref{map from K_0(CGamma) to K_0(N(Gamma))} above.
\end{lemma}
\proof If \mbox{$0 \longrightarrow M_0 \xrightarrow{i} M_1
\xrightarrow{p}
M_2 \longrightarrow 0$} is an exact  sequence of finitely generated
$\cc\Gamma$-modules we have to check in $K_0(\caln(\Gamma))$
\begin{eqnarray*}
[\bfP \caln(\Gamma) \otimes_{\cc\Gamma} M_0]
-
[\bfP \caln(\Gamma) \otimes_{\cc\Gamma} M_1]
+
[\bfP \caln(\Gamma) \otimes_{\cc\Gamma} M_2]
& = & 0.
\end{eqnarray*}
Consider the induced sequence
\mbox{$\bfP \caln(\Gamma)
\otimes_{\cc\Gamma} M_0 \xrightarrow{\overline{i}}
\bfP \caln(\Gamma) \otimes_{\cc\Gamma} M_1
\xrightarrow{\overline{p}}
\bfP \caln(\Gamma) \otimes_{\cc\Gamma}  M_2$}.
Obviously $\overline{p}$ is surjective as $p$ is surjective.
We conclude from Theorem
\ref{properties of new definition}.1
that $\ker(\overline{i})$ and $\ker(\overline{p})$
are a finitely generated projective $\caln(\Gamma)$-modules.
Theorem \ref{properties of new definition}.4 and
Theorem \ref{dimension of higher Tor-s vanish in the amenable case}
imply
\begin{eqnarray*}
\dim_{\caln(\Gamma)}\left(\ker(\overline{i})\right) & = & 0;
\\
\dim_{\caln(\Gamma)}\left(\im(\overline{i})\right)
& = & \dim_{\caln(\Gamma)}\left(\ker(\overline{p})\right).
\end{eqnarray*}
We conclude from Theorem \ref{properties of new definition} that
\mbox{$\overline{i}: \bfP\caln(\Gamma) \otimes_{\cc\Gamma} M_0
\longrightarrow \ker(\overline{p})$} is a weak isomorphism, i.e.
its kernel is trivial and
\mbox{$\overline{\im(\overline{i})} = \ker(\overline{p})$}.
Since the functor $\nu$ of Theorem \ref{functor nu}
respects weak exactness and the Polar Decomposition Theorem applied
to a weak isomorphism has an isomorphism as unitary part,
\mbox{$\bfP\caln(\Gamma) \otimes_{\cc\Gamma} M_0$} and
\mbox{$\ker(\overline{p})$} are
isomorphic as $\caln(\Gamma)$-modules. Since
\mbox{$\ker(\overline{p})
\oplus \bfP \caln(\Gamma) \otimes_{\cc\Gamma}  M_2$}
and \mbox{$\bfP \caln(\Gamma) \otimes_{\cc\Gamma}  M_1$} are isomorphic,
Lemma \ref{map from G_0(CGamma) to K_0(N(Gamma)} follows.\qed\par

If we regard \mbox{$\hom_{\cc\Gamma}(M,\caln(\Gamma))$}
as left $\caln(\Gamma)$-module by
\mbox{$(af)(x) = f(x) \cdot a^{\ast}$} for \mbox{$a \in \caln(\Gamma)$},
\mbox{$f \in \hom_{\cc\Gamma}(M,\caln(\Gamma))$} and \mbox{$x \in M$},
we obtain isomorphisms of $\caln(\Gamma)$-modules
\begin{eqnarray*}
\hom_{\cc\Gamma}(M,\caln(\Gamma)) & \xrightarrow{\cong} &
\left(\caln(\Gamma) \otimes_{\cc\Gamma} M\right)^{\ast};
\\
\left(\bfP \caln(\Gamma) \otimes_{\cc\Gamma} M\right)^{\ast}
 & \xrightarrow{\cong} &
\left(\caln(\Gamma) \otimes_{\cc\Gamma} M\right)^{\ast}.
\end{eqnarray*}
Since for a finitely generated projective $\caln(\Gamma)$-module
$P$ its dual $P^{\ast}$ is isomorphic to $P$, we conclude
for a finitely generated $\cc\Gamma$-module $M$
\begin{eqnarray}
j([M])
& = &
\left[\hom_{\cc\Gamma}(M,\caln(\Gamma))\right].
\label{two versions of dim^u: G_0 to class}
\end{eqnarray}

For a finitely generated projective $\caln(\Gamma)$-module $P$ let
\begin{eqnarray}
\dim_{\caln(\Gamma)}^u(P) & \in & \cent (\caln(\Gamma))
\label{center-valued von Neumann dimension}
\end{eqnarray}
be its {\em center-valued von Neumann dimension} which is given in terms
of the universal center-valued trace \mbox{$\tr_{\caln(\Gamma)}^u$}
\cite[Theorem 8.2.8 on page
517, Proposition 8.3.10 on page 525 and
Theorem 8.4.3. on page 532]{Kadison-Ringrose (1986)},
\cite[section 3]{Lueck (1995)}. The center-valued von Neumann dimension
is additive under direct sums and two finitely generated projective
$\caln(\Gamma)$-modules $P$ and $Q$ are ismorphic if and only if
\mbox{$\dim_{\caln(\Gamma)}^u(P) = \dim_{\caln(\Gamma)}^u(Q)$}.
We obtain an injection
\begin{eqnarray}
\dim_{\caln(\Gamma)}^u: ~ K_0(\caln(\Gamma)) ~ \rightarrow ~
\cent(\caln(\Gamma))^+
= \{a \in \cent(\caln(\Gamma)) \mid a = bb^{\ast} \mbox{ for }
b \in \caln(\Gamma))\}, & &
\label{map from K_0 to center}
\end{eqnarray}
which is an isomorphism if $\caln(\Gamma)$ is of type II, for instance if
$\Gamma$ is finitely generated and does not
contain an abelian subgroup of finite index
(\cite[Corollary 3.2 and Lemma 3.3]{Lueck (1995)}).

Next we investigate the relationship between $K_0(\cc\Gamma)$
and $G_0(\cc\Gamma)$ and between \mbox{$\dim_{\caln(\Gamma)}^u$}
and the Hattori-Stallings rank. Let $\con(\Gamma)$ be the set of conjugacy
classes of elements in $\Gamma$. Let
$\con(\Gamma)_f$ be the subset of $\con(\Gamma)$ of conjugacy
classes $(\gamma)$ for which each representative $\gamma$ has finite order.
Let $\con(\Gamma)_{cf}$ be the subset of $\con(\Gamma)$ of
conjugacy classes $(\gamma)$ which contain only finitely many
elements. We denote by $\class_0(\Gamma)$ and $\class_0(\Gamma)_f$
respectively the complex vector space with the set $\con(\Gamma)$ and
$\con(\Gamma)_f$ respectively as basis. We denote by
$\class(\Gamma)$ and $\class(\Gamma)_{cf}$
respectively the complex vector space
of functions from the set $\con(\Gamma)$ and $\con(\Gamma)_f$ respectively
to $\cc$. Notice that $\class_0(\Gamma)$ is the
complex vector space of class functions from $\Gamma$ to $\cc$ with
finite support. Define the
{\em universal $\cc\Gamma$-trace} of
\mbox{$\sum_{\gamma \in \Gamma} \lambda_{\gamma} \gamma \in \cc\Gamma$}
by
\begin{eqnarray}
\tr^u_{\cc\Gamma}
\left(\sum_{\gamma \in \Gamma} \lambda_{\gamma} \gamma\right)
& := & \sum_{\gamma \in \Gamma} \lambda_{\gamma} \cdot (\gamma)
\hspace{10mm} \in \class_0(\Gamma),
\label{universal cGamma trace for elements in CGamma}
\end{eqnarray}
This extends to square matrices in the usual way
\begin{eqnarray}
\tr_{\cc\Gamma}^u: M(n,n,\cc\Gamma) \longrightarrow \class_0(\Gamma)
& & A \mapsto \sum_{i=1}^n \tr_{\cc\Gamma}^u(a_{i,i}).
\label{universal cGamma trace for matrices over CGamma}
\end{eqnarray}
Let $P$ be a finitely generated projective $\cc\Gamma$-module.
Define its {Hattori-Stallings rank} by
\begin{eqnarray}
\HS(P) & := & \tr_{\cc\Gamma}^u(A) \hspace{10mm} \in \class_0(\Gamma),
\label{Hattori-Stallings rank of a module}
\end{eqnarray}
where $A$ is any element in $M(n,n,\cc\Gamma)$ with \mbox{$A^2 = A$}
such that the image of
the map \mbox{$\cc\Gamma^n \longrightarrow \cc\Gamma^n$}
given by right multiplication with $A$
is $\cc\Gamma$-isomorphic to $P$. This definition is independent of the
choice of $A$. The Hattori-Stallings rank defines a homomorphism
\begin{eqnarray}
\HS: K_0(\cc\Gamma) \longrightarrow \class_0(\Gamma) & &
[P] \mapsto \HS(P).
\label{Hattori-Stallings map}
\end{eqnarray}
Define a homomorphism
\begin{eqnarray}
\phi: \cent(\caln(\Gamma)) & \longrightarrow & \class(\Gamma)_{cf}
\label{map from center to class(Gamma_cf)}
\end{eqnarray}
by assigning to \mbox{$u \in \cent(\caln(\Gamma)$}
$$\phi(u): \con(\Gamma)_{cf} \longrightarrow \cc \hspace{10mm}
(\delta) \mapsto \tr_{\caln(\Gamma)}\left(u \cdot \sum_{\delta' \in (\delta)}
(\delta')^{-1}\right).$$

\begin{theorem} \label{K_0 and G_0 and HS and dim^u}
Suppose that $\Gamma$ is amenable.
Then the following diagram commutes
$$\begin{CD}
K_0(\cc\Gamma) & & @> \HS >> & & \class_0(\Gamma)
\\
@V i VV & & & & @V r VV
\\
G_0(\cc\Gamma) @> j >> K_0(\caln(\Gamma))
@> \dim_{\caln(\Gamma)}^u >> \cent(\caln(\Gamma))
@> \phi >> \class(\Gamma)_{cf}
\end{CD}$$
where $r$ is given by restriction
and the other maps have been defined in
\ref{map from K_0 to G_0},
Lemma \ref{map from G_0(CGamma) to K_0(N(Gamma)},
\ref{map from K_0 to center},
\ref{Hattori-Stallings map} and
\ref{map from center to class(Gamma_cf)}.
\end{theorem}
\proof
One has to show for an element \mbox{$A \in M(n,n,\cc\Gamma)$}
and \mbox{$\delta \in \Gamma$}  such that $(\delta)$ is finite
\begin{eqnarray}
\tr_{\cc\Gamma}^u(A)(\delta) & = &
(\phi \circ \tr_{\caln(\Gamma)}^u(A))(\delta). \label{eqn 6.200}
\end{eqnarray}
It suffices to show for \mbox{$\gamma \in \Gamma$}
and \mbox{$\delta \in \Gamma$} such that $(\delta)$ is finite
\begin{eqnarray}
\tr_{\cc\Gamma}^u(\gamma)(\delta) & = &
(\phi \circ \tr_{\caln(\Gamma)}^u(\gamma))(\delta). \label{eqn 6.201}
\end{eqnarray}
The universal center-valued von Neumann trace satisfies
for \mbox{$\gamma \in \Gamma$}
$$\tr_{\caln(\Gamma)}^u(\gamma) ~ = ~
\left\{\begin{array}{ll}
|(\gamma)|^{-1} \cdot \sum_{\gamma' \in (\gamma)} \gamma'
& \mbox { if } (\gamma) \mbox{ is finite }
\\
0 & \mbox{ otherwise }
\end{array}\right..$$
This follows from the facts  that
$\tr_{\caln(\Gamma)}^u$
is ultraweakly continuous and the
identity on the center of $\caln(\Gamma)$ and that
\mbox{$\tr_{\caln(\Gamma)}^u(\frac{1}{n} \cdot \sum_{i=1}^n \delta_n) =
\tr_{\caln(\Gamma)}^u(\delta)$} holds for
elements $\delta_1$, $\delta_2$, $\ldots $, $\delta_n$ in $(\delta)$.
Notice that \mbox{$\tr_{\caln(\Gamma)}(\delta)$}
is $1$ if $\delta = 1$ and $0$ otherwise and that
\mbox{$\tr_{\cc\Gamma}^u(\gamma)(\delta) = 0$} if $(\delta)$ is finite
and $(\gamma)$ is infinite.
Hence \ref{eqn 6.201} and thus \ref{eqn 6.200} follow
from the computation for \mbox{$\gamma , \delta \in \Gamma$}
such that $(\gamma)$ and $(\delta)$ are finite
\begin{eqnarray*}
\phi(|(\gamma)|^{-1} \cdot \sum_{\gamma' \in (\gamma)} \gamma')(\delta)
& = &
\tr_{\caln(\Gamma)}\left(
\left(|(\gamma)|^{-1} \cdot \sum_{\gamma' \in (\gamma)} \gamma'\right)
\cdot \left(\sum_{\delta' \in (\delta)} (\delta')^{-1}\right)\right)
\\
& = &
\sum_{\gamma' \in (\gamma)}\sum_{\delta' \in (\delta)}
|(\gamma)|^{-1} \cdot \tr_{\caln(\Gamma)}(\gamma' \cdot (\delta')^{-1})
\\
& = &
\sum_{\gamma' \in (\gamma),\delta' \in (\delta), \gamma' = \delta'}
|(\gamma)|^{-1}
\\
& = & \left\{\begin{array}{ll}
1 & \mbox{ if } (\gamma) = (\delta)
\\
0 & \mbox{ otherwise }
\end{array}\right.
\\
\\
 & = & \tr_{\cc\Gamma}^u(\gamma)(\delta).
\end{eqnarray*}
This finishes the proof of Theorem
\ref{K_0 and G_0 and HS and dim^u}.
\qed\par

\begin{lemma} \label{K_0(CG) and Hattori-Stallings rank}
Let $\Gamma$ be a discrete group. Then there is a commutative
diagram whose the left vertical arrow is an isomorphism
\begin{center}
$\begin{CD}
\left(\colim_{\Or(\Gamma,\calfin)} K_0(\cc H)\right) \otimes_{\zz}\cc 
@>l>> K_0(\cc\Gamma)\otimes_{\zz}\cc\\
@V h V \cong V @VV \HS V \\
\class_0(\Gamma)_f  @>>e > \class_0(\Gamma)
\end{CD}$
\end{center}
\end{lemma}
\proof
Firstly we explain the maps in the square.
The colimit is taken for the covariant functor
$$\Or(G,\calfin) \longrightarrow \ABEL \hspace{5mm}
G/H \mapsto K_0(\cc H)$$
to the category of abelian groups
which is given by induction. Here $\Or(G,\calfin)$
is the full subcategory of the orbit category $\Or(G)$
consisting of objects $G/H$ with finite $H$.
The map $l$ is induced by the universal
property of the colimit and the various maps
\mbox{$K_0(\cc H) \longrightarrow K_0(\cc\Gamma)$} induced
by the inclusions of finite subgroups $H$ of $\Gamma$ in $\Gamma$.
The map $e$ is given by the inclusion
\mbox{$\con (\Gamma)_f \longrightarrow \con (\Gamma)$}.\par

Define for a group homomorphism
\mbox{$\psi: \Gamma \longrightarrow \Gamma'$}
a map \mbox{$\psi_*: \con(\Gamma)  \longrightarrow  \con(\Gamma')$}
by sending $(h)$ to $\psi(h))$. It induces a homomorphism
\mbox{$\psi_*: \class_0(\Gamma) \longrightarrow \class_0(\Gamma')$}.
One easily checks that the following diagram commutes
\begin{eqnarray}
\begin{CD}
K_0(\cc\Gamma) @> \psi_* >> K_0(\cc\Gamma')
\\
@V \HS VV @V \HS VV
\\
\class_0(\Gamma)  @> \psi_* >> \class_0(\Gamma')
\end{CD} & &
\label{eqn 22.3}
\end{eqnarray}
There is a canonical isomorphism
\begin{eqnarray}
f_1: \left(\colim_{\Or(\Gamma,\calfin)}K_0(\cc H)\right)\otimes_{\zz}\cc
& \xrightarrow{\cong} &
\colim_{\Or(\Gamma,\calfin)}K_0(\cc H)\otimes_{\zz}\cc.
\label{eqn 22.4}
\end{eqnarray}
The Hattori-Stallings ranks for the various finite subgroups $H$ of $\Gamma$
induce an isomorphism
\begin{eqnarray}
f_2: \colim_{\Or(\Gamma,\calfin)}K_0(\cc H)\otimes_{\zz}\cc
& \xrightarrow{\cong} &
\colim_{\Or(\Gamma,\calfin)} \class(H).
\label{eqn 22.5}
\end{eqnarray}
Let \mbox{$f_3': \colim_{\Or(\Gamma,\calfin)} \con(H)
\longrightarrow \con(\Gamma)_f$} be the map induced by the
the inclusions of the finite subgroups $H$ of $\Gamma$.
Define a  map
\mbox{$f_4': \con(\Gamma)_f  \longrightarrow
\colim_{\Or(\Gamma,\calfin)} \con(H)$}
by sending \mbox{$(\gamma) \in \con(\Gamma)_f$} to the image of
\mbox{$(\gamma) \in \con(\langle \gamma \rangle)$}
under the canonical structure map from
\mbox{$\con(\langle \gamma \rangle )$} to
\mbox{$\colim_{\Or(\Gamma,\calfin)} \con(H)$}, where
\mbox{$\langle \gamma \rangle$} is the finite cyclic subgroup
generated by $\gamma$. One easily checks that this is independent
of the choice of the representative $\gamma$ in $(\gamma)$
and that $f_3'$ and $f_4'$ are inverse to one another.
The bijection $f_3'$ induces an isomorphism
\begin{eqnarray}
f_3: \colim_{\Or(\Gamma,\calfin)} \class(H) & \longrightarrow &
\class(\Gamma)_f,
\label{eqn 22.8}
\end{eqnarray}
because colimit and the functor sending a set to the complex vector space
with this set as basis commute. Now the isomorphism
$h$ is defined as the composition of
the isomorphisms $f_1$ from \ref{eqn 22.4}, $f_2$ from \ref{eqn 22.5}
and $f_3$ from \ref{eqn 22.8}.
It remains to check that the square in
Lemma \ref{K_0(CG) and Hattori-Stallings rank} commutes.
This follows from the commutativity of
\ref{eqn 22.3}. This finishes the proof of Lemma
\ref{K_0(CG) and Hattori-Stallings rank}. \qed

\begin{corollary} \label{image of G_0 to class in the amenable case}
Suppose that $\Gamma$ is amenable. Then the image of the
composition
$$G_0(\cc\Gamma)\otimes_{\zz}\cc \xrightarrow{j\otimes_{\zz} \cc}
 K_0(\caln(\Gamma))\otimes_{\zz} \cc \xrightarrow{\dim_{\caln(\Gamma)}^u}
\cent(\caln(\Gamma)) \xrightarrow{\phi}
\class(\Gamma)_{cf}$$
contains the complex vector space \mbox{$\class_0(\Gamma)_{f,cf}$}
with \mbox{$\con(\Gamma)_{f} \cap \con(\Gamma)_{cf}$} as basis.
\qed \end{corollary}

\begin{remark} \em
\label{conjecture about K_0(CGamma) and finite subgroups}
There is the conjecture that the canonical map
$$\colim_{\Or(\Gamma,\calfin)} K_0(\cc H) \xrightarrow{\cong}
K_0(\cc\Gamma)$$
is bijective for all groups $\Gamma$. In particular this would imply
by Lemma \ref{K_0(CG) and Hattori-Stallings rank}
that the Hattori-Stallings rank induces an isomorphism
\begin{eqnarray*}
\HS: K_0(\cc\Gamma)\otimes_{\zz} \cc &\xrightarrow{\cong} &
\class_0(\Gamma)_f. \qed
\end{eqnarray*}
\em
\end{remark}

\begin{theorem} \label{Conjecture about K_0(C Gamma)}
\begin{enumerate}
\item The map
$$l: \left(\colim_{\Or(\Gamma,\calfin)}K_0(\cc H)\right)\otimes_{\zz}\cc
\longrightarrow K_0(\cc\Gamma)\otimes_{\zz}\cc$$
is injective

\item If $\Gamma$ is virtually polycyclic, then
we obtain isomorphisms
\begin{eqnarray*}
\HS: K_0(\cc\Gamma)\otimes_{\zz} \cc & \xrightarrow{\cong} &
\class_0(\Gamma)_f;
\\
i: K_0(\cc\Gamma) & \xrightarrow{\cong} & G_0(\cc\Gamma).
\end{eqnarray*}
\end{enumerate}
\end{theorem}
\proof 1.) follows directly from Lemma
\ref{K_0(CG) and Hattori-Stallings rank}.\par\noindent
2.) Moody has shown
\cite{Moody(1989)} that the obvious map
\mbox{$\oplus_{H \in \calfin} G_0(\cc H) \longrightarrow G_0(\cc\Gamma)$}
given by induction is surjective.
Since $\Gamma$ is polycyclic the complex group ring
$\cc\Gamma$ is regular, i.e. noetherian and any $\cc\Gamma$-module
has a finite-dimensional projective resolution. Now 1.) and  Lemma
\ref{K_0(CG) and Hattori-Stallings rank} prove the claim.\qed\par

Theorem \ref{Conjecture about K_0(C Gamma)}.2 has already been proven in
\cite{Brown-Lorenz (1992)}.

\begin{remark} \label{remark on amenbale groups and [CGamma] in G_0}
\em In particular we get from Theorem \ref{K_0 and G_0 and HS and dim^u}
that the map
$$\iota: \zz \longrightarrow G_0(\cc\Gamma)
\hspace{10mm} n \mapsto [\cc\Gamma^n]$$
is split injective, provided that $\Gamma$ is amenable. It is likely
that this property characterizes amenable groups. At least
we can show for a group $\Gamma$ which contains the free group
$F_2$ in two letters as subgroup, that $\iota$ is trivial by the following
argument.\par

Induction with the inclusion \mbox{$F_2 \longrightarrow \Gamma$}
induces a homomorphism
\mbox{$G_0(\cc F_2) \longrightarrow G_0(\cc\Gamma)$}
which sends $[\cc F_2]$ to $[\cc\Gamma]$. Hence it suffices
to show \mbox{$[\cc F_2] = 0$} in \mbox{$G_0(\cc F_2)$}.
The cellular chain complex of the universal covering
of $S^1\vee S^1$ yields an exact sequence of $\cc\Gamma$-modules
\mbox{$0 \rightarrow (\cc F_2)^2  \rightarrow \cc F_2 \rightarrow \cc
\rightarrow 0$}, where $\cc$ is equipped with the trivial $F_2$-action.
Hence it suffices to show in \mbox{$[\cc] = 0$} in \mbox{$G_0(\cc F_2)$}.
Choose an epimorphism \mbox{$f: F_2 \longrightarrow \zz$}.
Restriction with $f$ defines a homomorphism
\mbox{$G_0(\cc\zz) \longrightarrow G_0(\cc F_2)$}. It
sends $\cc$ viewed as trivial $\cc\zz$-module to $\cc$ viewed as
trivial $\cc F_2$-module. Hence it remains to show
\mbox{$[\cc] = 0$} in \mbox{$G_0(\cc\zz)$}.
This follows from the exact sequence
\mbox{$0 \longrightarrow \cc\zz \xrightarrow{s-1} \cc\zz
\longrightarrow \cc \longrightarrow 0$} for $s$ a generator of $\zz$.
\em \qed\end{remark}


\typeout{--------------------- section 7 -----------------------}

\tit{Groups with vanishing $L^2$-Betti numbers}
\label{Groups with vanishing L^2-Betti numbers}

In this section we investigate the following class of groups
\begin{definition}
\label{class bnull of groups}
Define the class of groups
\begin{eqnarray*}
\bnull_{d}  & := &
\{\Gamma \mid b_p^{(2)}(\Gamma) = 0
\mbox{ for } 0 \le p \le d\};
\\
\bnull_{\infty}  & := & \{\Gamma \mid b_p^{(2)}(\Gamma) = 0
\mbox{ for } 0 \le p\}.\qed
\end{eqnarray*}
\end{definition}

Notice that $\bnull_0$ is the class of infinite groups
by Theorem \ref{L^2-invariants in dimension 0}.
Definition \ref{class bnull of groups}
is motivated among other things by
Corollary \ref{vanishing of chi_{virt}}
and the following result.

\begin{theorem} \label{consequence from Gamma in bnull}
Let $1 \longrightarrow \Delta \longrightarrow \Gamma \longrightarrow
\pi$ be an exact sequence of groups. Suppose that $\Gamma$ is
finitely presented  and one of the following conditions is satisfied.
\begin{enumerate}
\item
$|\Delta| = \infty$, $b_1^{(2)}(\Delta) < \infty$ and
$\pi$ contains an element of infinite order or
contains finite subgroups of arbitrary large order;

\item
The ordinary first Betti number of $\Delta$ satisfies
$b_1(\Delta)< \infty$  and $\pi$ belongs to $\bnull_1$.
\end{enumerate}

Then:
\begin{enumerate}

\item Let $M$ be a closed oriented $4$-manifold with $\Gamma$ as
fundamental group. Then
$$ |\sign (M)| ~ \le ~ \ch(M);$$

\item Let $\defi(\Gamma)$ be the deficiency,
i.e. the maximum $g(P)-r(P)$ for all presentations $P$
where $g(P)$ is the number of generators and $r(P)$ the number
of relations. Then
$$\defi(\Gamma) \le 1;$$

\end{enumerate}
\end{theorem}

\proof
If the first condition is satisfied, then $\Gamma$ belongs to
$\bnull_1$ by Theorem \ref{properties of bnull}.5.
Now apply
\cite[Theorem 5.1 and Theorem 6.1 on page 212]{Lueck (1994b)} .\par

Suppose that the second condition is satisfied.
Let \mbox{$p: \overline{M} \longrightarrow M$}
be the regular covering associated to
$\Delta$. There is a universal coefficient
spectral sequence converging to
\mbox{$H_{p+q}^{\pi}(\overline{M};\caln(\pi))$}
with
\mbox{$E^2_{p,q} ~ = ~
\Tor^{\cc\pi}_p(H_q(\overline{M};\cc),\caln(\pi))$}
\cite[Theorem 5.6.4 on page 143]{Weibel(1994)}.
Since $H_q(\overline{M};\cc)$ is $\cc$ with the trivial
$\pi$-action for $q=0$ and
finite-dimensional as complex vector space
by assumption for $q = 1$,
Theorem \ref{properties of new definition}.4
and Lemma \ref{invariants of caln(Gamma) otimes_{zzGamma} zz}.3
imply \mbox{$\dim(E^2_{p,q}) = 0$} for \mbox{$p + q = 1$}
and hence \mbox{$b_1^{(2)}(\overline{M};\caln(\pi)) = 0$}.
The arguments in
\cite[Theorem 5.1 and Theorem 6.1 on page 212]{Lueck (1994b)}
for the universal covering of $M$ apply also to
$\overline{M}$.
\qed

The idea to take another covering than the universal covering
in the proof of Theorem \ref{consequence from Gamma in bnull}
is taken from \cite[Corollary 5.2 on page 391]{Eckmann(1992)}.
More information about results like
Theorem \ref{consequence from Gamma in bnull} are given
can be found in \cite{Eckmann(1997)}.
\begin{theorem}
\label{properties of bnull} Let $d$ be a non-negative integer
or \mbox{$d = \infty$}. Then:
\begin{enumerate}
\item
The class $\bnull_{\infty}$
contains all infinite amenable groups;

\item If $\Gamma$ contains a normal subgroup $\Delta$
with \mbox{$\Delta \in \bnull_d$}, then
\mbox{$\Gamma \in \bnull_d$};

\item  If \mbox{$\Gamma$} is the union of a directed system of
subgroups \mbox{$\{\Gamma_i \mid i \in I\}$} such that each
$\Gamma_i$ belongs to $\bnull_d$, then \mbox{$\Gamma \in \bnull_d$};

\item Let \textshortexactsequenceone{\Delta}{i}{\Gamma}{p}{\pi}
be an exact sequence of groups such that \mbox{$b_p^{(2)}(\Delta)$}
is finite for all \mbox{$p \le d$}. Suppose that $B\pi$
has finite $d$-skeleton  and that there is an injective endomorphism
\mbox{$j: \pi \longrightarrow \pi$} whose image has finite index, but is
not equal to $\pi$ (for example $\pi = \zz^n$). Then
\mbox{$\Gamma \in \bnull_{d}$};

\item  Let \textshortexactsequenceone{\Delta}{i}{\Gamma}{p}{\pi}
be an exact sequence of groups such that
$\Delta \in \bnull_{d-1}$, $b_d^{(2)}(\Delta) < \infty$ and
$\Gamma$ contains an element of infinite order or a finite subgroup
of arbitrary large order. Then $\Gamma \in \calb_d$;

\item Suppose that there are groups
$\Gamma_1$ and $\Gamma_2$ and group homomorphisms
\mbox{$\phi_i: \Gamma_0 \longrightarrow \Gamma_i$} for
\mbox{$i=1,2$} such that $\phi_0$ and $\phi_1$ are injective, $\Gamma_0$
belongs to $\bnull_{d-1}$,$\Gamma_1$ and $\Gamma_2$
belong to $\bnull_d$ and
$\Gamma$ is the amalgamated product
\mbox{$\Gamma_1 \ast_{\Gamma_0} \Gamma_2$}
with respect to $\phi_1$ and $\phi_2$. Then $\Gamma$ belongs
to $\bnull_d$.

\end{enumerate}
\end{theorem}
\proof
1.) see Corollary
\ref{vanishing of L^2-Betti numbers for amenable groups}.
\par\noindent
2.) We obtain a fibration
\mbox{$B\Delta \longrightarrow B\Gamma \longrightarrow B\pi$}
for \mbox{$\pi = \Gamma/\Delta$}. There is the
Leray-Serre spectral sequence
converging to
\mbox{$H_{p+q}^{\Gamma}(E\Gamma,\caln(\Gamma))$}
with
\mbox{$E_{p,q}^1 =
H_q^{\Delta}(E\Delta;\caln(\Gamma))\otimes_{\zz\pi} C_p(E\pi)$}
for an appropriate $\zz\pi$-action on
\mbox{$H_q^{\Delta}(E\Delta;\caln(\Gamma))$} coming from
the fiber transport.  Because of Additivity
(see Theorem \ref{properties of new definition}.4)
it suffices to show for \mbox{$p + q \le d$}
\begin{eqnarray}
\dim_{\caln(\Gamma)}\left(E_{p,q}^1\right) & = & 0.
\label{eqn 7.1}
\end{eqnarray}
Since \mbox{$b_q^{(2)}(\Delta) =
\dim_{\caln(\Gamma)}\left(H_q^{\Delta}(E\Delta;\caln(\Gamma))\right)$}
by Theorem \ref{exactness and dimension under induction}
and \mbox{$C_p(E\pi)$} is a direct sum of copies of $\zz\pi$,
Cofinality (see Theorem \ref{properties of new definition}.4)
proves \ref{eqn 7.1}.
\par\noindent
3.) Using for instance the bar-resolution model for $E\Gamma$,
one gets that $E\Gamma$ is the colimit of a directed system
of subspaces of the form \mbox{$E\Gamma_i \times_{\Gamma_i} \Gamma$}
directed by $I$. Hence
$$H_p^{\Gamma}(E\Gamma;\caln(\Gamma)) ~ = ~
\colim_{i \in I}
H_p^{\Gamma}(E\Gamma_i \times_{\Gamma_i} \Gamma;\caln(\Gamma)).$$
Since
\mbox{$\dim_{\caln(\Gamma)}\left(
H_p^{\Gamma}(E\Gamma_i \times_{\Gamma_i} \Gamma;\caln(\Gamma))\right) =
b_p^{(2)}(\Gamma_i)$} by Theorem
\ref{induction for spaces}
the claim follows from Theorem
\ref{dimension of colimits for arbitrary index sets}.
\par \noindent
4.) Fix an integer \mbox{$n \ge 1$}. Put
\mbox{$\Gamma' = p^{-1}(\im (j^n))$}. If $k$
is the index of $\im(j)$ in $\pi$, then
$k^n$ is the index of $\im(j^n)$ in $\pi$ and of $\Gamma'$ in $\Gamma$
and we conclude
\begin{eqnarray}
b_p^{(2)}(\Gamma) & = & \frac{b_p^{(2)}(\Gamma')}{k^n}.
\label{eqn 7.2}
\end{eqnarray}
Since $\im(j^n)$ is isomorphic to $\pi$,
we have an exact sequence
\textshortexactsequenceone{\Delta}{}{\Gamma'}{}{\pi}.
Let $i_p$ be the number of $p$-cells in $B\pi$.
We get from the
Leray-Serre spectral sequence and Additivity
(see Theorem \ref{properties of new definition})
\begin{eqnarray}
b_p^{(2)}(\Gamma') & \le &
\sum_{n = 0}^p b_q^{(2)}(\Delta) \cdot i_{p-q}.
\label{eqn 7.3}
\end{eqnarray}
Equations \ref{eqn 7.2} and \ref{eqn 7.3} imply
\begin{eqnarray}
b_p^{(2)}(\Gamma) & = &
\frac{\sum_{q = 0}^p b_q^{(2)}(\Delta) \cdot i_{p-q}}{k^n}.
\label{eqn 7.4}
\end{eqnarray}
Since \mbox{$k > 1$} and \ref{eqn 7.4} holds for all $n \ge1$
and \mbox{$\sum_{q = 0}^p b_q^{(2)}(\Delta) \cdot i_{p-q}$} is finite
for $p \le d$ by assumption, the claim follows.
\par\noindent
5.) Using the spectral sequence which converges to
\mbox{$H_{p+q}^{\Gamma}(E\Gamma;\caln(\Gamma))$}
and has as $E^2$-term
\mbox{$E^2_{p,q} =
H_p^{\pi}(E\pi;H_q^{\Delta}(E\Delta;\caln(\Gamma)))$}
the proof of assertion 5.) is reduced to the proof of
\begin{eqnarray}
\dim_{\caln(\Gamma)}\left(
H_0^{\pi}(E\pi;H_d^{\Delta}(E\Delta;\caln(\Gamma)))\right)
& = & 0,
\label{eqn 7.7}
\end{eqnarray}
since \mbox{$\dim_{\caln(\Gamma)}\left(
H_q^{\Delta}(E\Delta;\caln(\Gamma))\right) =
b_q^{(2)}(\Delta)$} by
Theorem \ref{exactness and dimension under induction}
and hence vanishes for $q < d$ by assumption.
Let \mbox{$\pi'\subset \pi$} be a subgroup (not necessarily normal). Let
\mbox{$\Gamma' \subset \Gamma$} be the
preimage of $\pi$ under the canonical projection
\mbox{$\Gamma \longrightarrow \pi$}.
Then we obtain an exact sequence
\textshortexactsequenceone{\Delta}{}{\Gamma'}{}{\pi'}.
We have
\begin{eqnarray*}
H_0^{\pi'}(E\pi';H_d^{\Delta}(E\Delta;\caln(\Gamma')))
& = &
H_d^{\Delta}(E\Delta;\caln(\Gamma')) \otimes_{\cc[\pi']} \cc;
\\
H_0^{\pi}(E\pi;H_d^{\Delta}(E\Delta;\caln(\Gamma)))
& = &
H_d^{\Delta}(E\Delta;\caln(\Gamma)) \otimes_{\cc[\pi]} \cc;
\end{eqnarray*}
Since
\mbox{$H_d^{\Delta}(E\Delta;\caln(\Gamma))
\otimes_{\cc[\pi]} \cc$}
is a quotient of
\mbox{$H_d^{\Delta}(E\Delta;\caln(\Gamma))
\otimes_{\cc[\pi']} \cc$}
we conclude from
Additivity (see Theorem \ref{properties of new definition})
and from Theorem
\ref{exactness and dimension under induction}
\begin{eqnarray*}
\dim_{\caln(\Gamma)}\left(
H_d^{\Delta}(E\Delta;\caln(\Gamma))
\otimes_{\cc[\pi]} \cc\right)
 & \le &
\dim_{\caln(\Gamma)}\left(
H_d^{\Delta}(E\Delta;\caln(\Gamma))
\otimes_{\cc[\pi']} \cc\right);
\\
\dim_{\caln(\Gamma')}\left(
H_d^{\Delta}(E\Delta;\caln(\Gamma'))
\otimes_{\cc[\pi']} \cc\right)
& = &
\dim_{\caln(\Gamma)}\left(
H_d^{\Delta}(E\Delta;\caln(\Gamma))
\otimes_{\cc[\pi']} \cc\right).
\end{eqnarray*}
This implies
\begin{eqnarray*}
\dim_{\caln(\Gamma)}\left(
H_0^{\pi}(E\pi;H_d^{\Delta}(E\Delta;\caln(\Gamma)))\right)
& \le &
\dim_{\caln(\Gamma')}\left(
H_0^{\pi'}(E\pi';H_d^{\Delta}(E\Delta;\caln(\Gamma')))\right).
\end{eqnarray*}
Hence \ref{eqn 7.7} would follow if we can find for
each $\epsilon < 0$ a subgroup
\mbox{$\pi' \subset \pi$} satisfying
\begin{eqnarray}
\dim_{\caln(\Gamma')}\left(
H_0^{\pi'}(E\pi';H_d^{\Delta}(E\Delta;\caln(\Gamma')))\right)
& \le & \epsilon.
\label{eqn 7.32323}
\end{eqnarray}

We begin with the case where $\pi'$ is $\zz$.
From assertion 4 we conclude
\begin{eqnarray}
\dim_{\caln(\Gamma')}\left(
H_p^{\Gamma'}(E\Gamma';\caln(\Gamma'))\right)
& = & 0 \hspace{10mm} \mbox{ for } p \ge n.
\label{eqn 7.8}
\end{eqnarray}
The Leray-Serre spectral sequence associated to
\mbox{$1 \longrightarrow \Delta \longrightarrow \Gamma'
\longrightarrow \zz \longrightarrow 1$}
has an $E^2$-term which satisfies
\mbox{$E^2_{p,q} = 0$} for \mbox{$q \not= 0,1$}
since $B\zz$ has the $1$-dimensional model $S^1$.
Since it converges to
\mbox{$H_{p+q}^{\Gamma'}(E\Gamma';\caln(\Gamma'))$},
we conclude \ref{eqn 7.32323} for $\epsilon = 0$
from \ref{eqn 7.8} and
Additivity (see Theorem \ref{properties of new definition}).
Now suppose $\pi'$ is finite. Then we get
\begin{eqnarray*}
\dim_{\caln(\Gamma')}\left(
H_0^{\pi'}(E\pi';H_d^{\Delta}(E\Delta;\caln(\Gamma')))\right)
& = &
\dim_{\caln(\Gamma')}\left(
H_d^{\Gamma'}(E\Gamma';\caln(\Gamma'))\right)
\\
& = & b_d^{(2)}(\Gamma')
\\
& = & \frac{b_d^{(2)}(\Delta)}{|\pi'|}.
\end{eqnarray*}
If we can find $\pi'$ with arbitrary large $|\pi'|$
we get \ref{eqn 7.32323}.\par\noindent
6.) One easily checks using the Seifert-van Kampen Theorem, that there is
a $\Gamma$-push out
\comsquare{\Gamma \times_{\Gamma_0} E\Gamma_0}{}
{\Gamma \times_{\Gamma_1} E\Gamma_1}{}{}
{\Gamma \times_{\Gamma_2} E\Gamma_2}{}
{E\Gamma}
We conclude from Theorem
\ref{induction for spaces}
\mbox{$\dim\left(H_p^{\Gamma}(E\Gamma_i \times_{\Gamma_i} \Gamma
;\caln(\Gamma))\right) = b_p^{(2)}(\Gamma_i)$}.
Now the claim
follows from Additivity
(see Theorem \ref{properties of new definition}.4)
and the a long exact homology sequence for
\mbox{$H_{\ast}^{\Gamma}(-,\caln(\Gamma))$}.
This finishes the proof of Theorem
\ref{properties of bnull}.\qed\par

So far we have no example with negative answer
to the following question
and can give an affirmative answer in some special cases.

\begin{question} \em
\label{question bnull and extensions}
Let \textshortexactsequenceone{\Delta}{}{\Gamma}{}{\pi}
be an exact sequence such that
\mbox{$b_p^{(2)}(\Delta) < \infty$}
for all \mbox{$p \ge 0$} and $\pi$ belongs to $\bnull_{\infty}$.
Does then $\Gamma$ belong to $\bnull_{\infty}$? \par

More generally, if
\mbox{$F \longrightarrow E \longrightarrow B$} is a fibration
such that \mbox{$b_p^{(2)}(F;\caln(\pi_1(E))) < \infty $}
and \mbox{$b_p^{(2)}(B;\caln(\pi_1(B))) = 0$} holds for
\mbox{$p \ge 0$}, does then
\mbox{$b_p^{(2)}(E;\caln(\pi_1(E))) = 0$} hold for
\mbox{$p \ge 0$}?
\qed \em
\end{question}

\begin{remark}
\label{further references}\em
Compact $3$-manifolds whose fundamental groups belong to $\bnull_{\infty}$
are characterized in \cite[Proposition 6.5 on page 54]{Lott-Lueck (1995)}.
The generalized Singer-Conjecture says that for an aspherical closed
manifold $M$ all the $L^2$-Betti numbers of its universal
covering vanish possibly except in the middle dimension. In particular
it implies that the fundamental group of an aspherical closed
odd-dimensional manifold belongs to $\bnull_{\infty}$. Thompson's
group $F$ belongs to $\bnull_{\infty}$ \cite[Theorem 0.8]{Lueck (1995)}.
More information about the class $\bnull_1$ is given in
\cite{Bekka-Valette(1994)}. \qed
\em \end{remark}


\typeout{--------------------- section 8 -----------------------}

\tit{$L^2$-Euler characteristics and the Burnside group}
\label{L^2-Euler characteristics and the Burnside group}

In this section we extend some of the results
\cite{Cheeger-Gromov (1986)} about $L^2$-Euler characteristics
and investigate the Burnside group of a discrete group $\Gamma$.
This extends the classical notions of the Burnside ring,
Burnside ring congruences and equivariant
Euler characteristics for finite groups.

If $X$ is a $\Gamma$-$CW$-complex, denote by
\mbox{$I(X)$} the set of its equivariant cells.
For a cell \mbox{$c \in I(X)$} let $(\Gamma_c)$ be
the conjugacy class of subgroups of $\Gamma$ given by its orbit type
and let $\dim(c)$ be its dimension.
Denote by $|\Gamma_c|^{-1}$ the inverse of the order of any representative
of $(\Gamma_c)$, where $|\Gamma_c|^{-1}$ is to be understood to be zero
if the order is infinite.

\begin{definition}
\label{i(X,V), h(X,V) and chi^{(2)}(X;V)}
Let $X$ be a (left) $\Gamma$-space and
$V$ be a $\cala$-$\zz\Gamma$-module. Define
\begin{eqnarray*}
h(X;V) & := & \sum_{p \ge 0} b_p^{(2)}(X;V) ~  \in [0,\infty];
\\
\chi^{(2)}(X;V) & := & \sum_{p \ge 0} (-1)^p \cdot b_p^{(2)}(X;V)
~ \in \rr \hspace{10mm}\mbox{ ,if } h(X;V) < \infty;
\\m(X) & := & \sum_{c \in I(X)} |\Gamma_c|^{-1} ~  \in [0,\infty],
\hspace{10mm}\mbox{ ,if } X
\mbox{ is a } \Gamma\mbox{-}CW\mbox{-complex}. \qed
\end{eqnarray*}
\end{definition}

The condition \mbox{$h(X;V) < \infty$} ensures that the sum defining
\mbox{$\chi^{(2)}(X;V)$} converges and that
\mbox{$\chi^{(2)}(X;V)$} satisfies
the usual additivity formula, i.e. for a $\Gamma$-$CW$-complex
$X$ with $\Gamma$-$CW$-subcomplexes $X_0$, $X_1$ and $X_2$ satisfying
\mbox{$X = X_1 \cup X_2$}, \mbox{$X_0 = X_1 \cap X_2$}
and \mbox{$h(X_k;V) < \infty$} for \mbox{$k= 0,1,2$}
one has
\begin{eqnarray}
h(X;V) & < & \infty;
\label{additivity for h(X;V)}
\\
\chi^{(2)}(X;V) & = &
\chi^{(2)}(X_1;V) + \chi^{(2)}(X_2;V) - \chi^{(2)}(X_0;V).
\label{additivity for the L^2-Euler characteristic}
\end{eqnarray}
The next theorem
generalizes \cite[Theorem 0.3 on page 191]{Cheeger-Gromov (1986)}

\begin{theorem}
\label{chi^{(2)} and cells}
Let $X$ and $Y$ be $\Gamma$-$CW$-complexes such that
\mbox{$m(X) < \infty$} and \mbox{$m(Y) < \infty$} holds. Then
\begin{enumerate}
\item \begin{eqnarray*}
h(X;\caln(\Gamma)) & < & \infty;
\\
\sum_{c \in I(X)} (-1)^{\dim(c)} \cdot |\Gamma_c|^{-1} & = &
\chi^{(2)}(X;\caln(\Gamma)).
\end{eqnarray*}

\item Suppose that $\Gamma$ is amenable. Then
\begin{eqnarray*}
\sum_{c \in I(X)} (-1)^{\dim(c)} \cdot |\Gamma_c|^{-1} & = &
\sum_{p \ge 0} (-1)^p \cdot
\dim\left(\caln(\Gamma)\otimes_{\cc\Gamma} H_p(X;\cc)\right),
\end{eqnarray*}
where \mbox{$H_p(X;\cc)$} is the cellular or the singular homology of $X$
with complex coefficients. In particular
\mbox{$\sum_{c \in I(X)} (-1)^{\dim(c)} \cdot |\Gamma_c|^{-1}$}
depends only the $\cc\Gamma$-isomorphism
class of the $\cc\Gamma$-modules
\mbox{$H_n(X;\cc)$} for all \mbox{$n \ge 0$};

\item If for all \mbox{$c \in I(X)$}
the group $\Gamma_c$ is finite or belongs to the class
$\bnull_{\infty}$, then
\begin{eqnarray*}
b_p^{(2)}(X;\caln(\Gamma)) & = &
b_p^{(2)}(E\Gamma \times X;\caln(\Gamma)) \hspace{10mm}
\mbox{ for }p \ge 0;
\\
\chi^{(2)}(X;\caln(\Gamma)) & = &
\chi^{(2)}(E\Gamma \times X;\caln(\Gamma));
\\
\sum_{c \in I(X)} (-1)^{\dim(c)} \cdot |\Gamma_c|^{-1}& = &
\chi^{(2)}(E\Gamma \times X;\caln(\Gamma)).
\end{eqnarray*}

\item Suppose that
\mbox{$f: X \longrightarrow Y$} is a $\Gamma$-equivariant map,
such that the induced map \mbox{$H_p(f;\cc)$} on the
singular or cellular homology with complex coefficients is bijective.
Suppose that for all  \mbox{$c \in I(X)$}
and \mbox{$c \in I(Y)$}
the group $\Gamma_c$ is finite or belongs to the class
$\bnull_{\infty}$. Then
\begin{eqnarray*}
\sum_{c \in I(X)} (-1)^{\dim(c)} \cdot |\Gamma_c|^{-1}
& = &
\sum_{c \in I(Y)} (-1)^{\dim(c)} \cdot |\Gamma_c|^{-1}.
\end{eqnarray*}
\end{enumerate}
\end{theorem}
\proof
1.) Additivity and Cofinality
(see Theorem \ref{properties of new definition})
and Lemma \ref{invariants of caln(Gamma) otimes_{zzGamma} zz}.1  imply
\begin{eqnarray*}
\dim(C_p(X;\caln(\Gamma)))
& = &
\sum_{c \in I(X), \dim(c) = p} |\Gamma_c|^{-1};
\\
\dim(H_p^{\Gamma}(X;\caln(\Gamma)))
& \le &
\dim(C_p(X;\caln(\Gamma)));
\\
\sum_{p \ge 0} (-1)^p \cdot \dim(C_p(X;\caln(\Gamma)))
& = &
\chi^{(2)}(X;\caln(\Gamma))).
\end{eqnarray*}
2.) follows from the first assertion and Theorem
\ref{L^2-Betti numbers and homology in the amenable case}.
\par\noindent
3.) Because of the first equation it suffices to prove
that the dimension of the kernel and the cokernel of the map induced by the
projection
\begin{eqnarray*}
\pr_{\ast}:
H_p^{\Gamma}(E\Gamma \times X;\caln(\Gamma))
\longrightarrow H_p^{\Gamma}(X;\caln(\Gamma))
\end{eqnarray*}
are trivial.
Notice that $X$ is the colimit of its finite $\Gamma$-subcomplexes.
Since \mbox{$H_p^{\Gamma}(-,\caln(\Gamma)$} is compatible with colimits
and colimit preserves exact sequences, we can assume by Theorem
\ref{dimension of colimits for arbitrary index sets}
and Additivity (see Theorem \ref{properties of new definition})
that $X$ itself is finite.
By induction over the number of equivariant
cells, the long exact homology
sequence and Additivity (see Theorem \ref{properties of new definition})
the claim reduces to the case where $X$ is of the shape
$\Gamma/H$. Because of Theorem
\ref{induction for spaces}
it suffices to prove for the map
\mbox{$\pr_{\ast}: H_p^{H}(EH;\caln(H))
\longrightarrow H_p^{H}(\{\ast\};\caln(H))$}
that its kernel and cokernel
have trivial dimension, provided that $H$ is finite or
belongs to $\bnull_{\infty}$.
This is obvious for finite $H$ and follows
for $H \in \bnull_{\infty}$ from the definition of $\bnull_{\infty}$ and
Theorem
\ref{L^2-invariants in dimension 0}.
\par\noindent
4.) Since $E\Gamma \times X$ is free and the map
\mbox{$\id \times f:  E\Gamma \times X \longrightarrow
E\Gamma \times Y$} induces an isomorphism on singular homology
it induces an isomorphism
\mbox{$H_p^{\Gamma}(E\Gamma \times X;\caln(\Gamma)) \longrightarrow
H_p^{\Gamma}(E\Gamma \times Y;\caln(\Gamma))$} and we
conclude
\mbox{$\chi^{(2)}(E\Gamma \times X;\caln(\Gamma))  =
\chi^{(2)}(E\Gamma \times Y;\caln(\Gamma))$}.
Now assertion 4.) follows from assertion 3.).
This finishes the proof of
Theorem \ref{chi^{(2)} and cells}.\qed\par

As explained in \cite[Proposition 0.4 on page 192]{Cheeger-Gromov (1986)}
the $L^2$-Euler characteristic extends the notion of the
virtual Euler characteristic which is due to
Wall \cite{Wall (1961)}.  Information about this notion 
can be found  for instance in
\cite[chaper IX]{Brown (1982)}. The next result
generalizes \cite[Corollary 0.6 on page 193]{Cheeger-Gromov (1986)}

\begin{corollary} \label{vanishing of chi_{virt}}
Let $\Gamma$ be a group belonging to $\bnull_{\infty}$.
Then $\chi^2(E\Gamma;\caln(\Gamma))$ is defined and vanishes.
If its virtual Euler characteristic $\chi_{\virt}(\Gamma)$ is defined,
then it vanishes. In particular $\chi(B\Gamma)$ vanishes
if $B\Gamma$ can be choosen to be a finite $CW$-complex. \qed
\end{corollary}

Next we introduce the Burnside group and the equivariant
Euler characteristic. The elementary proof the following
lemma is left to the reader.

\begin{lemma} \label{orbits and Weyl groups}
Let $H$ and $K$ be subgroups of $\Gamma$. Then

\begin{enumerate}

\item $\Gamma/H^K ~ = ~ \{gH \mid g^{-1}Kg \subset H\};$

\item The map
$$\phi: \Gamma/H^K \longrightarrow \consub(H)
\hspace*{10mm} gH \mapsto g^{-1}Kg$$
induces an injection
$$W\!K\backslash(\Gamma/H^K) \longrightarrow  \consub(H),$$
where $\consub(H)$ is the set of conjugacy classes in $H$
of subgroups of $H$;

\item The $W\!K$-isotropy group of \mbox{$gH \in \Gamma/H^K$} is
\mbox{$(gHg^{-1}\cap N\!K)/K \subset N\!K/K =W\!K$},
where $N\!K$ is the normalizer of $K$ in $\Gamma$ and $W\!K$ is
$N\!K/K$;

\item If $H$ is finite, then $\Gamma/H^K$ is a
finite union of $W\!K$-orbits of the shape $W\!K/L$ for
finite subgroups \mbox{$L \subset W\!K$}. \qed

\end{enumerate}
\end{lemma}

\begin{definition} \label{definition of Burnside group}
Define the {\em Burnside group} $A(\Gamma)$ by the Grothendieck group
of the abelian monoid under disjoint union
of $\Gamma$-isomorphism classes of proper cocompact $\Gamma$-sets $S$,
i.e. $\Gamma$-sets $S$ for which the isotropy group of each element in $S$
and the quotient $\Gamma\backslash S$ are finite. \qed
\end{definition}

Notice that $A(\Gamma)$ is the free abelian group generated by
$\Gamma$-isomorphism classes of orbits $\Gamma/H$ for finite
subgroups $H \subset \Gamma$ and that $\Gamma/H$ and $\Gamma/K$
are $\Gamma$-isomorphic if and only if $H$ and $K$ are
conjugate in $\Gamma$. If $\Gamma$ is a finite group,
$A(\Gamma)$ is the classical Burnside ring
\cite[section 5]{Dieck(1979)},
\cite[chapter IV]{tom Dieck (1987)}. If $\Gamma$ is infinite,
then the cartesian product of two proper cocompact $\Gamma$-sets
with the diagonal action is not cocompact any more so that the cartesian
product does not induce a ring structure on $A(\Gamma)$.
At least there is a bilinear map induced by the cartesian product
\mbox{$A(\Gamma_1) \otimes A(\Gamma_2)  \longrightarrow
A(\Gamma_1 \times \Gamma_2)$}.

\begin{definition} \label{universal equivariant Euler characteristic}
Let $X$ be a proper finite $\Gamma$-$CW$-complex. Define its
{\em equivariant Euler characteristic}
$$\chi^{\Gamma}(X) ~ := ~ \sum_{c \in I(X)} (-1)^{\dim(c)}\cdot
[\Gamma/\Gamma_c]
 \hspace{10mm}~ \in ~ A(\Gamma).\qed$$
\end{definition}

An {\em additive invariant} $(A,a)$ for proper finite
$\Gamma$-$CW$-complexes $X$ consists of an abelian group $A$ and
a function $a$ which assigns to any proper finite $\Gamma$-$CW$-complex
$X$ an element \mbox{$a(X) \in A$} such that the following three conditions
hold, i.) if $X$ and $Y$ are $\Gamma$-homotopy
equivalent, then $a(X) = a(Y)$, ii.)
if $X_0$,$X_1$ and $X_2$ are $\Gamma$-$CW$-subcomplexes of $X$ with
\mbox{$X = X_1 \cup X_2$} and \mbox{$X_0 = X_1 \cap X_2$}, then
\mbox{$a(X) = a(X_1) + a(X_2) - a(X_0)$}, and iii.) \mbox{$a(\emptyset) =
0$}.
We call an additive invariant $(U,u)$ {\em universal},
 if for any additive invariant
$(A,a)$ there is precisely one homomorphism \mbox{$\psi: U \longrightarrow
A$}
such that \mbox{$\psi(u(X)) = a(X))$} holds for all proper finite
$\Gamma$-$CW$-complexes. One easily checks using induction over the number
of equivariant cells

\begin{lemma} \label{(A(Gamma),chi^{Gamma}) is universal additive invariant}
\mbox{$(A(\Gamma),\chi^{\Gamma})$}
is the universal additive invariant for finite proper
$\Gamma$-$CW$-complexes.
\qed
\end{lemma}

\begin{definition} \label{L^2-character map for the Burnside group}
Define for a finite subgroup \mbox{$K \subset \Gamma$}
the {\em  $L^2$-character map}
$$\ch_K^{\Gamma}: A(\Gamma) \longrightarrow \qq
\hspace{10mm} [S] \mapsto \sum_{i=1}^r |L_i|^{-1}$$
if $W\!K/L_1$, $W\!K/L_2$,\ldots $W\!K/L_r$ are the $W\!K$-orbits of $S^K$.
Define the {\em global $L^2$-character map} by
$$\ch^{\Gamma} := \prod_{(K)} \ch_K^{\Gamma}: ~
A(\Gamma) \longrightarrow \prod_{(K)} \qq$$
where $(K)$ runs over the conjugacy classes of finite subgroups
of $\Gamma$.  \qed
\end{definition}

\begin{lemma} \label{chi_K^{Gamma}(X) for finite proper G-CW-complex X}
Let $X$ be a finite proper $\Gamma$-$CW$-complex
and \mbox{$K \subset \Gamma$}
be a finite subgroup. Then $X^K$ is a finite proper $W\!K$-$CW$-complex
and
$$\chi^{(2)}(X^K;\caln(W\!K)) ~ = ~ \ch_K^{\Gamma}(\chi^{\Gamma}(X)).$$
\end{lemma}
\proof
The $W\!K$-space $X^K$ is a finite proper $W\!K$-$CW$-complex
because for finite \mbox{$H \subset \Gamma$} the
$W\!K$-set $\Gamma/H^K$ is proper and cocompact by Lemma
\ref{orbits and Weyl groups}. Since
the assignment which associates to a finite proper $\Gamma$-$CW$-complex
$X$ the element \mbox{$\chi^{(2)}(X^K;\caln(W\!K))$} in $\qq$
is an additive invariant, it suffices by Lemma
\ref{(A(Gamma),chi^{Gamma}) is universal additive invariant}
to check the claim for $X= \Gamma/H$ for finite \mbox{$H \subset \Gamma$}.
Then the claim follows from the fact that
\mbox{$\chi^{(2)}(W\!K/L;\caln(W\!K)) = |L|^{-1}$} holds
for finite  \mbox{$L \subset W\!K$}. \qed

Notice that one gets from Lemma
\ref{orbits and Weyl groups} the following explicit formula
for the value of $\ch^{\Gamma}_K(\Gamma/H)$. Namely,
define
\begin{eqnarray*}
{\cal L}_K(H) & := & \{(L) \in \consub(H) \mid
L \mbox{ conjugated to } K \mbox{ in }\Gamma\}.
\end{eqnarray*}
For $(L) \in {\cal L}_K(H)$ choose \mbox{$L \in (L)$}
and \mbox{$g \in \Gamma$}
with \mbox{$g^{-1}Kg = L$}. Then
\begin{eqnarray*}
g(H \cap N\!L)g^{-1} & = & gHg^{-1} \cap N\!K;
\\
|(gHg^{-1} \cap N\!K)/K|^{-1} & = & \frac{|K|}{|H \cap N\!L|}.
\end{eqnarray*}
This implies
\begin{eqnarray}
\label{value of chi^{Gamma}_K(Gamma/H))}
\ch_K^{\Gamma}(\Gamma/H) & = & \sum_{L \in {\cal L}_K(H)}
\frac{|K|}{|H \cap N\!L|}.
\label{eqn 8.59965}
\end{eqnarray}

\begin{lemma} \label{injectivity resp bijectivity of the global
character map}
The global $L^2$-character map of
Definition \ref{L^2-character map for the Burnside group}
induces a map denoted by
$$\ch^{\Gamma}  \otimes_{\zz} \qq: ~
A(\Gamma) \otimes_{\zz} \qq \longrightarrow \prod_{(K)} \qq.$$
It is injective. If $\Gamma$ has only finitely many conjugacy classes
of finite subgroups, then it is bijective.
\end{lemma}
\proof
Consider an element \mbox{$\sum_{i=1}^n r_i \cdot [\Gamma/H_i]$}
in the kernel of \mbox{$\ch^{\Gamma}\otimes_{\zz}\qq$}. We show by
induction over $n$ that the element must be trivial. The begin $n=0$ is
trivial, the induction step done as follows.We can choose the
numeration such that $H_i$ subconjugated to $H_j$ implies
\mbox{$i \ge j$}.
We get from \ref{value of chi^{Gamma}_K(Gamma/H))}
\begin{eqnarray*}
\ch_K^{\Gamma}(\Gamma/H) & = & 1 \hspace{10mm} \mbox{, if } H = K;
\\
\ch_K^{\Gamma}(\Gamma/H) & = & 0 \hspace{10mm} \mbox{, if } K
\mbox{  is not subconjugated to } H \mbox{ in } \Gamma.
\end{eqnarray*}
This implies
$$\ch^{\Gamma}_{H_1}\left(\sum_{i=1}^n r_i \cdot [\Gamma/H_i]\right) =
r_1$$
and hence \mbox{$r_1 = 0$}. Hence the global $L^2$-character map is
injective.
If $\Gamma$ has only finitely many conjugacy classes of finite
subgroups, then the source and target of
\mbox{$\ch^{\Gamma} \otimes_{\zz}\qq$}
are rational vector spaces of the same finite dimension and
\mbox{$\ch^{\Gamma}\otimes_{\zz}\qq$} must be bijective.
\qed

\begin{remark} \label{explicite inverse of the global L^2-character map} \em
Suppose that there are only finitely many
conjugacy classes $(H_1)$, $(H_2)$, $\ldots$, $(H_r)$
of finite subgroups in $\Gamma$. Without loss
of generality we can assume that $H_i$ subconjugated to $H_j$ implies
\mbox{$i \ge j$}. With respect to the obvious ordered basis
for the source and target the map
\mbox{$\ch^{\Gamma} \otimes_{\zz}\qq$} is described
by an upper triangular matrix $A$ with ones on the diagonal.
One can get an explicit inverse $A^{-1}$
which again has ones on the diagonal.  This leads
to a characterization of the image of $A(\Gamma)$ under
the global $L^2$-character map $\chi^{\Gamma}$. Namely, an element in
\mbox{$\eta \in \prod_{i=1}^r\qq$}
lies in \mbox{$\ch^{\Gamma}(A(\Gamma))$} if and only if
the following {\em Burnside integrality conditions} are satisfied
\begin{eqnarray}
A^{-1}\eta & \in & \prod_{i = 1}^r \zz.
\label{Burnside group congruences}
\end{eqnarray}

Now suppose that $\Gamma$ is finite. Then the global
$L^2$-character map is related to the classical character
map by the factor $|W\!K|^{-1}$, i.e. we have for each subgroup
$K$ of $\Gamma$ and any finite $\Gamma$-set $S$
\begin{eqnarray}
\ch^{\Gamma}_K(S) & = & |W\!K|^{-1} \cdot |S^K|.
\label{relation between L^2-character map and classical character map}
\end{eqnarray}
One easily checks that under the identification
\ref{relation between L^2-character map and classical character map}
the integrality conditions \ref{Burnside group congruences}
correspond to the classical Burnside ring congruences for finite groups
\cite[section 5.8]{Dieck(1979)},
\cite[section IV.5]{tom Dieck (1987)}.
\qed \em
\end{remark}

Let $E(\Gamma,\calfin)$ be the classifying $\Gamma$-space
for the family $\calfin$ of finite subgroups. This $\Gamma$-$CW$-complex
is characterized up to $\Gamma$-homotopy by the property that its
$H$-fixed point set is contractible if \mbox{$H \subset \Gamma$} is finite
and empty otherwise. It is also called the classifying space for proper
$\Gamma$-spaces and denoted by $\underline{E}\Gamma$ in the literature.
For more information about $E(\Gamma,\calfin)$ we refer for instance to
\cite{Baum-Connes-Higson(1994)}, \cite[section 7]{Davis-Lueck (1996)},
\cite[section I.6]{tom Dieck (1987)},

\begin{lemma}
\label{amenable groups and chi^{Gamma}_K(E(Gamma,calfin)))}
\label{ch^{Gamma}_K(chi^{Gamma}(E(Gamma,calfin)))}
Suppose that there is model for \mbox{$E(\Gamma,\calfin)$} which
is a finite $\Gamma$-$CW$-complex. Then there are only
finitely many conjugacy classes of finite subgroups and
for a finite subgroup \mbox{$K\subset \Gamma$}
$$\ch^{\Gamma}_K(\chi^{\Gamma}(E(\Gamma,\calfin))) ~ = ~
\chi^{(2)}(W\!K).$$
If $\Gamma$ is amenable, then we get for a finite subgroup
\mbox{$K\subset \Gamma$}
$$\ch^{\Gamma}_K(\chi^{\Gamma}(E(\Gamma,\calfin))) ~ = ~
|W\!K|^{-1},$$
where $|W\!K|^{-1}$ is to be understood as $0$ for infinite $W\!K$.
\end{lemma}
\proof We get from Lemma
\ref{chi_K^{Gamma}(X) for finite proper G-CW-complex X}
since $E(\Gamma,\calfin)^K$ is a model for $E(W\!K,\calfin)$
\begin{eqnarray*}
\ch^{\Gamma}_K(\chi^{\Gamma}(E(\Gamma,\calfin ))) & = &
\chi^{(2)}(E(W\!K,\calfin);\caln(W\!K)).
\end{eqnarray*}
Now apply Theorem \ref{chi^{(2)} and cells}.3.
and Corollary \ref{vanishing of L^2-Betti numbers for amenable groups}.
\qed

\begin{example} \label{example of an extension of zz^n by zz/p}
\em Let
\mbox{$1 \longrightarrow \zz^n \longrightarrow \Gamma
\longrightarrow\zz/p \longrightarrow 1$}
be an extension of groups for \mbox{$ n \ge 1$}
and a prime number $p$. Then $E(\Gamma,\calfin)$
can be choosen as a finite $\Gamma$-$CW$-complex
because only the following cases can occur.
If $\Gamma$ contains a finite subgroup, then
$\Gamma$ is a semi-direct product of $\zz^n$ and $\zz/l$
and one can construct a finite $\Gamma$-$CW$-complex as
model for $E(\Gamma,\calfin)$ with $\rr^n$ as underlying space.
If the group $\Gamma$ contains no finite subgroup,
then $\Gamma$ is an extention of finitely generated
free abelian groups and hence $B\Gamma$ can be chosen as a finite
$CW$-complex. We want to compute
\mbox{$\chi^{\Gamma}(E(\Gamma,\calfin))$}. The conjugation action of
$\Gamma$ on the normal subgroup $\zz^n$ factorizes through
the projection \mbox{$\Gamma\longrightarrow \zz/p$}
into a operation $\rho$ of $\zz/p$ onto $\zz^n$. If this operation
has a non-trivial fixed point, then $W\!H$ is infinite for any finite
subgroup $H$ of $\Gamma$ and we conclude from Lemma
\ref{injectivity resp bijectivity of the global
character map} and Theorem
\ref{amenable groups and chi^{Gamma}_K(E(Gamma,calfin)))} that
\begin{eqnarray*}
\chi^{\Gamma}(E(\Gamma,\calfin)) & = & 0.
\end{eqnarray*}
Now suppose that this operation $\rho$ has no non-trivial fixed points.
Let $H_0$ be the trivial subgroup and $H_1$, $H_2$,
$\ldots $, $H_r$ be a complete set of representatives of the
conjugacy classes of finite subgroups. Each $H_i$ is
isomorphic to $\zz/l$. One easily checks that there is a bijection
$$H^1(\zz /p;\zz^n_{\rho}) ~ \longrightarrow ~
\{(H) \mid H \subset \Gamma, 1 < |H| < \infty \},$$
if $\zz^n_{\rho}$ denotes the $\zz[\zz /l]$-module given by
$\zz^n$ and $\rho$. We compute using \ref{eqn 8.59965}
$$\begin{array}{lcllll}
\ch^{\Gamma}_{H_0}(\Gamma/H_0) & = & 1; & & &
\\
\ch^{\Gamma}_{H_0}(\Gamma/H_j) & = & \frac{1}{p} & & & j = 1,2, \ldots r;
\\
\ch^{\Gamma}_{H_i}(\Gamma/H_j) & = & 1 & & &
i = j, ~ i,j = 1,2, \ldots r;
\\
\ch^{\Gamma}_{H_i}(\Gamma/H_j) & = & 0 & & &
i \not= j, ~ i,j = 1,2, \ldots r;.
\\
\end{array}$$
We conclude
\begin{eqnarray*}
\chi^{\Gamma}(E(\Gamma,\calfin)) & = &
- \frac{r}{p}\cdot [\Gamma/H_0] + \sum_{i=1}^r ~ [\Gamma/H_i].
\end{eqnarray*}
The integrality conditions of \ref{Burnside group congruences}
become in this case
\begin{eqnarray*}
\eta_0 - \frac{1}{p} \cdot \sum_{i = 1}^r \eta_i & \in & \zz;
\\
\eta_i & \in &\zz \hspace{10mm} \mbox{i = 1, 2 ,\ldots r}. \qed
\end{eqnarray*}
\em \end{example}


\typeout{--------------------- section 9 -----------------------}

\tit{Values of $L^2$-Betti numbers}
\label{Values of L^2-Betti numbers}

In this section we investigate the possible
values of $L^2$-Betti numbers.\par

\begin{conjecture} \label{Atiyah Conjecture}
Let $\Gamma$ be a group and let $X$ be a free finite
$\Gamma$-$CW$-complex. Then
$$b_p^{(2)}(X;\caln(\Gamma)) ~ \in \qq.$$
If $d$ is a positive integer
such that the order of any finite subgroup of
$\Gamma$ divides $d$, then
$$d \cdot b_p^{(2)}(X;\caln(\Gamma))
~ \in \zz.\qed$$
\end{conjecture}

The significance of Conjecture \ref{Atiyah Conjecture}
and its relation to a question of
Atiyah \cite[page 72]{Atiyah (1976)}
about the rationality of the analytic
$L^2$-Betti numbers of \ref{analytic p- th L^2-Betti number} are
explained in \cite[Section 2]{Lueck (1997)}.
Let $\calc$ be the smallest class of groups
which contains all free groups, is closed
under directed unions and satisfies \mbox{$G \in {\cal C}$}
whenever $G$ contains a normal subgroup $H$ such that $H$ belongs
to $\calc$ and $G/H$ is elementary amenable. Conjecture
\ref{Atiyah Conjecture} has been proven for groups $\Gamma \in \calc$ 
by Linnell \cite{Linnell (1993)}, provided that 
there is an upper bound on the order of finite subgroups of $\Gamma$.\par

\begin{theorem} \label{possible values of the L^2-Betti numbers}

\begin{enumerate}

\item Let $\Gamma$ be a group such that there is no bound
on the order of finite subgroups. Then
\begin{enumerate}

\item Given $\beta \in [0,\infty]$,
there is a countably generated projective
$\zz\Gamma$-module $P$ satisfying
$$\dim_{\caln(\Gamma)}(\caln(\Gamma) \otimes_{\zz\Gamma} P) ~ = ~ \beta;$$

\item Given a sequence $\beta_3$, $\beta_4$, $\ldots$ of elements in
$[0,\infty]$, there is a free $\Gamma$-$CW$-complex $X$ satisfying
$$b_p^{(2)}(X;\caln(\Gamma)) ~ = ~ \beta_p
\hspace{10mm}  \mbox{ for }  p \ge 3;$$
If $\Gamma$ is countably presented, one can arrange that $X$ has
countably many $\Gamma$-equivariant cells;
\end{enumerate}

\item Let $\Gamma$ be a group such that there is a bound
on the order of finite subgroups. Let $d$ be the least common
multiple of the orders of finite subgroups of $\Gamma$.
Suppose that Conjecture \ref{Atiyah Conjecture} holds for
$\Gamma$. Then we get for any $\Gamma$-space
$X$ and $p \ge 0$
$$d \cdot b_p^{(2)}(X;\caln(\Gamma)) ~ \in \zz \cup \{\infty\};$$

\item Given a sequence of elements
$\beta_1$, $\beta_2$, $\ldots$ $[0,\infty]$,
there is a countable group $\Gamma$
with \mbox{$b_p^{(2)}(\Gamma) = \beta_p$}
for $p \ge 1$. If $\beta_1$ is rational, $\Gamma$ can be chosen to be
finitely generated.

\end{enumerate}
\end{theorem}
\proof
1.a.) Since there is no bound on the order of finite subgroups,
we can find a sequence of finite subgroups
$H_1$, $H_2$, $\ldots$ of $\Gamma$ such that
\mbox{$\beta = \sum_{i = 1}^{\infty} |H_i|^{-1}$}. Then
$$P ~ = ~ \oplus_{i = 1}^{\infty}  \cc[\Gamma/H_i]$$
is the desired module by
Theorem \ref{invariants of caln(Gamma) otimes_{zzGamma} zz}.1
and Additivity and Cofinality (see Theorem
\ref{properties of new definition}.4.).\par\noindent
1.b.)  By assertion 1.a) we can  choose
a sequence of countably generated projective $\cc\Gamma$-modules
$P_3$, $P_4$, $\ldots$ such that for $p \ge 3$
\begin{eqnarray}
\dim_{\caln(\Gamma)}(\caln(\Gamma) \otimes_{\zz\Gamma} P_p) & = & \beta_p.
\label{eqn 9.1}
\end{eqnarray}
Next we construct inductively a nested sequence
\mbox{$X_2 \subset X_3 \subset \ldots $}
of $\Gamma$-$CW$-complexes together
with $\Gamma$-retractions \mbox{$r_p: X_p \longrightarrow X_{p-1}$}
for \mbox{$p \ge 3$}
such that $X_2$ is the $2$-skeleton of a model for
$E\Gamma$, $X_p$ is obtained from $X_{p-1}$ by attaching
countably many free $\Gamma$-equivariant $p$-cells and $p+1$-cells and
\begin{eqnarray}
H_n(X_p,X_{p-1}) & = & \left\{ \begin{array}{ll}
P_p & n = p\\ 0 & n \not= p \end{array}\right. .
\label{eqn 9.2}
\end{eqnarray}
The Eilenberg-swindle yields a split-exact sequence
\textshortexactsequence{C_{p+1}}{c_{p+1}}{C_p}{}{P_p}
of $\cc\Gamma$-modules
such that $C_{p+1}$ and $C_p$ are countably generated free
$\zz\Gamma$-modules with
a basis.  Now one attaches for each element of the basis
of $C_p$ trivially a free $\Gamma$-equivariant $p$-cell to $X$.
Then one attaches for each element of the basis
of $C_{p+1}$ a free $\Gamma$-equivariant $p+1$-cell to $X$,
where the attaching maps
are choosen such that the cellular $\cc\Gamma$-chain complex
of $(Y,X)$ is just the $\zz\Gamma$-chain complex
which is concentrated in dimension $p+1$ and $p$ and given there
by \mbox{$C_{p+1} \xrightarrow{c_{p+1}} C_p$}.
Details of the construction of the $\Gamma$-$CW$-complexes
$X_p$ and $\Gamma$-retractions $r_p$
can be found in
\cite[Theorem 2.2, page 201]{Mislin (1980)}, \cite{Wall(1965a)}.
Now define
\mbox{$Y  ~ = ~ \colim_{p \to \infty} X_p$}.
One easily checks for $p\ge 3$
\begin{eqnarray*}
H_p^{\Gamma}(Y;\caln(\Gamma)) & = & \caln(\Gamma)\otimes_{\zz\Gamma} P_p;
\\
b_p^{(2)}(Y;\caln(\Gamma)) & = & \beta_p.
\end{eqnarray*}
2.) Let \mbox{$f: Y \longrightarrow X$} be a
$\Gamma$-$CW$-approximation of $X$ \cite[page 35]{Lueck (1989)},
i.e. a $\Gamma$-$CW$-complex $Y$ with a $\Gamma$-map $f$ such that $f^H$
is a weak homotopy equivalence and hence a weak homology equivalence
\cite[Theorem IV.7.15 on page 182]{Whitehead (1978)}
for $H \subset \Gamma$. Then
\mbox{$b_p^{(2)}(Y;\caln(\Gamma)) = b_p^{(2)}(X;\caln(\Gamma))$}
for all \mbox{$p \ge 0$} by Lemma
\ref{weak homotopy invariance}.2. Hence we can assume without loss of
generality that $X$ is a $\Gamma$-$CW$-complex.

Conjecture \ref{Atiyah Conjecture}
is equivalent to the statement that for any finitely presented
$\zz\Gamma$-module $M$
\begin{eqnarray}
d \cdot \dim_{\caln(\Gamma)}(\caln(\Gamma) \otimes_{\zz\Gamma} M)
& \in  & \zz.
\label{eqn 9.6775}
\end{eqnarray}
This follows essentially from \cite[Lemma 2.2]{Lueck (1997)}.
Now let $D_*$ be any $\zz\Gamma$-chain complex such that
$D_p$ is isomorphic to \mbox{$\oplus_{i=1}^r \zz[\Gamma/H_i]$}
for some non-negative integer $r$ and finite subgroups
$H_i$. Then the cokernel of each of the differentials $d_p$
is a finitely presented $\zz\Gamma$-module and $\ref{eqn 9.6775}$
yields for all $p \ge 0$
using Additivity (see Theorem \ref{properties of new definition}.4)
and Lemma \ref{invariants of caln(Gamma) otimes_{zzGamma} zz}.1
\begin{eqnarray}
d \cdot \dim_{\caln\Gamma)}
\left(\cok\left(\id_{\caln(\Gamma)} \otimes_{\zz\Gamma} d_p\right)\right)
& \in  & \zz; \nonumber
\\
d \cdot \dim_{\caln\Gamma)}
\left(\im\left(\id_{\caln(\Gamma)} \otimes_{\zz\Gamma} d_p\right)\right)
& \in  & \zz; \nonumber
\\
d \cdot \dim_{\caln\Gamma)}
\left(\ker\left(\id_{\caln(\Gamma)} \otimes_{\zz\Gamma} d_p\right)\right)
& \in  & \zz;\nonumber
\\
d \cdot \dim_{\caln\Gamma)}\left(
H_p(\caln(\Gamma) \otimes_{\zz\Gamma} D_*)\right)
& \in  & \zz.
\label{eqn 9.561}
\end{eqnarray}
Let $X^{\infty}$ be the $\Gamma$-$CW$-subcomplex of $X$ consisting
of points whose isotropy groups are infinite.
The sequence
\textshortexactsequence{C_*(X^{\infty})}{}{C_*(X)}{}{C_*(X,X^{\infty})}
of cellular $\zz\Gamma$-chain complexes is exact.
Since it is $\zz\Gamma$-split exact in each dimension, the sequence
obtained by tensoring with $\caln(\Gamma)$ is still exact.
The associated long exact homology sequence,
Additivity (see Theorem \ref{properties of new definition}.4)
and Lemma \ref{invariants of caln(Gamma) otimes_{zzGamma} zz}.1 imply
for $p \ge 0$
\begin{eqnarray}
\dim_{\caln(\Gamma)}(\caln(\Gamma) \otimes_{\zz\Gamma} C_p(X^{\infty}))
& = & 0;
\\
b_p^{(2)}(X;\caln(\Gamma)) & = &
\dim_{\caln(\Gamma)}
\left(H_p(\caln(\Gamma) \otimes_{\zz\Gamma} C_*(X,X^{\infty}))\right).
\label{eqn 9.766}
\end{eqnarray}
Notice that $C_p(X,X^{\infty})$ is a sum of $\zz\Gamma$-modules
of the shape \mbox{$\zz[\Gamma/H]$} for finite groups
\mbox{$H \subset \Gamma$}. Hence $C_*(X,X^{\infty})$
is a colimit over a directed set $I$ of subcomplexes
$D_*[i]$ (directed by inclusion)
such that each $D_p[i]$ is isomorphic to
\mbox{$\oplus_{i=1}^r \zz[\Gamma/H_i]$} for some non-negative
integer $r$ and finite subgroups $H_i$.
Since homology commutes with colimits
we conclude from \ref{eqn 9.766} and Theorem
\ref{dimension of colimits for arbitrary index sets}.2
\begin{eqnarray*}
b_p^{(2)}(X;\caln(\Gamma)) & = &
\sup \{\inf \{ \dim_{\caln(\Gamma)}(
\im(H_p(\caln(\Gamma) \otimes_{\zz\Gamma} D_*[i]) \longrightarrow
H_p(\caln(\Gamma) \otimes_{\zz\Gamma} D_*[j])))
\\
& &  \hspace*{20mm}
~ |~ j \in I, i \le j\} ~ |~ i \in I\}.
\end{eqnarray*}
Since the set \mbox{$\{r \in \rr \mid d \cdot r \in \zz\}$}
is discrete in $\rr$, it suffices to show for
each inclusion \mbox{$\iota: D_*[i] \longrightarrow D_*[j]$}
and all \mbox{$p \ge 0$}
\begin{eqnarray}
d \cdot \dim_{\caln(\Gamma)}\left(\im(
\iota_*: H_p(\caln(\Gamma)\otimes_{\zz\Gamma} D_*[i]) \longrightarrow
H_p(\caln(\Gamma)\otimes_{\zz\Gamma} D_*[j]))\right)
& \in \zz.
\label{eqn 9.6}
\end{eqnarray}
Let $F_*$ be any acyclic $\caln(\Gamma)$-chain complex with
\mbox{$F_p = 0$} for \mbox{$p < 0$} such that
\mbox{$d \cdot \dim(F_p) \in \zz$}
holds for all \mbox{$p \ge 0$}. Then we get
\mbox{$d \cdot \dim(\im(f_p: F_p \longrightarrow F_{p-1})) \in \zz$}
for all \mbox{$p \ge 0$} since we have the short exact sequences
\textshortexactsequence{\im(f_{p+1})}{}{F_p}{}{\im(f_p)}
and \mbox{$\im(f_1) = F_0$}.
Hence we obtain \ref{eqn 9.6} from \ref{eqn 9.561}
and the conclusion above for the case
where $F$ is the long exact homology sequence of the pair
\mbox{$(\cyl(\iota),D_*[i]$)} since there is a
$\zz\Gamma$-chain homotopy equivalence from the mapping cylinder
$\cyl(\iota)$ to $D_*[j]$ whose composition with the inclusion of
$D_*[i]$ in $\cyl(\iota)$ is $\iota$.
\par\noindent
3.) is proven in \cite[section 4]{Cheeger-Gromov (1986)}.
This finishes the proof of Theorem
\ref{possible values of the L^2-Betti numbers}.\qed

\begin{remark}
\label{groups and L^2-Betti numbers and homology} \em
The group \mbox{$\Gamma = \prod_{i = 1}^{\infty} \zz \ast\zz$}
satisfies \mbox{$H_p^{\Gamma}(E\Gamma;\caln(\Gamma)) = 0$}
for all $p \ge 0$. This is interesting in connection
with the zero-in-the-spectrum conjecture
(\cite{Lott (1997)}, \cite[section 11]{Lueck (1997)}.\em \qed
\end{remark}


\typeout{--------------------- references -----------------------}

\begin{center}
 Current address\\
Wolfgang L\"uck \\
Fachbereich Mathematik und Informatik\\
Westf\"alische Wilhelms-Universit\"at M\"unster\\
Einsteinstr. 62\\
48149 M\"unster\\
Bundesrepublik Deutschland\\
email:  lueck@math.uni-muenster.de\\
FAX: 0251 8338370\\
internet: http://wwwmath.uni-muenster.de/math/u/lueck/
\end{center}

\begin{center} Version of \today \end{center}


\end{document}